\documentclass[12pt]{article}

\newif\ifpdf
\ifx\pdfoutput\undefined
    \pdffalse       
\else
    \pdfoutput=1    
    \pdftrue
\fi

\usepackage{amsmath,amssymb,amsthm}

\ifpdf
    \usepackage[pdftex]{graphicx}
    \usepackage[backref, colorlinks=true, pdfstartview=FitH, linkcolor=blue,
        citecolor=blue, urlcolor=blue]{hyperref}
\else
    \usepackage{graphicx}
    \usepackage{hyperref}
\fi

\title{Continuous Ramsey Theory and Sidon Sets}
\author{Greg Martin (gerg@math.ubc.ca)\\ and \\ Kevin O'Bryant (kobryant@math.ucsd.edu)}
\date{September 28, 2002}

\newcommand{\Z}{{\mathbb Z}}

\newcommand{\R}{{\mathbb R}}

\newcommand{\K}{{\mathbb K}}
\newcommand{\T}{{\mathbb T}}

\newcommand{\A}{A}
\newcommand{\Prob}{\mathop{\rm Pr}}

\newcommand{\vecd}{{\bf d}}
\newcommand{\vecz}{{\bf z}}
\newcommand{\e}{\varepsilon}
\newcommand{\De}{\Delta(\varepsilon)}

\newcommand{\ff}{f\ast f}
\newcommand{\Bg}{\mbox{$B^{\ast}[g]$} }
\newcommand{\Bgm}{\mbox{$B^{\ast}[g]\pmod{n}$} }
\newcommand{\ffi}{\|\ff\|_\infty}

\newcommand{\Dconstant}{0.591389}
\newcommand{\twotimesDconstant}{1.182778} 
\newcommand{\Rconstant}{1.30036} 
\newcommand{\Deltaonehalfconstant}{0.1496} 
\newcommand{\TrivialLowerBoundconstant}{0.61522} 
\newcommand{\fstarftwonormconstant}{1.14915} 

\newcommand{\bigO}[1]{O\left(#1\right)}

\newcommand{\floor}[1]{\left\lfloor #1 \right\rfloor}
\newcommand{\ceiling}[1]{\big\lceil #1 \big\rceil}
\newcommand{\lnorm}[3]{{}_{#1}\| #3 \|_{#2}}
\newcommand{\field}[1]{\text{GF}\!\left( #1 \right)}
\newcommand{\sdr}[1]{ {#1}^{\text{sdr}}}

 \DeclareMathOperator{\supp}{supp}

\newtheorem{thm}{Theorem}[section]
\newtheorem{lem}[thm]{Lemma}

\newtheorem{cor}[thm]{Corollary}

\newtheorem{cnj}[thm]{Conjecture}
\newtheorem{prop}[thm]{Proposition}

\newenvironment{proofof}[1]{\medskip\noindent{\em Proof of #1.}}{\qed\medskip}

\begin{document}
    \maketitle

\begin{abstract}
A symmetric subset of the reals is one that remains invariant under some reflection
$x\mapsto c-x$. Given $0<\e\le1$, there exists a real number $\De$ with the following
property: if $0\le\delta<\De$, then every subset of $[0,1]$ with measure $\e$ contains a
symmetric subset with measure $\delta$, while if $\delta>\De$, then there exists a
subset of $[0,1]$ with measure $\e$ that does not contain a symmetric subset with
measure $\delta$. In this paper we establish upper and lower bounds for $\De$ of the
same order of magnitude: for example, we prove that $\De=2\e-1$ for
$\frac{11}{16}\le\e\le1$ and that $0.59\e^2<\De<0.8\e^2$ for $0<\e\le\frac{11}{16}$.

This continuous problem is intimately connected with a corresponding discrete problem. A
set $S$ of integers is called a $\Bg$ set if for any given $m$ there are at most $g$
ordered pairs $(s_1,s_2)\in S \times S$ with $s_1+s_2=m$; in the case $g=2$, these are
better known as Sidon sets. We also establish upper and lower bounds of the same order
of magnitude for the maximal possible size of a $\Bg$ set contained in $\{1,\dots,n\}$,
which we denote by $R(g,n)$. For example, we prove that $R(g,n)<1.31\sqrt{gn}$ for all
$n\ge g\ge2$, while $R(g,n) > 0.79\sqrt{gn}$ for sufficiently large integers $g$ and
$n$.

These two problems are so interconnected that both continuous and discrete tools can be
applied to each problem with surprising effectiveness. The harmonic analysis methods and
inequalities among various $L^p$ norms we use to derive lower bounds for $\De$ also
provide uniform upper bounds for $R(g,n)$, while the techniques from combinatorial and
probabilistic number theory that we employ to obtain constructions of large $\Bg$ sets
yield strong upper bounds for $\De$.
\end{abstract}

\thispagestyle{empty} \pagebreak \tableofcontents \pagebreak

\section{Introduction}

A set $C\subseteq\R$ is {\em symmetric} if there exists a number $c$ (the center of $C$)
such that $c+x\in C$ if and only if $c-x\in C$. Given a set $A\subseteq[0,1)$ of
positive measure, is there necessarily a symmetric subset $C\subseteq A$ of positive
measure? The answer turns out to be ``yes'', and the main topic of this paper is to
determine how large, in terms of the Lebesgue measure of $A$, one may take the symmetric
set $C$. In other words, for each $\e>0$ we are interested in
\begin{equation}
    \De := \sup \left\{\delta \colon\quad
    \begin{matrix}
        \text{every measurable subset of $[0,1)$ with measure $\e$} \\
        \text{contains a symmetric subset with measure $\delta$}
    \end{matrix}
        \right\}.
\label{De.definition}
\end{equation}
It is not immediately obvious that $\De>0$.

We have dubbed this sort of question ``continuous Ramsey theory'', and we direct the
reader to Section 2 for problems with a similar flavor; some of these have appeared in
the literature and some are given here for the first time.

We determine a lower bound for $\De$ using tools and ideas from harmonic analysis,
nonstandard analysis, and the theory of wavelets. We also construct sets without large
symmetric subsets using results from the theory of finite fields and probabilistic
number theory. These two lines of attack complement each other, and our bounds yield new
results in additive number theory as well.

Consider first the analogous discrete problem: given a set $A\subseteq\{1,2,\dots,n\}$,
how large is the largest symmetric subset of $A$? There are $\sim \tfrac12 |A|^2$ pairs
of distinct elements of $A$ (where $|A|$ denotes the cardinality of $A$), and each pair
$(a,b)$ has a center $\tfrac{a+b}{2}$ which is among the $\sim 2n$ values
$\{\tfrac32,2,\tfrac52,3,\dots,\tfrac{2n-1}{2}\}$. By the pigeonhole principle, there is
some $c$ that is the center of at least $\sim \frac12 |A|^2/(2n)$
pairs of elements of $A$. The union of those pairs is a symmetric set, i.e., there is a
symmetric set $C\subseteq A$ with
    $$
    \frac{|C|}{n} \gtrsim \frac2n \frac{\frac12 |A|^2}{2n}
    =\frac12 \left(\frac{|A|}{n}\right)^2.
    $$
In other words, the density of $C$ is roughly at least one-half
the square of the density of~$A$.

A {\em Sidon set} is a set $A$ of integers with the property that distinct pairs of
elements have distinct sums: if $a,b,c,d \in A$ and $a+b=c+d$, then $\{a,b\}=\{c,d\}$.
This is equivalent to asserting that $A$ has no symmetric subsets with more than 2
elements. It is known \cite{1938.Singer} that there is a Sidon set $A$ contained in
$\{1,2,\dots,n\}$ with $|A| \sim \sqrt{n}$. Thus, if $C\subseteq A$ is a symmetric set,
then
    $$
    \frac{|C|}{n} \leq \frac 2n \sim 2 \left(\frac{|A|}{n}\right)^2.
    $$
In other words, the density of $C$ is roughly at most twice the square of
the density of~$A$.

For the discrete version of the problem, at least, we see that one can guarantee a
``quadratically large'' symmetric subset, and that one cannot do
better in general.

\subsection{Continuous Results}\label{Intro.Continuous.sect}

Let $\lambda$ denote Lebesgue measure on $\R$. We are led by analogy with the discrete
problem to guess that every subset $A\subseteq[0,1)$ has a symmetric subset with measure
$\frac12 \lambda(A)^2$, and that there are subsets $A\subseteq [0,1)$ that do not have
symmetric subsets with measure larger than $2\lambda(A)^2$. That is, we are led to guess
that $\tfrac12 \e^2 \leq \De \leq 2 \e^2$.

We find the following equivalent definition of $\De$ easier to work with than
the definition given in Eq.~\eqref{De.definition}: if we define
    \begin{equation}
    D(A) := \sup\{ \lambda(C)\colon\quad C\subseteq A,\, \text{$C$ is symmetric}\},
    \label{Ddef}
    \end{equation}
then
    \begin{equation}
    \De := \inf\{ D(A)\colon\quad A\subseteq [0,1),\, \lambda(A)=\e\}.
    \label{Deltadef}
    \end{equation}

Write $A(x)$ for the indicator function of a set $A$, and define the sumset
$A+A:=\{a_1+a_2\colon a_i\in A\}$. We notice first that the maximal symmetric subset of
$A$ with center $c$ is $A \cap (2c-A)$, where $2c-A = \{2c-a\colon a\in A\}$; this
intersection has measure $\int A(x)A(2c-x)\,dx$, which can be written as the value of
the convolution $A\ast A(2c)$. This means that $D(A)$ simply equals the supremum norm
$\| A \ast A\|_\infty$. Since $\supp(A\ast A)$, the support of the function $A\ast A$,
equals $A+A$ (up to a set of measure zero) and is thus contained in $[0,2)$, we see that
    $$
    \|A\ast A\|_\infty \geq \frac{\|A\ast A\|_1}{\lambda(\supp(A\ast A))}
        = \frac{\lambda(A)^2}{\lambda(A+A)}\geq \frac12 \lambda(A)^2,
    $$
and hence $\De \geq \frac12 \e^2$ (as we had guessed from the analogy with the discrete
case). This lower bound, which we shall call the trivial lower bound on $\De$, is not so
far from the best we can derive! In fact, the bulk of this paper is devoted to improving
the constant in this lower bound from $\frac12$ to $\Dconstant$. Moreover, we are able
to establish a complementary upper bound for $\De$ in a manner that we shall discuss in
the next section.

In addition to a heavy reliance on Fourier analysis, we make use of wavelets, albeit
only in a rather naive manner. Although much of our argument (and indeed the initial
problem of bounding $\De$) was initially motivated by nonstandard analysis, we do not
make direct use of it here.

\begin{figure}
    \begin{center}
    \begin{picture}(384,240)
        \put(12,230){$\De/\e^2$}
        \put(370,12){$\e$}
    \ifpdf
        \put(0,-126){\includegraphics{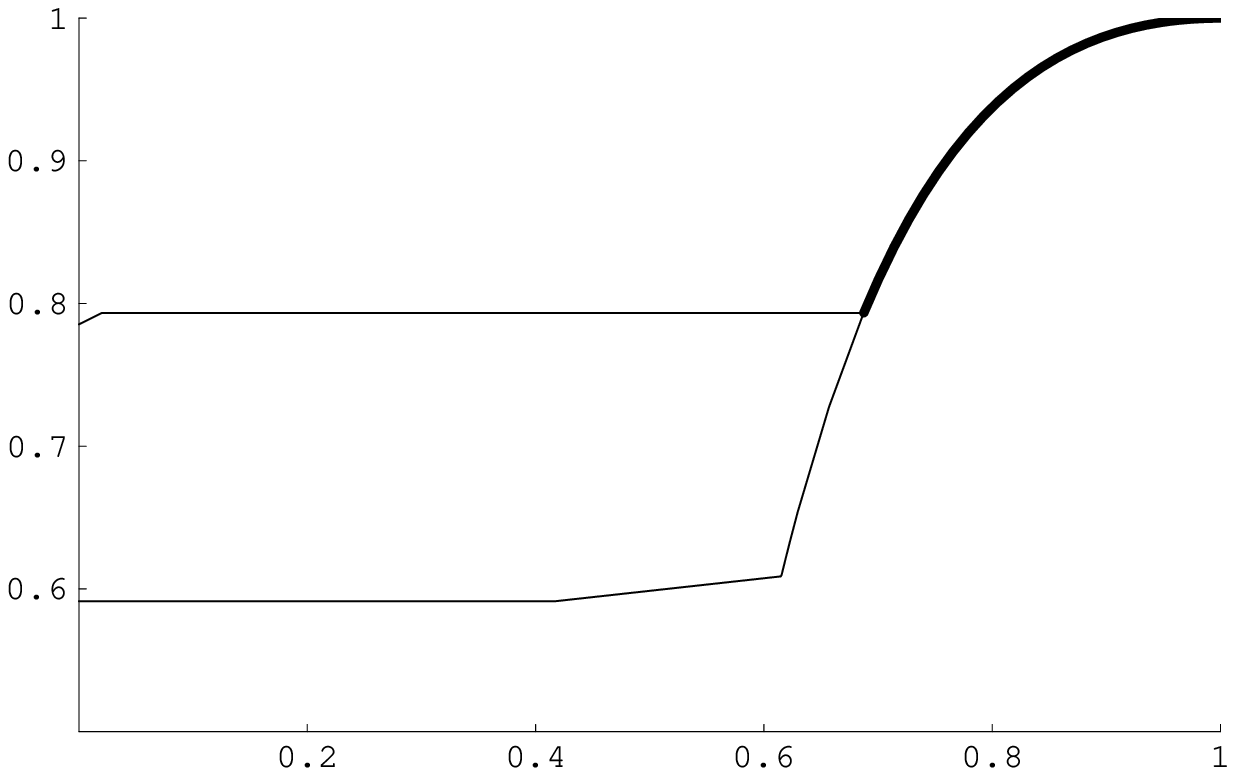}}
    \else
        \put(0,0){\includegraphics{BothBoundsPicture}}
    \fi
    \end{picture}
    \end{center}
    \caption{Upper and Lower Bounds for $\De/\e^2$}
    \label{overallpic}
\end{figure}

Figure~\ref{overallpic} shows the precise upper and lower bounds we obtain for
$\De/\e^2$ as functions of $\e$, and Theorem~\ref{Delta.Summary.thm} gives the
highlights:
\begin{thm}
\label{Delta.Summary.thm} We have:
\renewcommand{\theenumi}{\roman{enumi}}
    \begin{enumerate}
        \item $\Delta(\e)=2\e-1$ for $\frac{11}{16}\le\e\le1$, and $\De \ge2\e-1$ for
                 $\frac12\le\e\le\frac{11}{16}$;
        \item $\De \ge \Dconstant \e^2$ for all $0<\e\le1$;
        \item $\De \ge 0.5546\e^2 + 0.088079\e^3$ for all $0<\e \le1$;
        \item $\De \le \frac{96}{121}\e^2 < 0.7934 \e^2$ for
$0<\e\le\frac{11}{16}$;
        \item $\De \le \frac{\pi\e^2}{(1+\sqrt{1-\e})^2} = \frac\pi4\e^2 +
O(\e^3)$ for all $0<\e\le1$.
    \end{enumerate}
\end{thm}

Note that $\frac\pi4<0.7854$. The upper bound in part (v) of the theorem is superior to
the one in part (iv) in the range $0<\e<\frac{11}{96}(8\sqrt{6\pi}-11\pi) \doteq
0.0201$. The five parts of the theorem are proved separately in
Proposition~\ref{Line.2e-1.prop}, Proposition~\ref{the.full.De.bound.prop},
Proposition~\ref{Delta.one.half.thm}, Corollary~\ref{part.iii.cor}, and
Proposition~\ref{part.iv.prop}, respectively.

Thus we have established that $\e^2$ is the correct order of magnitude for the function
$\Delta$, and we have improved upon both the constants $\frac12$ and 2 which appeared in
our heuristic. In the following subsection, we discuss how we might improve the
corresponding constants for the discrete problem as well.

\subsection{Connecting the Continuous to the Discrete}

It turns out that the upper bounds given in Theorem~\ref{Delta.Summary.thm}(iv)--(v) are
derived from number-theoretic considerations. A set $S$ of integers is called a $B_2[g]$
set if for any given $m$ there are at most $2g$ ordered pairs $(s_1,s_2)\in S \times S$
with $s_1+s_2=m$ (equivalently, if the coefficients of $\left(\sum_{n\in S}
z^n\right)^2$ are bounded by $2g$). Sidon used $B_2[1]$ sets, which are the Sidon sets
mentioned earlier, as a tool in his study of Fourier series. It is perhaps fitting that
we now use Fourier analysis as a tool in our study of $B_2[g]$ sets.

Many constructions of $B_2[g]$
sets have appeared in recent years. Our constructions of sets of reals without large
symmetric subsets are based on known constructions of $B_2[g]$ sets. This geometric
aspect of $B_2[g]$ sets has led us to generalize and optimize these constructions
further and to give improved upper bounds on the density of $B_2[g]$ sets contained in
$\{1,2,\dots,n\}$.

Given a $B_2[g]$ set $S\subseteq\{1,2,\dots,n\}$ and any integer $m$, the union of all
pairs $(s_1,s_2)\in S \times S$ with $s_1+s_2=m$ is the largest symmetric subset of $S$
with center $\frac m2$, and all symmetric subsets of $S$ arise in this way. From the
definition of a $B_2[g]$ set, we see that $S$ contains no symmetric subset with more
than $2g$ integers. In light of this, it is not surprising (though it certainly requires
proof---see Proposition \ref{Delta.from.R.Connection.prop}) that the largest symmetric
subset of the set of real numbers
    \begin{equation}
    A(S) := \bigcup_{s\in S} \Big[ \frac{s-1}{n},\frac{s}{n}
    \Big)\subseteq[0,1)
    \label{ASdefinition}
    \end{equation}
has measure at most $\frac{2g}{n}$. In other words,
$\Delta\big(\frac{|S|}{n}\big) \leq D(A(S)) \le \frac{2g}{n}$.

For technical reasons, it is more convenient for us to speak of $B^\ast[g]$ sets
rather than $B_2[g]$ sets. A set $S$ of integers is called a $B^\ast[g]$ set if
for any given $m$ there are at most $g$ ordered pairs $(s_1,s_2)\in S \times S$
with $s_1+s_2=m$. Note that the definitions of $B_2[g]$ sets and $B^\ast[2g]$
sets coincide, but we shall also consider $B^\ast[g]$ sets with $g$ odd. We
introduce the function
    \begin{equation}
    R(g,n) := \max\{ |S|\colon\quad
            S\subseteq \{1,2,\dots,n\},\, \text{$S$ is a $B^\ast[g]$ set}\}.
    \label{Rgn.definition}
    \end{equation}
Since there are $|S|^2$ sums of pairs from $S$, while there are fewer than
$2n$ possible sums each of which can be realized at most $g$ times, we
immediately deduce the upper bound $R(g,n)\le\sqrt{2gn}$, which we shall call
the trivial upper bound for $R(g,n)$.

The construction of $A(S)$ indicated above thus gives the bound
$\Delta\big(\frac{R(g,n)}{n}\big) \leq \frac{g}{n},$ and we can use this inequality to
give an {\em upper} bound on $R(g,n)$. The trivial lower bound on $\Delta(\e)$ gives
$\frac12 \big(\frac{R(g,n)}{n}\big)^2 \leq \Delta\big(\frac{R(g,n)}{n}\big) \leq \frac
gn$, whence $R(g,n)\leq \sqrt{2gn}$. Thus, the trivial lower bound on $\De$ implies the
trivial upper bound on $R(g,n)$, and any improvement we can make on the trivial lower
bound for $\De$ will immediately provide stronger upper bounds for $R(g,n)$. In the same
way, we shall use lower bounds on $R(g,n)$ to derive upper bounds on $\De$. All this is
made rigorous in Sections~\ref{upper.bounds.Rgn.from.De.section}
and~\ref{upper.bounds.for.De.section}.

In fact much more than $\Delta\big(\frac{R(g,n)}{n}\big) \leq \frac{g}{n}$
is true. Proposition~\ref{Delta.as.infimum.prop} below states that
    $$
    \De = \inf\{ {\textstyle\frac gn} \colon n\ge g\ge1,\, R(g,n)\ge n\e \}.
    $$
Thus a sufficient understanding of the dependence of $R(g,n)$ on $g$ and $n$
would completely solve our continuous problem. Unfortunately, this understanding
is still somewhat lacking.

\subsection{Discrete Results}

The true size of Sidon sets is essentially known. Erd\H os and
Tur{\'a}n~\cite{1941.Erdos.Turan} exploited the fact that the pairwise differences $s_1-s_2$
from a Sidon set are distinct to establish the improved upper bound $R(2,n) \lesssim
\sqrt{n}$. (By the notation $f(n) \lesssim g(n)$, we mean that $\limsup_{n\to\infty}
\frac{f(n)}{g(n)} \le 1$.) Ruzsa~\cite{1993.Ruzsa} has observed that the Erd\H os/Tur{\'a}n
argument can be modified to give also $R(3,n)\lesssim \sqrt{n}$. A construction of
Singer~\cite{1938.Singer} (see Proposition~\ref{Modular.Constructions.prop}(iii) below)
yields the complementary lower bounds $R(3,n) \geq R(2,n) \gtrsim \sqrt n$.

Present knowledge of $R(g,n)$ for $g>3$ is much less impressive. Constructions of large
$B^\ast[g]$ sets have so far yielded only moderate lower bounds on $R(g,n)$. The trivial
upper bound on $R(g,n)$ has been improved for general $g$, but only quite recently. In
this paper, we present the strongest upper bound to date on $R(g,n)$ for all $g\geq 21$
as well as for all odd $g\ge5$, and we also improve or match the best known lower bounds
on $R(g,n)$.

The seminal paper of Cilleruelo, Ruzsa, and Trujillo~\cite{Cilleruelo.Ruzsa.Trujillo}
gave the first upper bound for $R(g,n)$ that was nontrivial for infinitely many $g$,
namely $R(2g,n) \lesssim 1.3181 \sqrt{2gn}$, and Green~\cite{2001.Green} improved this
to $R(2g,n)\lesssim 1.3038 \sqrt{2gn}$. Green also proved
\begin{equation}
R(2g,n) \lesssim \sqrt{\tfrac74(2g-1)n},
\label{Greens.bound}
\end{equation}
which is stronger than our results for even integers $g\le20$. (In
both these papers, the results were stated in terms of the function $F_2(g,n)$, a
notational difference only as $F_2(g,n) = R(2g,n)$ for all $g$ and $n$.)

It seems likely that $\lim_{n\to\infty} \frac{R(g,n)}{\sqrt{gn}}$ exists for
every $g$, but this has been proved only for $g=2$ and $g=3$. We define
    \begin{align*}
    \underline{\rho}(g) &:= \liminf_{n\to\infty} \frac{R(g,n)}{\sqrt{gn}}, \\
    \overline{\rho}(g)  &:= \limsup_{n\to\infty} \frac{R(g,n)}{\sqrt{gn}}, \\
    \rho(g)             &:= \lim_{n\to\infty} \frac{R(g,n)}{\sqrt{gn}}.
    \end{align*}
Thus $\rho(2)=\frac{1}{\sqrt{2}}$, $\rho(3)=\tfrac{1}{\sqrt{3}}$, and $\rho(g)$ is not
known to be well-defined for $g\geq 4$. Green's result~\cite{2001.Green} is equivalent
to $\overline{\rho}(g) \le \sqrt{\frac74}\sqrt{1-\frac1g}$.

The following theorem gives the bounds on $R(g,n)$ that we prove in this
work. Note that for even $g\leq20$, Green's bound \eqref{Greens.bound} is superior.

\begin{thm}\label{R.Upper.Bound.thm}
$R(g,n) \leq \Rconstant \sqrt{gn}$ for all $g$ and $R(g,n)\leq
1.31925 \sqrt{(g-1)n}+\frac13$ if $g$ is odd. Further, we have the following upper bounds for $\overline{\rho}(g)^2$:
$$\overline{\rho}(g)^2 \leq
\left\{%
\begin{array}{ll}
    1.74043-\frac{1.00483}g,
    & \hbox{$g\le8$ and even;} \\
    1.58337-\frac{0.026335}g + \sqrt{0.011572-\frac{0.083397}g+\frac{0.00069356}{g^2}},
    & \hbox{$g\ge 10$ and even;}\\
    1.74043-\frac{4.75492}g,
    & \hbox{$g\le 23$ and odd;} \\
    1.58337-\frac{0.071949}g + \sqrt{0.011572-\frac{0.22784}g+\frac{0.0051768}{g^2}},
    & \hbox{$g\ge 25$ and odd.}
\end{array}%
\right.
$$
\end{thm}

We comment that our result is an improvement in two aspects: we improve the numerical
constants in the bound on $\overline{\rho}(g)$ and we replace ``$\lesssim$'' with
``$\leq$'' in the bounds on $R(g,n)$. These minor improvements aside, we believe that
our geometrical approach is easier to follow and promises future improvements.
Accordingly, throughout this paper we discuss the quality of the results, what could
possibly be improved and what is best possible. In addition, to the best of our
knowledge, Theorem~\ref{R.Upper.Bound.thm} gives the first general upper bound for
$R(2g-1,n)$ that improves upon the trivial inequality $R(2g-1,n) \le R(2g,n)$. Theorem
\ref{R.Upper.Bound.thm} is proved in Corollary \ref{Delta.R.Basic.cor} and Proposition
\ref{R.from.Delta.Connection.cor}.

We turn now to lower bounds for $R(g,n)$. Constructions of ``$B^\ast[2] \pmod n$'' sets
have been given by Singer \cite{1938.Singer}, Bose~\cite{1942.Bose}, and
Ruzsa~\cite{1993.Ruzsa}. These constructions were extended to $B^\ast[g]$ sets
in~\cite{Cilleruelo.Ruzsa.Trujillo}. In Propositions~\ref{Modular.Constructions.prop}
and \ref{Modular.Construction.prop} below, we generalize the first three constructions
and further optimize the extension given in \cite{Cilleruelo.Ruzsa.Trujillo}.

Theorem~\ref{R.Lower.Bounds.thm} presents our improved lower bounds for $R(g,n)$, stated
in terms of the function $\underline{\rho}(g)$.

\begin{thm}
\label{R.Lower.Bounds.thm} We have
$$
\begin{array}{r@{{}\ge{}}l@{{}>{}}l}
    \underline{\rho}(4)  & \tfrac{2}{\sqrt7}          & 0.755, \\ \vspace{1mm}
    \underline{\rho}(6)  & \tfrac{2\sqrt2}{\sqrt{15}} & 0.730, \\ \vspace{1mm}
    \underline{\rho}(8)  & \tfrac2{\sqrt7}            & 0.755, \\ \vspace{1mm}
    \underline{\rho}(10) & \tfrac7{3\sqrt{10}}        & 0.737, \\ \vspace{1mm}
    \underline{\rho}(12) & \tfrac{\sqrt3}{\sqrt5}     & 0.774,
\end{array}
\qquad
\begin{array}{r@{{}\ge{}}l@{{}>{}}l}
    \underline{\rho}(14) & \tfrac{11}{\sqrt{210}}     & 0.759, \\ \vspace{1mm}
    \underline{\rho}(16) & \tfrac{17}{4\sqrt{30}}     & 0.775, \\ \vspace{1mm}
    \underline{\rho}(18) & \tfrac{4}{3\sqrt{3}}       & 0.769, \\ \vspace{1mm}
    \underline{\rho}(20) & \tfrac{2\sqrt{5}}{\sqrt{33}}& 0.778,\\ \vspace{1mm}
    \underline{\rho}(22) & \tfrac{18}{5\sqrt{22}}     & 0.767,
\end{array}
$$
and for any $g\ge12$,
    $$
    \underline{\rho}(2g) \ge \frac{g+2\floor{g/3}+\floor{g/6}}
    {\sqrt{6g^2-2g\floor{g/3}+2g}}\,.
    $$
In particular, for any $\delta>0$ we have $R(g,n)>(\frac{11}{8\sqrt3}-\delta)\sqrt{gn}$
if both $g$ and $\frac ng$ are sufficiently large in terms of~$\delta$.
\end{thm}

\begin{figure}
    \begin{center}
    \begin{picture}(384,240)
    \put(376,18){$g$}
    \ifpdf
        \put(0,-126){\includegraphics{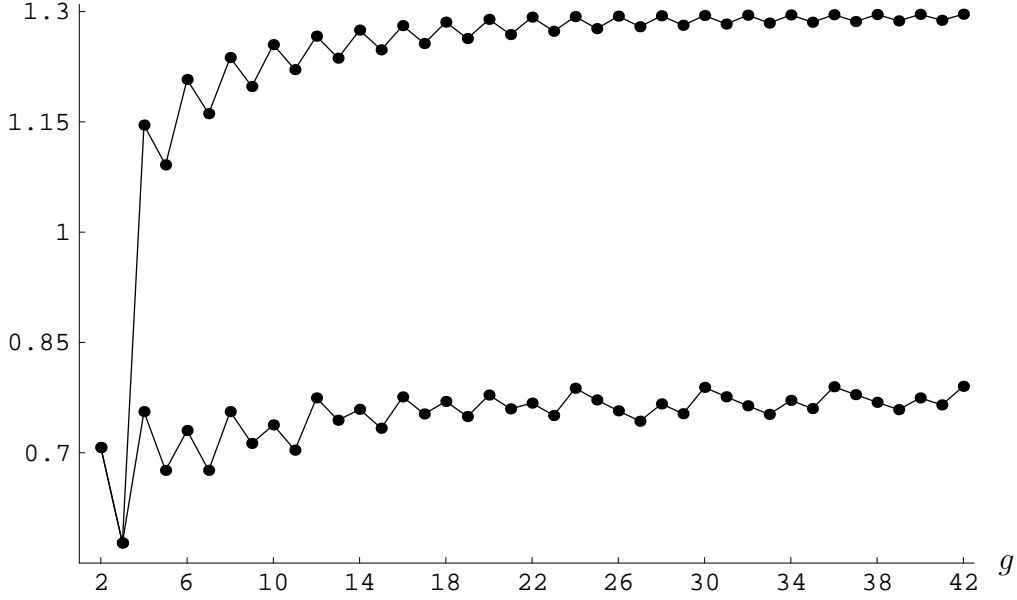}}
    \else
        \put(0,0){\includegraphics{BothBoundsOnR}}
    \fi
    \end{picture}
    \end{center}
    \caption{Lower bounds on $\underline{\rho}(g)$ and
                upper bounds on $\overline{\rho}(g)$}
    \label{BothBoundsOnR.pic}
\end{figure}

We note that $\frac{11}{8\sqrt3} > 0.7938$. The lower bound for $\underline{\rho}(4)$
reproduces a result of Habsieger and Plagne~\cite{Habsieger.Plagne}, while the lower
bounds for $\underline{\rho}(6)$ and $\underline{\rho}(10)$ reproduce results
of~\cite{Cilleruelo.Ruzsa.Trujillo}; for other even $g$, our lower bounds are new. These
lower bounds on $\underline\rho$, together with the strongest known upper bounds on
$\overline{\rho}$ including those derived from Theorem~\ref{R.Upper.Bound.thm}, are
plotted for $2\leq g \leq 42$ in Figure~\ref{BothBoundsOnR.pic}.

To the authors' knowledge, nothing more is known about lower bounds for \linebreak
$R(2g+1,n)$ for general $n$ and $g$ than the obvious inequality $R(2g,n)\leq R(2g+1,n)$.
In particular,
    $$
    \sqrt{\frac{2g}{2g+1}}\,\frac{R(2g,n)}{\sqrt{2gn}} \leq
    \frac{R(2g+1,n)}{\sqrt{(2g+1)n}},
    $$
which implies that
    $$
    \underline{\rho}(2g+1) \geq \sqrt{\frac{2g}{2g+1}} \,\underline{\rho}(2g).
    $$
It seems likely that this inequality is actually an equality (i.e., when $g$ is fixed,
the quotient $R(2g+1,n)/R(2g,n)$ should tend to 1 as $n$ tends to infinity), but this is
known only for $g=1$~\cite{1993.Ruzsa}.

Theorem~\ref{R.Lower.Bounds.thm} contains our best results for fixed $g$, but we can
obtain a better constant in the lower bound if we allow $g$ to grow slowly with $n$:

\begin{thm}
For any $\delta>0$, we have $R(g,n)>\big( \frac2{\sqrt\pi} -\delta \big)
\sqrt{gn}$ if both $\frac g{\log n}$ and $\frac ng$ are sufficiently large in
terms of~$\delta$.
\label{g.growing.lower.bound.thm}
\end{thm}

We note that $2/\sqrt{\pi} > 1.128$. Theorems~\ref{R.Lower.Bounds.thm}
and~\ref{g.growing.lower.bound.thm} are both important for ob\-tain\-ing the upper
bounds on $\De$ given in Theorem~\ref{Delta.Summary.thm}.
Theorem~\ref{R.Lower.Bounds.thm} is proved in Section~\ref{lower.bounds.bg.section}
through explicit constructions, while Theorem~\ref{g.growing.lower.bound.thm} is proved
in Section~\ref{probabilistic.bg.section} by a probabilistic argument.

As mentioned above, in Section 2 we describe several other problems in continuous Ramsey
theory, along with their discrete analogues. In Section 3 we list some easy-to-derive
properties of the function $\De$, while in Section 4 we describe how we establish our
strongest lower bounds for $\De$. Section 4 is the most involved part of this paper,
employing many techniques from harmonic analysis, most notably Fourier series. The
connection between the continuous and discrete problems measured by $\De$ and $R(g,n)$,
respectively, is given in Section 5, and our upper bounds for $R(g,n)$ are derived
therein. In that section, we also establish in Theorem \ref{ubiquity.thm} that large \Bg
sets must have many pairwise sums that repeat at least $\alpha g$ times for
suitable~$\alpha$, and in Theorem~\ref{UD.Hypothesis.thm} we show a stronger upper bound
for \Bg sets that are uniformly distributed in $\{1,\dots,n\}$ (which we conjecture is
typical for the largest possible \Bg sets). Section 6 is devoted to lower bounds for
$R(g,n)$, which rely upon improved constructions of \Bg sets and their ``mod $n$''
counterparts. In Section 7 we apply the results of Section 6 to bound $\De$ from above.
Finally, in Section 8 we collect together several open questions and problems relating
to our methods.

\pagebreak \section{Some Problems in Continuous Ramsey Theory}

A ``coloring Ramsey theorem'' has the form:
\begin{quote}
Given a sufficiently large set of mathematical objects colored with a finite number of
colors, there is a highly structured monochromatic subset.
\end{quote}
The prototypical example is Ramsey's Theorem itself: However one colors the edges of the
complete graph $K_n$ with $r$ colors, there is a monochromatic complete subgraph on $t$
vertices, provided that $n$ is sufficiently large in terms of $r$ and $t$. Another
example is van der Waerden's Theorem: However one colors the integers $\{1,2,\dots,n\}$
with $r$ colors, there is a monochromatic arithmetic progression with $t$ terms,
provided that $n$ is sufficiently large in terms of $r$ and $t$.

In many cases, the coloring aspect of a Ramsey-type theorem is a ruse and one may prove
a stronger statement with the form:
\begin{quote}
Given a sufficiently large set $R$ of mathematical objects, any large subset of $R$
contains a highly structured subset.
\end{quote}
Such a result is called a ``density Ramsey theorem.'' For example, van der Waerden's
Theorem is a special case of the density Ramsey theorem of Szemer{\'e}di: Every subset of
$\{1,2,\dots,n\}$ with cardinality at least $n\delta$ contains a $t$-term arithmetic
progression, provided that $n$ is sufficiently large in terms of $\delta$ and $t$.

Ramsey theory is the study of such theorems on different types of structures. By
``continuous Ramsey theory'' we refer to Ramsey-type problems on continuous measure
spaces. In particular, this thesis is concerned with a density-Ramsey problem on the
structure $[0,1)\subseteq\R$ with Lebesgue measure. The type of substructure we focus on
are symmetric subsets.

\begin{sloppypar}
There is frequently a connection between problems concerning $[0,1)$ with Lebesgue
measure and problems concerning $\{1,2,\dots,n\}$ with a uniform probability measure. We
attempt to make this connection more explicit in our descriptions of problems in the
remainder of this section.
\end{sloppypar}

\subsection{The Correlation Problem}
Erd\H{o}s asked for upper and lower bounds on $M_n(\tfrac12)$, where
    $$
    M_n(\alpha) := \inf_{\substack{|S| = \floor{\alpha n} \\ T= \{1,2,\dots,n\}
    \backslash S}} \sup_{c\in\Z} \frac{\left| (S+c)\cap T \right| }{n}.
    $$
In essence, Erd\H{o}s asked for a formalization of the observation: a large subset of
$\{1,2,\dots,n\}$ can always be translated so as to have a large intersection with its
complement.

S. {\'S}wierczkowski~\cite{1958.Swierczkowski} showed that $\lim_{n\to\infty}M_n(\alpha)
= M(\alpha)$, where
    $$
    M(\alpha)   :=\inf_{\substack{\lambda(S)=\alpha \\ T
                = [0,1) \backslash S}} \sup_{c\in\R} \lambda((S+c)\cap T).
    $$
Working in the continuous setting, {\'S}wierczkowski showed that
$$M(\alpha)>\frac{2-\sqrt{4-10\alpha(1-\alpha)}}{5},$$ which improved on the asymptotic bounds
then known for $M_n(\alpha)$.

It is known (see \cite[problem C17]{1994.Guy}) that $0.178 < M(\tfrac12) < 0.2$, but the
precise value remains unknown.

\subsection{The Convolution Problem}

Banakh, Verbitsky, and Vorobets~\cite{2000.Banakh.Verbitsky.Vorobets} considered the
coloring version of the problem considered in the present paper: given a measurable
finite coloring of $[0,1)$, how large a monochromatic symmetric set is guaranteed? They
claim that if 2 colors are used, then there is necessarily a monochromatic symmetric
subset with measure $\frac{1}{4+\sqrt{6}}\approx 0.155$, and if $r>2$ colors are used
then there is a monochromatic subset with measure $\frac{1/2}{r^2}$. Unfortunately, the
proofs of several crucial lemmas of this interesting paper appear only in earlier
Russian-language articles, and we have been unable to verify their proofs. If $r$ colors
are used, then there is a monochromatic set with measure at least $\frac 1r$. In
Section~\ref{full.bound.section}, we show that {\em any} subset of $[0,1)$ with measure
at least $\frac 1r$ contains a symmetric subset with measure $\frac{\Dconstant}{r^2}$,
and in Section~\ref{delta.half.section}, we show further that any subset with measure at
least $\frac 12$ has a symmetric subset with measure $\Deltaonehalfconstant$. This paper
is concerned with the density version of their coloring problem.

We call this the convolution problem because, if $f$ is the indicator function of $E
\subseteq \R$, then $f \ast f(c) := \int_{\R} f(x)f(c-x)\,dx$ is the measure of the
maximal symmetric subset of $E$ with center $c/2$. Therefore the
Banakh--Verbitsky--Vorobets problem is equivalent to finding the infimum of the possible
values of $$\max\{ \|f*f\|_\infty, \|(1-f)*(1-f)\|_\infty \}$$ as $f$ ranges over all
indicator functions of measurable subsets of $[0,1)$. We exploit this connection to
Fourier analysis to give lower bounds for the function $\De$ defined in
Eq.~(\ref{Deltadef}).

We note that the discrete versions of these problems are closely related to the
continuous versions. Specifically, define $MS(n,r)$ to be the maximum of integers $M$
such that for any coloring of $\{1,2,\dots,n\}$ with $r$ colors, there is a
monochromatic symmetric subset with cardinality $M$. Also, define $MS([0,1],r)$ to the
supremum of real $\mu$ such that for any measurable coloring of $[0,1]$ with $r$ colors,
there is a monochromatic subset with measure $\mu$. It is shown
in~\cite{2000.Banakh.Verbitsky.Vorobets} that $MS([0,1],r) = \lim_{n\to\infty}
MS(n,r)/n$. Proposition~\ref{Delta.as.infimum.prop} below gives the analogous connection
between $\De$ and the function $R(g,n)$ defined in Eq.~(\ref{Rgn.definition}).

\subsection{The Linear Equation Problem}

Following Chung, Erd\H{o}s, and Graham~\cite{Chung.Erdos.Graham}, for any matrix $A$, we
say that $S\subseteq\{1,2,\dots,n\}$ is $A$-hitting if for every vector
$\bar{x}=(x_1,x_2,\dots,x_s)$ with $x_i \in \{1,2,\dots,n\}$ and $\bar{x}A=\bar{0}$,
either there is some $x_i \in S$ or $x_1=x_2=\dots=x_s$ (this second condition is
included to avoid certain trivialities). Set $\delta_n$ to be the minimum density of an
$A$-hitting set, i.e.,
    $$
    \delta_n := \min \left\{\tfrac{|S|}{n}\colon\quad
                        S\subseteq\{1,2,\dots,n\},\,S\text{ is $A$-hitting}\right\}.
    $$
If, for example, $A=(1,-2,1)$, then $\delta_n$ is the minimum density of a subset of
$\{1,2,\dots,n\}$ that intersects every 3-term arithmetic progression in
$\{1,2,\dots,n\}$.

We also say that $S\subseteq[0,1]$ is $A$-hitting if for every vector
$\bar{x}=(x_1,x_2,\dots,x_s)$ with $x_i \in [0,1]$ and $\bar{x}A=\bar{0}$, either
$x_1=x_2=\dots=x_s$ or there is some $x_i \in S$. Set $\delta$ to be the infimum of the
measures of $A$-hitting sets, i.e.,
    $$
    \delta := \inf \left\{ \lambda(S) \colon\quad
                        S\subseteq[0,1],\, S\text{ is $A$-hitting} \right\}.
    $$

Given our experience with the above correlation and convolution problems, one might
guess that $\liminf_{n} \delta_n = \delta$. In~\cite{Chung.Erdos.Graham}, the values of
both $\delta$ and $\liminf_{n} \delta_n$ are derived for several matrices $A$, and
indeed equality seems to be the typical case. Surprisingly, however, there are
non-trivial examples where this is not the case. For
    $$
    A = \begin{pmatrix}
            2 & 3 \\
            -1 & 0 \\
            0 & -1
        \end{pmatrix},
    $$
for example, it is shown in~\cite{Chung.Erdos.Graham} that $\delta=\tfrac15$ whereas
$0.199 < \liminf_n \delta_n < 0.1997$.

\pagebreak\section{Easy Bounds for $\De$}

We now turn our attention to the investigation of the function $\De$
defined in Eq.~(\ref{Deltadef}). In this section we establish several simple
lemmas describing basic properties of $\Delta$.

\begin{lem}
\label{Delta.Trivial.Bounds.lem} $\De \ge 2\e-1$ for all $0\le\e\le1$.
\end{lem}

\begin{proof}
For $A\subseteq[0,1)$, the centrally symmetric set $A\cap(1-A)$ has measure equal to
    $$
    \lambda(A)+\lambda(1-A)-\lambda(A\cup(1-A))
        = 2\lambda(A)-\lambda(A\cup(1-A)) \geq 2\lambda(A)-1.
    $$
Therefore $D(A)\ge2\lambda(A)-1$ from the definition (\ref{Ddef}) of the function $D$.
Taking the infimum over all subsets $A$ of $[0,1)$ with measure $\e$, this becomes
$\De\ge2\e-1$ as claimed.
\end{proof}

While this bound may seem obvious, it is in many situations the state of the art. As we
show in Proposition~\ref{Line.2e-1.prop} below, $\De$ actually equals $2\e-1$ for
$\tfrac{11}{16}\le\e\le1$; and $\De \geq 2\e-1$ is the best lower bound of which we are
aware in the range $\TrivialLowerBoundconstant \leq \e <\tfrac{11}{16}\doteq 0.6875$.

One is tempted to try to sharpen the bound $\De \geq 2\e-1$ by considering the symmetric
subsets with center $1/3$, $1/2$, or $2/3$, for example, instead of merely $1/2$.
Unfortunately, it can be shown that given any $\e\ge\frac12$ and any finite set
$\{c_1,c_2,\dots,c_n\}$, one can construct a sequence $S_k$ of sets, each with measure
$\e$, that satisfies
    $$
    \lim_{k\to\infty}\left(\max_{1\leq i \leq n}
            \{\lambda\left(S_k\cap(2c_i-S_k)\right)\}\right)
                        = 2\e-1.
    $$
Thus, no improvement is possible with this sort of argument.

\begin{lem}[Trivial Lower Bound]
$\Delta(\e) \ge \frac12\e^2$ for all $0\le\e\le1$.
\label{Delta.Trivial.Lower.Bound.lem}
\end{lem}

\begin{proof}
We repeat the argument given briefly in Section~\ref{Intro.Continuous.sect}. Given a
subset $A$ of $[0,1)$ of measure $\e$, let $A(x)$ denote the indicator function of $A$,
so that the integral of $A(x)$ over any interval containing $[0,1)$ equals $\e$. If we
define $f(c) := \int_{-\infty}^\infty A(x) A(2c-x) \, dx$, then $f(c)$ is the measure of
the largest symmetric subset of $A$ with center $c$, and we seek to maximize $f(c)$. But
$f$ is clearly supported on $[0,1)$, and so
    \begin{multline*}
    D(A)
    = \max_c f(c) \ge \int_0^1 f(c)\,dc
    = \int_{-\infty}^\infty \int_0^1 A(x) A(2c-x) \, dc \, dx \\
    = \int_{-\infty}^\infty A(x) \bigg( \frac12 \int_{-x}^{2-x} A(w)\, dw \bigg) \, dx
    = \tfrac12 \e^2.
    \end{multline*}
Since $A$ was an arbitrary subset of $[0,1)$ of measure $\e$, we have shown that
\mbox{$\Delta(\e)\ge\frac12\e^2$}.
\end{proof}

\begin{lem}
\label{Delta.Line.To.Origin.lem} $\Delta(\e) \leq \frac{\Delta(x)}{x}\e$ for all $0\leq
\e \leq x\le1$. In particular, $\De \le \e$.
\end{lem}

\begin{proof}
If $tA := \{ta\colon a\in A\}$ is a scaled copy of a set $A$, then clearly
$D(tA)=tD(A)$. Applying this with any set $A\subseteq[0,1)$ of measure $x$ and with
$t=\frac\e{x}\le1$, we see that $\frac\e{x}A$ is a subset of $[0,1)$ with measure $\e$,
and so by the definition of $\Delta$ we have $\Delta(\e) \le D(\frac\e{x}A) =
\frac\e{x}D(A)$. Taking the infimum over all sets $A\subseteq[0,1)$ of measure $x$, we
conclude that $\Delta(\e) \le \frac\e{x}\Delta(x)$. The second assertion of the lemma is
obvious, and indeed it follows from the first assertion in light of the trivial value
$\Delta(1)=1$.
\end{proof}

It is obvious from the definition of $\Delta$ that $\De$ is an increasing function;
\linebreak Lemma~\ref{Delta.Line.To.Origin.lem} shows that $\frac{\De}\e$ is also an
increasing function. Later in this paper (see
Proposition~\ref{De.over.e2.increasing.prop}), we will show that in fact even
$\frac{\De}{\e^2}$ is an increasing function.

Let $$S\diamond T := (S\setminus T)\cup(T\setminus S)$$ denote the symmetric difference
of $S$ and $T$. (While this operation is more commonly denoted with a triangle rather
than with a diamond, we would rather avoid any potential confusion with the function
$\Delta$ featured prominently in this paper.)

\begin{lem}
If $S$ and $T$ are two sets of real numbers, then $|D(S)-D(T)| \le
2\lambda(S\diamond T)$. \label{Diamond.lem}
\end{lem}

\begin{proof}
Let $E$ be any symmetric subset of $S$, and let $c$ be the center of $E$, so
that $E=2c-E$. Define $F=E\cap T\cap(2c-T)$, which is a symmetric subset of
$T$ with center $2c$. We can write $\lambda(F)$ using the inclusion-exclusion
formula
    \begin{multline*}
    \lambda(F) = \lambda(E) + \lambda(T) + \lambda(2c-T) \\ - \lambda(E\cup T) -
    \lambda(E\cup(2c-T)) - \lambda(T\cup(2c-T)) + \lambda(E\cup T\cup(2c-T)).
    \end{multline*}
Rearranging terms, and noting that $T\cup(2c-T) \subseteq E\cup T\cup(2c-T)$, we see
that
    \begin{equation*}
    \lambda(E) - \lambda(F) \le -\lambda(T)-\lambda(2c-T)+\lambda(E\cup T)
            +\lambda(E\cup (2c-T))
    \end{equation*}
Because reflecting a set in the point $c$ does not change its measure, this is
the same as
    \begin{equation*}
    \begin{split}
    \lambda(E) - \lambda(F) &\le -\lambda(T)-\lambda(T)+\lambda(E\cup T)
            +\lambda(E\cup (2c-T)) \\
    &= -\lambda(T)-\lambda(T)+\lambda(E\cup T)+\lambda((2c-E)\cup T) \\
    &= 2\big( \lambda(E\cup T) - \lambda(T) \big) \\
    &\le 2\big( \lambda(S\cup T) - \lambda(T) \big)
    = 2\lambda(S\setminus T) \le 2\lambda( S\diamond T),
    \end{split}
    \end{equation*}
Therefore, since $F$ is a symmetric subset of $T$,
\begin{equation*}
\lambda(E) \le \lambda(F) + 2\lambda( S\diamond T) \le D(T) + 2\lambda(
S\diamond T).
\end{equation*}
Taking the supremum over all symmetric subsets $E$ of $S$, we conclude that $D(S) \le
D(T) + 2\lambda( S\diamond T)$. If we now exchange the roles of $S$ and $T$, we see that
the proof is complete.
\end{proof}

\begin{lem}
The function $\Delta$ satisfies the Lipschitz condition
$|\Delta(x)-\Delta(y)|\leq 2|x-y|$ for all $x$ and $y$ in $[0,1]$. In
particular, $\Delta$ is continuous. \label{Lipschitz.Condition.lem}
\end{lem}

\begin{proof}
Without loss of generality assume $y<x$. In light of the monotonicity $\Delta(y)\le
\Delta(x)$, it suffices to show that $\Delta(y) \ge \Delta(x) - 2(x-y)$. Let
$S\subseteq[0,1)$ have measure $y$. Choose any set $R\subseteq[0,1)\setminus S$ with
measure $x-y$, and set $T=S\cup R$. Then $S\diamond T=R$, and so by
Lemma~\ref{Diamond.lem}, $D(T)-D(S) \le 2\lambda(R) = 2(x-y)$. Therefore $D(S) \ge D(T)
- 2(x-y) \ge \Delta(x) - 2(x-y)$ by the definition of $\Delta$. Taking the infimum over
all sets $S\subseteq[0,1)$ of measure $y$ yields $\Delta(y) \ge \Delta(x) - 2(x-y)$ as
desired.
\end{proof}

%

\pagebreak \section{Lower Bounds for $\De$}

Section~\ref{introduction.section} below makes explicit the connection between $\De$ and
harmonic analysis. Section~\ref{basic.argument.section} gives a simple, but quite good,
lower bound on $\De$. In Section~\ref{main.bound.section}, we give a more general form
of the argument of Section~\ref{basic.argument.section}. Using an analytic inequality
established in Section~\ref{useful.inequalities.section}, we investigate in
Section~\ref{fourier.coefficients.section} the connection between $\ffi$ and the Fourier
coefficients of $f$. In Section~\ref{full.bound.section}, we combine the results of
Sections~\ref{basic.argument.section} and~\ref{fourier.coefficients.section} to show
$\De \geq \Dconstant \e^2$. The bound given in Section~\ref{basic.argument.section} and
improved in Section~\ref{full.bound.section} depends on a kernel function with certain
properties; in Section~\ref{kernel.problem.section} we discuss how we chose our kernel.
In Section~\ref{delta.half.section}, we use a different approach to derive a lower bound
on $\De$ which is superior for $\frac38<\e<\frac58$.

\subsection{Introduction and Notation}\label{introduction.section}

There are many ways to define the basic objects of Fourier analysis; we follow
\cite{1984.Folland}. Unless specifically noted otherwise, all integrals are over the
circle group $\T:=\R/\Z$; for example, $L^1$ denotes the class of functions $f$ for
which $\int_\T |f(x)|\,dx$ is finite. For each integer $j$, we define $\hat{f}(j):=\int
f(x) e^{-2\pi ijx}\,dx$, so that for any function $f\in L^1$, we have
$f(x)=\sum_{j=-\infty}^\infty \hat{f}(j)e^{2\pi ij x}$ almost everywhere. We define the
convolution $f*g(c)=\int f(x)g(c-x) \,dx$, and we note that $\widehat{f\ast
g}(j)=\hat{f}(j)\hat g(j)$ for every integer $j$; in particular,
$\widehat{\ff}(j)=\hat{f}(j)^2$.

We define the usual $L^p$ norms $$\textstyle\|f\|_p = \big( \int |f(x)|^p\,
dx\big)^{1/p}$$ and
    $$
   \|f\|_\infty    = \lim_{p\to\infty} \|f\|_p
                    = \sup\big\{y\colon\quad \lambda(\{x\colon |f(x)|>y\})>0\big\}.
    $$
With these definitions, H{\"o}lder's Inequality is valid: if $p$ and $q$ are conjugate
ex\-po\-nents---that is, $\frac1p+\frac1q=1$---then $\|fg\|_1 \leq \|f\|_p \|g\|_q$. We
also note that $\|f*g\|_1 = \|f\|_1\|g\|_1$; in particular,
$\|\ff\|_1=\|f\|_1^2=\hat{f}(0)^2$. We shall also employ the $\ell^p$ norms for
bi-infinite sequences: if $a=\{a_j\}_{j\in\Z}$, then $\|a\|_p = \big( \sum_{j\in\Z}
|a_j|^p \big)^{1/p}$ and $\|a\|_\infty = \lim_{p\to\infty} \|a\|_p = \sup_{j\in\Z}
|a_j|$. Although we use the same notation for the $L^p$ and $\ell^p$ norms, no confusion
should arise, as the object inside the norm symbol will either be a function on $\T$ or
its sequence of Fourier coefficients, respectively. With this notation, we recall
Parseval's identity $$\int f(x)g(x)\,dx = \sum \hat f(j)\hat g(-j)$$ (assuming the
integral and sum both converge); in particular, if $f=g$ is real-valued (so that $\hat
f(-j)$ is the conjugate of $\hat f(j)$ for all $j$), this becomes $\|f\|_2 =
\|\hat{f}\|_2$. The Hausdorff-Young inequality, $\|\hat{f}\|_q \leq \| f \|_p$ whenever
$p$ and $q$ are conjugate exponents with $1\le p\le 2\le q\le\infty$, can be thought of
as a generalization of this latter version of Parseval's identity. We also require the
definition
    \begin{equation}
    \lnorm{m}{p}{a} = \left(\sum_{|j|\geq m} |a(j)|^p \right)^{1/p}
    \label{lnorm.def}
    \end{equation}
for any sequence $a=\{a_j\}_{j\in\Z}$, so that $\lnorm0pa = \|a\|_p$, for example.

We recall that the Fourier coefficients of any function $f\in L^1$ satisfy the estimate
$\hat f(j) = \bigO{\frac1j}$; in particular, $\|\hat f\|_p$ is finite for all $p>1$.
Moreover, if $f\in L^1$ is continuous, then $\hat f(j)=\bigO{\frac 1{j^2}}$. We also
note that for any fixed sequence $a=\{a_j\}_{j\in\Z}$, the $\ell^p$-norm $\|a\|_p$ is a
decreasing function of $p$. To see this, suppose that $1\le p\le q<\infty$ and
$a\in\ell^p$. Then $|a_j| \le \|a\|_p$ for all $j\in\Z$, whence $|a_j|^{q-p} \le
\|a\|_p^{q-p}$ (since $q-p\ge0$) and so $|a_j|^q \le \|a\|_p^{q-p} |a_j|^p$. Summing
both sides over all $j\in\Z$ yields $\|a\|_q^q \le \|a\|_p^{q-p} \|a\|_p^p = \|a\|_p^q$,
and taking $q$th roots gives the desired inequality $\|a\|_q \le \|a\|_p$.

Finally, we define a ``pdf'', short for ``probability density function'', to be a
nonnegative function in $L^2$ whose $L^1$-norm (which is necessarily finite, since $\T$
is a finite measure space) equals~1. Also, we single out a special type of pdf called an
``nif'', short for ``normalized indicator function'', which is a pdf that only takes one
nonzero value, that value necessarily being the reciprocal of the measure of the support
of the function. (We exclude the possibility that an nif takes the value 0 almost
everywhere.) Specifically, we define for each $E\subseteq\T$ the nif
    $$
    f_E(x):=
        \begin{cases}
            \lambda(E)^{-1} & x\in E, \\
            0 & x\not\in E.
        \end{cases}
    $$
Note that if $f$ is a pdf, then $1=\hat{f}(0)=\hat{f}(0)^2=\|f\|_1^2=\|\ff\|_1$.

We are now ready to reformulate the function $\De$ in terms of this notation.

\begin{lem}
\label{Fourier.Connection.lem} We have
    $$ \textstyle
    \frac12 \e^2 \inf_g \|g*g\|_\infty \le \frac12 \e^2 \inf_f \ffi = \De,
    $$
the first infimum being taken over all pdfs $g$ that are supported on
$[-\tfrac14,\tfrac14]$, and the second infimum being taken over all nifs $f$ whose
support is a subset of $[-\frac14,\frac14]$ of measure $\frac\e2$.
\end{lem}

\begin{proof}
The inequality is trivial, since every nif is a pdf; it remains to prove the equality.

For each measurable $A\subseteq[0,1)$, define $E_A:=\{\frac12 (a-\frac12) \colon a\in
A\}\subseteq[-\frac14,\frac14]$. The sets $A$ and $E_A$ differ only by translation and
scaling, so that $\lambda(A)=2\lambda(E_A)$ and $D(A)=2D(E_A)$. Thus
    \begin{align*}
    \De     &:= \inf\{ D(A) \colon A\subseteq[0,1),\, \lambda(A)=\e\} \\
            &=  \e^2 \inf\left\{ \frac{D(A)}{\lambda(A)^2}
                            \colon A\subseteq[0,1),\, \lambda(A)=\e\right\} \\
            &=  \e^2 \inf\left\{ \frac{2D(E_A)}{(2\lambda(E_A))^2}
                            \colon A\subseteq[0,1),\, \lambda(A)=\e\right\} \\
            &=  \frac12 \e^2 \inf\left\{ \frac{D(E)}{\lambda(E)^2}
                            \colon E\subseteq[-\tfrac14,\tfrac14],\,
                            \lambda(E)=\tfrac{\e}{2}              \right\}.
    \end{align*}
For each $E\subseteq[-\tfrac14,\tfrac14]$ with $\lambda(E)=\frac{\e}{2}$, the function
$f_E(x)$ is an nif supported on a subset of $[-\tfrac14,\tfrac14]$ with measure
$\frac{\e}{2}$, and it is clear that every such nif arises from some set $E$. Thus, it
remains only to show that $\frac{D(E)}{\lambda(E)^2}=\| f_E \ast f_E \|_\infty$, i.e.,
that $D(E)=\lambda(E)^2 \| f_E \ast f_E \|_\infty$.

Fix $E\subseteq[-\tfrac14,\tfrac14]$, and let $E(x)$ be the indicator function of $E$.
Note that $f_E(x) = \lambda(E)^{-1}E(x)$. The maximal symmetric subset of $E$ with
center $c$ is $E\cap(2c-E)$, and this has measure $\int E(x) E(2c-x)\,dx$. Thus
    \begin{align*}
    D(E)    &:=\sup\{ \lambda(C)\colon\quad C\subseteq E,\, \text{$C$ is symmetric}\}\\
            &= \sup_c \left( \int E(x)E(2c-x)\,dx \right)\\
            &= \sup_c \left(
                    \int\lambda(E)f_E(x)\lambda(E)f_E(2c-x)\,dx \right)\\
            &= \lambda(E)^2 \sup_c \left( \int f_E(x)f_E(2c-x)\,dx \right)\\
            &= \lambda(E)^2 \sup_c f_E \ast f_E (2c) \\
            &= \lambda(E)^2 \left\| f_E \ast f_E \right\|_\infty,
    \end{align*}
as desired.
\end{proof}

The convolution in Lemma~\ref{Fourier.Connection.lem} may be taken over $\R$ or over
$\T$, the two settings being equivalent since $\ff$ is supported on an interval of
length 1. In fact, the reason we scale $f$ to be supported on an interval of length
$1/2$ is so that we may replace convolution over $\R$, which is the natural place to
study $\De$, with convolution over $\T$, which is the natural place to do harmonic
analysis.

\subsection{The Basic Argument}\label{basic.argument.section}

We begin the process of improving upon the trivial lower bound for $\De$ by stating a
simple version of our method that illustrates the ideas and techniques involved.

\begin{prop}
\label{Cf.First.Result.prop} Let $K$ be any continuous function on $\T$ satisfying
$K(x)\geq1$ when $x\in[-\tfrac14,\tfrac14]$, and let $f$ be a pdf supported on
$[-\tfrac14,\tfrac14]$. Then $$\ffi \geq \|\ff\|_2^2 \geq \|\hat{K}\|_{4/3}^{-4}.$$
\end{prop}

\begin{proof}
We have
$$
1=\int f(x)\,dx \leq \int f(x)K(x)\,dx = \sum_j \hat{f}(j)\hat{K}(-j)
$$
by Parseval's identity. H{\"o}lder's Inequality now gives $1\leq \|\hat{f}\|_4
\|\hat{K}\|_{4/3}$, which we restate as the inequality $\|\hat{K}\|_{4/3}^{-4} \leq
\|\hat{f}\|_4^4$.

Now $\|\hat{f}\|_4^4 = \sum_j |\hat{f}(j)|^4=\sum_j |\widehat{\ff}(j)|^2 = \|\ff\|_2^2$
by another application of Parseval's identity. Since $(\ff)^2\leq \ffi(\ff)$,
integration yields $\|\ff\|_2^2 \leq \ffi\|\ff\|_1 = \ffi$. Combining the last three
sentences, we see that
    $  \|\hat{K}\|_{4/3}^{-4}
                \leq \|\hat{f}\|_4^4
                = \| \ff \|_2^2
                \leq \ffi$
as claimed.
\end{proof}

This reasonably simple theorem already allows us to give a nontrivial lower bound
for~$\De$.

The step function
    $$
    K_1(x) :=
        \begin{cases}
            1   &   0\le|x|\le \frac14\\
            1 - \frac{2\pi^4}{\pi^4 + 24{\zeta(\frac{4}{3})}^3
                \left( 5 + 2^{4/3} - 2^{8/3} \right)}
                     &   \frac14<|x|\le \frac12
        \end{cases}
    $$
has $\|\hat{K_1}\|_{4/3}^{-4} = 1 + \frac{\pi^4} {8\,{\left( 2^{4/3} -1\right) }^3\,
{\zeta(\frac{4}{3})}^3} > 1.074 $ (the elaborate constant used in the definition of
$K_1$ was chosen to minimize $\|\hat{K_1}\|_{4/3}$); a careful reader may complain that
$K_1$ is not continuous. The continuity condition is not essential, however, as we may
approximate $K_1$ by a continuous function $L$ with $\|L\|_{4/3}$ arbitrarily close to
$\|K_1\|_{4/3}$. Green \cite{2001.Green} used a discretization of the kernel function
    $$
    K_2(x) :=
        \begin{cases}
            1   &   0\le|x|\le \frac14\\
            1 - \alpha  +  \alpha \left( 40 (2x - 1)^4 - \frac 32 \right)
                     &   \frac14<|x|\le \frac12
        \end{cases}
    $$
with a suitably chosen $\alpha$ to get $\|\hat{K_2}\|_{4/3}^{-4}> \frac87>1.142$. We get
a slightly larger value of $\|\hat{K}\|_{4/3}^{-4}$ in Corollary~\ref{First.cor} with a
much more complicated kernel. See Section~\ref{kernel.problem.section} for a discussion
of how we came to find our kernel.

\begin{cor}\label{First.cor}
If $f$ is a pdf supported on $[-\frac14,\frac14]$, then
    $$\|\ff\|_2^2 \geq \fstarftwonormconstant.$$
Consequently, $\De \geq 0.574575 \e^2$ for all $0\le\e\le1$.
\end{cor}

\begin{proof}
Set
    \begin{equation}
    K_3(x):=
    \begin{cases}
        1           &   0\leq |x| \leq \tfrac14, \\
        0.6644+0.3356\left(\tfrac{2}{\pi}
        \tan^{-1}\left(\tfrac{1-2x}{\sqrt{4x-1}}\right)\right)^{1.2015}
                    &   \tfrac14 \le |x| \le \tfrac12.
    \end{cases}
    \label{easy.K.def}
    \end{equation}
$K_3(x)$ is pictured in Figure~\ref{GoodK}.
    \begin{figure}
    \begin{center}
    \begin{picture}(384,128)
    \put(176,120){$K_3(x)$}
    \put(350,10){$x$}
    \ifpdf
        \put(0,-240){\includegraphics{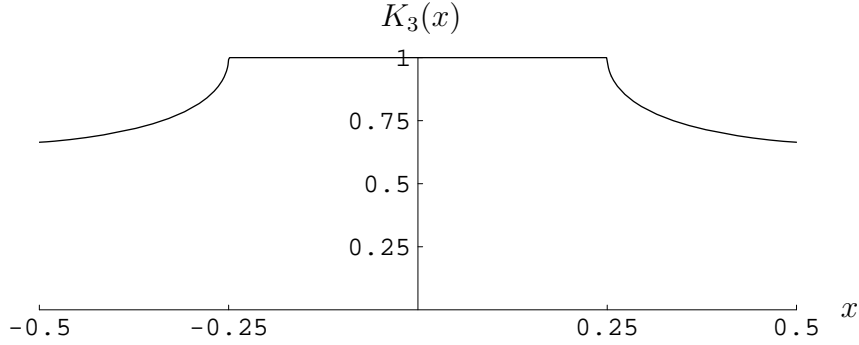}}
    \else
        \put(0,0){\includegraphics{GoodK}}
    \fi
    \end{picture}
    \end{center}
    \caption{The function $K_3(x)$} \label{GoodK}
    \end{figure}

We do not know how to rigorously bound $\|\hat{K_3}\|_{4/3}$, but we can rigorously
bound $\|\hat{K_4}\|_{4/3}$ where $K_4$ is a piecewise linear function `close' to $K_3$.
Specifically, let $K_4(x)$ the even piecewise linear function with corners at
    $$(0,1),\left(\frac14,1\right), \left(\frac14+\frac{t}{4\times 10^4},
    K_3\Big(\frac14+\frac{t}{4\times10^4}\Big)\right) \quad (t=0,1,\dots,10^4).$$
We calculate (using Proposition~\ref{Piecewise.Linear.prop} below) that
$\|\hat{K_4}\|_{4/3} < 0.9658413$. Therefore, by Proposition~\ref{Cf.First.Result.prop}
we have
    $$\|\ff\|_2^2 \geq (0.9658413)^{-4} > \fstarftwonormconstant.$$
Using Lemma~\ref{Fourier.Connection.lem}, we now have $\De >
\frac12\e^2(\fstarftwonormconstant) > 0.574575 \e^2$.
\end{proof}

The constants in the definition \eqref{easy.K.def} of $K_3(x)$ were numerically
optimized to minimize $\|\hat{K_4}\|_{4/3}$ and otherwise have no special significance.
The definition of $K_3(x)$ is certainly not obvious; there are much simpler kernels that
do give nontrivial bounds. In Section~\ref{kernel.problem.section} below, we indicate
the experiments that led to our choice.

We note that the function
    $$
    b(x):=
        \left\{%
        \begin{array}{ll}
            \frac{4/\pi}{\sqrt{1-16x^2}}, & -1/4<x<1/4, \\
            0, & \hbox{otherwise} \\
        \end{array}%
        \right.
    $$
has $\int b=1$ and $\|b\ast b\|_2^2 <1.14939$. Although $b$ is not a pdf (it is not in
$L^2$), it gives strong testimony that the bound on $\|\ff\|_2^2$ given in
Corollary~\ref{First.cor} is not far from best possible.

Note that this corollary, along with Corollary~\ref{R.from.Delta.Connection.cor}, gives
$R(g,n) \leq 1.31925 \sqrt{gn}$, only slightly worse than the bound given
in~\cite{Cilleruelo.Ruzsa.Trujillo}.

This bound on $\|\ff\|_2^2$ may be nearly correct, but the resulting bound on $\ffi$ is
not: we prove below that $\ffi \geq \twotimesDconstant$, and believe that $\ffi \geq
\pi/2$. We have tried to improve the argument given in
Proposition~\ref{Cf.First.Result.prop} in the following four ways:
\begin{enumerate}
\item Instead of considering the sum $\sum_{j} \hat{f}(j)\hat{K}(-j)$ as a whole, we
    single out the central terms (which dominate the sum) and the tails (which
    contribute essentially nothing), and deal with the three resulting sums separately.
    This generalized form of the above argument is expounded in the next section. The
    success of this generalization relies on certain inequalities restricting the possible
    values of these central coefficients; establishing these restrictions is the goal of
    Sections~\ref{useful.inequalities.section} and~\ref{fourier.coefficients.section}. The
    final lower bound derived from these methods is given in
    Section~\ref{full.bound.section}.
\item We can try to find more advantageous kernel functions $K(x)$ for which we can
    compute $\|\hat{K}\|_{4/3}$ in an accurate way. A detailed discussion of our search for
    the best kernel functions is in Section~\ref{kernel.problem.section}.
\item The application of Parseval's identity can be replaced with the Hausdorff-Young
    inequality, which leads to the conclusion $\|\ff\|_\infty \geq \|\hat{K}\|_p^{-q}$,
    where $p\leq \frac{4}{3}$ and $q\geq4$ are conjugate exponents. Numerically, the values
    $(p,q)=(\frac{4}{3},4)$ appear to be optimal. However, Beckner's sharpening
    \cite{Beckner} of the Hausdorff-Young inequality leads to the stronger conclusion
    $\|\ff\|_\infty \geq C(q)\|\hat{K}\|_p^{-q}$ where $C(q) = \frac q2(1-\frac2q)^{q/2-1} =
    \frac q{2e}+O(1)$. We have not experimented to see whether a larger lower bound can be
    obtained from this stronger inequality by taking $q>4$.
\item Notice that we used the inequality $\| g \|_2^2 \leq \|g\|_\infty \|g\|_1$ with
    the function $g=\ff$. This inequality is sharp exactly when the function $g$ takes only
    one nonzero value (i.e., when $g$ is a nif), but the convolution $\ff$ never seems to
    behave that way. Perhaps for these autoconvolutions, an analogous inequality with a
    stronger constant than 1 could be established. Unfortunately, we have not been able to
    realize any success with this idea, although we believe
    Conjecture~\ref{infinity.two.cnj} below. If true, Conjecture~\ref{infinity.two.cnj} implies
    the bound $\De \geq 0.651 \e^2$.
\end{enumerate}
    \begin{cnj}\label{infinity.two.cnj}
    If $f$ is a pdf supported on $[-\frac14,\frac14]$, then
        $$\frac{\ffi }{\|\ff\|_2^2} \geq \frac{\pi/2}{\log 4},$$
    with equality only if either $f(x)$ or $f(-x)$ equals $\sqrt{\frac2{4x+1}}$ on the
    interval $|x|\le\frac14$.
    \end{cnj}

We remark that Proposition~\ref{Cf.First.Result.prop} can be extended from a twofold
convolution in one dimension to an $h$-fold convolution in $d$ dimensions.

\begin{prop}
Let $K$ be any continuous function on $\T^d$ satisfying $K(\bar{x})\geq1$ when
$\bar{x}\in[-\tfrac1{2h},\tfrac1{2h}]^d$, and let $f$ be a pdf supported on
$[-\tfrac1{2h},\tfrac1{2h}]^d$. Then
    $$\| f^{\ast h}\|_\infty \geq \|f^{\ast h}\|_2^2 \geq \|\hat{K}\|_{2h/(2h-1)}^{-2h}.$$
Every subset of $[0,1]^d$ with measure $\e$ contains a symmetric subset with measure
$(0.574575)^d \e^2$.
\end{prop}

\begin{proof}
The proof proceeds as above, with the conjugate exponents $(\frac{2h}{2h-1},2h)$ in
place of $(\frac43,4)$, and the kernel function $K(x_1,x_2,\dots,x_d)=K(x_1)K(x_2)\cdots
K(x_d)$ in place of the kernel function $K(x)$ defined in the proof of
Corollary~\ref{First.cor}. The second assertion of the proposition follows on taking $h=2$.
\end{proof}

\subsection{The Main Bound}\label{main.bound.section}

We now present a more subtle version of Proposition \ref{Cf.First.Result.prop}. Recall
that the notation $\lnorm npa$ was defined in Eq.~(\ref{lnorm.def}). We also use $\Re z$
to denote the real part of the complex number~$z$.

\begin{prop}
\label{Main.Bound.prop} Let $1\leq m < m^{\prime} \leq \infty$. Suppose that $f$ is a
pdf supported on $[-\tfrac14,\tfrac14]$ and that $K$ is even, continuous, satisfies
$K(x) = 1$ for $-\tfrac14 \le x <\frac14$, and $\lnorm{m'}{p}{\hat{K}}>0$. Set $M := 1-
\hat{K}(0) - 2\sum_{j=1}^{m-1} \hat{K}(j) \Re\hat{f}(j)$. Then
\begin{equation}\label{qfinite}
  \|\ff\|_2^2 \geq \sum_{|j| \le m^{\prime}} |\hat{f}(j)|^4 \geq
            1+\bigg(\frac{M}{\lnorm{m}{4/3}{\hat{K}}}\bigg)^4
            +2 \sum_{j=1}^{m-1} |\Re\hat{f}(j)|^4 - o(1)
\end{equation}
as $m^{\prime}\to\infty$.
\end{prop}

\begin{proof}
The first inequality follows from Parseval's formula
    $$\|\ff\|_2^2 = \sum_j |\widehat{\ff}(j)|^2 = \sum_j |\hat{f}(j)|^4
    \geq\sum_{|j| \le m^{\prime}} |\hat{f}(j)|^4.$$

As in the proof of Proposition~\ref{Cf.First.Result.prop}, we have
    $$1         =   \int f(x)K(x)\,dx
                =   \sum_j \hat{f}(j)\hat{K}(-j)
                =   \sum_{|j|<m} \hat{f}(j)\hat{K}(-j) +
                         \sum_{|j|\geq m} \hat{f}(j)\hat{K}(-j).$$
Since $K$ is even, $\hat{K}(-j)=\hat{K}(j)$ is real, and since $f$ is real valued,
$\hat{f}(-j)=\overline{\hat{f}(j)}$. We have
    $$1      =   \hat{K}(0)+2\sum_{j=1}^{m-1} \hat{K}(j) \Re\hat{f}(j)
                    +\sum_{|j|\geq m} \hat{f}(j) \hat{K}(j),$$
which we can also write as $M=\sum_{|j|\geq m} \hat{f}(j)\hat{K}(j)$. Set
$\eta:=\sum_{|j|> m^{\prime}} |\hat{K}(j)|$, and note that $\eta=o(1)$ as $m'\to\infty$ since $K$ is
continuous. Taking absolute values and applying H{\"o}lder's inequality, we have
    \begin{align*}
    |M|     &\leq \sum_{|j|\geq m} |\hat{f}(j)\hat{K}(j)| \\
            &\leq \sum_{m^{\prime}\ge |j|\geq m} |\hat{f}(j)\hat{K}(j)|
                    +\sum_{|j|> m^{\prime}} |\hat{K}(j)| \\
            &\leq \left(\sum_{m^{\prime}\ge|j|\geq m} |\hat{f}(j)|^4\right)^{1/4}
                 \left(\sum_{m^{\prime}\ge|j|\geq m} |\hat{K}(j)|^{4/3}\right)^{3/4}+\eta \\
            &\leq \left(\sum_{m^{\prime}\ge|j|\geq m} |\hat{f}(j)|^4\right)^{1/4}
                    \lnorm{m}{4/3}{\hat{K}}+o(1),
    \end{align*}
which we recast in the form
    $$\sum_{m^{\prime}\ge|j|\geq m} |\hat{f}(j)|^4 \geq
            \bigg(\frac{|M|-o(1)}{\lnorm{m}{4/3}{\hat{K}}} \bigg)^4
            =\bigg(\frac{M}{\lnorm{m}{4/3}{\hat{K}}} \bigg)^4-o(1).$$
We add $\sum_{|j|<m} |\hat{f}(j)|^4$ to both sides and observe that $\hat{f}(0)=1$ and
$|\hat{f}(j)| \geq |\Re \hat{f}(j)|$ to finish the proof of the second inequality.
\end{proof}

%

\begin{cor}\label{overall.twonorm.bound}
If $f$ is a pdf supported on $[-\frac14,\frac14]$, then
    $$\sum_{|j|\le m^{\prime}} |\hat{f}(j)|^4 \geq \fstarftwonormconstant-o(1)$$
 as $m^{\prime}\to\infty$
 \end{cor}

\begin{proof}
The kernel function used in the proof of Corollary~\ref{First.cor} has $\hat{K_4}(0)
\doteq 0.870250799$ and $\lnorm{1}{4/3}{\hat{K_4}} \doteq 0.208784534$. Now apply
Proposition~\ref{Main.Bound.prop} with $m=1.$
\end{proof}

\begin{cor}
    \label{Cf.Main.Bound.cor}
Let $f$ be a pdf supported on $[-\frac14,\frac14]$, and set $x_1:=\Re\hat{f}(1)$. Then
as $m^{\prime}\to\infty$
    \begin{equation*}
    \|\ff\|_2^2 \geq \sum_{|j|\le m^{\prime}} |\hat{f}(j)|^4
        \ge 1 +2 x_1^4
        +\bigg( \frac{0.368067372 - 0.541553784 x_1}
        {0.239175395}\bigg)^4-o(1).
    \end{equation*}
\end{cor}

\begin{proof}
Set
    $$
    K_5(x)=
    \begin{cases}
        1 & |x|\leq\tfrac14, \\
        1 - (1 - (4(\tfrac12 - x))^{1.61707})^{0.546335} & \tfrac14 < |x|\leq \tfrac12.
    \end{cases}
    $$
Denote by $K_6(x)$ the even piecewise linear function with corners at
 $$(0,1), \left(\frac14,1\right),\left(\frac14+\frac{t}{4\times10^4},
    K_5(\frac14+\frac{t}{4\times10^4})\right) \text{(where $t=0,1,\dots,10^4$)}.$$
We find (using Proposition~\ref{Piecewise.Linear.prop}) that $\hat{K_6}(0)\doteq
0.631932628$, $\hat{K_6}(1)\doteq 0.270776892$, and $\lnorm{2}{4/3}{\hat{K_6}}\doteq
0.239175395$. Apply Proposition~\ref{Main.Bound.prop} with $m=2$ to finish the proof.
\end{proof}

With $m^{\prime}=\infty$ and $x_j=\Re \hat{f}(j)$, the bound of
Proposition~\ref{Main.Bound.prop} becomes
    $$ 1+\bigg(\frac{1-\hat{K}(0)-2\sum_{j=1}^{m-1}\hat{K}(j)x_j}
           {\lnorm{m}{4/3}{\hat{K}}} \bigg)^4 + 2 \sum_{j=1}^{m-1} x_j^4.$$
This is a quartic polynomial in the $x_j$, and consequently it is not difficult to
minimize, giving an absolute lower bound on $\| \ff \|_2^2$. This minimum occurs at
    $$x_j=\frac{(\hat{K}(j))^{1/3}
    \left(1-\hat{K}(0)-
    2\sum_{i=1}^{j-1}\hat{K}(i)x_i\right)}{\lnorm{j}{4/3}{\hat{K}}^{4/3}},$$
where $(\hat{K}(j))^{1/3}$ is the real cube root of $\hat{K}(j)$. Consequently,
 $$\inf_{x_j\in\R} \left\{1+\bigg(\frac{1-\hat{K}(0)-2\sum_{j=1}^{m-1}\hat{K}(j)x_j}
                    {\lnorm{m}{4/3}{\hat{K}}}     \bigg)^4
            +2 \sum_{j=1}^{m-1} x_j^4 \right\}=
            1+\left( \frac{1-\hat{K}(0)}{\lnorm{1}{4/3}{\hat{K}}} \right)^4,$$
which is nothing more than the bound that Proposition~\ref{Main.Bound.prop} gives with
$m=1$. Moreover,
    $$ 1+\left( \frac{1-\hat{K}(0)}{\lnorm{1}{4/3}{\hat{K}}} \right)^4
        = \sup_{0\le \alpha \le 1} \| (\alpha+(1-\alpha)K)^{\wedge} \|_{4/3}^{-4}, $$
(the details of this calculation are given in Section~\ref{adhoc.section}) so that
Proposition~\ref{Main.Bound.prop} by itself does not give a different bound on $\| \ff
\|_2^2$ than Proposition~\ref{Cf.First.Result.prop}. However, we shall obtain additional
information on $\hat{f}(j)$ in terms of $\ffi$ in
subsection~\ref{fourier.coefficients.section} below, and this information can be
combined with Proposition~\ref{Main.Bound.prop} to provide a stronger lower bound on
$\ffi$ than that given by Proposition~\ref{Cf.First.Result.prop}.

\subsection{Some Useful Inequalities}\label{useful.inequalities.section}

Hardy, Littlewood, and P{\'o}lya~\cite{1988.Hardy.Littlewood.Polya} call a function $u(x)$
{\it symmetric decreasing\/} if $u(x)=u(-x)$ and $u(x) \geq u(y)$ for all $0\leq x\leq
y$, and they call
 $$\sdr{f}(x):= \inf\left\{ y \colon\quad
      \lambda\left(\left\{t \colon f(t)\geq y \right\}\right)\le 2|x| \right\}$$
the {\it symmetric decreasing rearrangement\/} of $f$. For example, if $f$ is the
indicator function of a set with measure $\mu$, then $\sdr{f}$ is simply the indicator
function of the interval $(-\frac\mu2,\frac\mu2)$. Another example is any function $f$
defined on an interval $[-a,a]$ and is periodic with period $\frac{2a}n$, where $n$ is a
positive integer, and that is symmetric decreasing on the subinterval $[-\frac an,\frac
an]$; then $\sdr{f}(x) = f(\frac xn)$ for all $x\in[-a,a]$. In particular, on the
interval $[-\frac14,\frac14]$, we have $\sdr{\cos}(2\pi j x)=\cos(2\pi x)$ for any
nonzero integer $j$. We shall need the following result.

\begin{lem}
$$\int f(x)u(x)\,dx \leq \int \sdr{f}(x)\sdr{u}(x)\,dx.$$
    \label{1988.Hardy.Littlewood.Polya.lem}
\end{lem}

\begin{proof}
This is Theorem 378 of~\cite{1988.Hardy.Littlewood.Polya}.
\end{proof}

We say that {\em $\bar{f}$ is more focused than $f$} (and $f$ is less focused than
$\bar{f}$) if for all $z\in[0,\frac12]$ and all $r\in\T$ we have
 $$\int_{r-z}^{r+z} f \leq \int_{-z}^z \bar{f}.$$
For example, $\sdr{f}$ is more focused than $f$. In fact, we introduce this terminology
because it refines the notion of symmetric decreasing rearrangement in a way that is
useful for us. To give another example, if $f$ is a nonnegative function, set $\bar f$
to be $\|f\|_\infty$ times the indicator function of the interval
$[-\frac1{2\|f\|_\infty},\frac1{2\|f\|_\infty}]$; then $\bar f$ is more focused than
$f$.

\begin{lem}
    \label{Focus.lem}
Let $u(x)$ be a symmetric decreasing function, and let
$h,\bar{h}$ be pdfs with $\bar{h}$ more focused than $h$. Then for all $r\in\T$,
    $$ \int h(x-r) u(x) \, dx \leq \int \bar{h}(x) u(x) \, dx.$$
\end{lem}

\begin{proof}
Without loss of generality we may assume that $r=0$, since if $\bar h(x)$ is more
focused than $h(x)$, then it is also more focused than $h(x-r)$. Also, without loss of
generality we may assume that $h,\bar h$ are continuous and strictly positive on $\T$,
since any nonnegative function in $L^1$ can be $L^1$-approximated by such.

Define $H(z) = \int_{-z}^z h(t)\,dt$ and $\bar H(z) = \int_{-z}^z \bar h(t)\,dt$, so
that $H(\frac12) = \bar H(\frac12) = 1$, and note that the more-focused hypothesis
implies that $H(z)\le\bar H(z)$ for all $z\in[0,\frac12]$. Now $h$ is continuous and
strictly positive, which implies that $H$ is differentiable and strictly increasing on
$[0,\frac12]$ since $H'(z)=h(z)+h(-z)$. Therefore $H^{-1}$ exists as a function from
$[0,1]$ to $[0,\frac12]$. Similar comments hold for $\bar H^{-1}$.

Since $H\le\bar H$, we see that $\bar H^{-1}(s) \le H^{-1}(s)$ for all $s\in[0,1]$.
Then, since $H^{-1}(s)$ and $H^{-1}(s)$ are positive and $u$ is decreasing for positive
arguments, we conclude that $u(H^{-1}(s)) \le u(\bar H^{-1}(s))$, and so
    \begin{equation}
    \int_0^1 u(H^{-1}(s))\,ds \le \int_0^1 u(\bar H^{-1}(s))\,ds.
    \label{H.inverse.inequality}
    \end{equation}
On the other hand, making the change of variables $s=H(t)$, we see that
    $$
    \int_0^1 u(H^{-1}(s))\,ds = \int_0^{H^{-1}(1)} u(t) H'(t)\,dt = \int_0^{1/2}
    u(t) (h(t)+h(-t))\,dt = \int_{\T} u(t)h(t)\,dt
    $$
since $u$ is symmetric. Similarly $\int_0^1 u(\bar H^{-1}(s))\,ds = \int_{\T} u(t)\bar
h(t)\,dt$, and so inequality (\ref{H.inverse.inequality}) becomes $\int u(t)h(t)\,dt \le
\int u(t)\bar h(t)\,dt$ as desired.
\end{proof}

\subsection{Fourier Coefficients of Density Functions}\label{fourier.coefficients.section}

To use Proposition~\ref{Main.Bound.prop} to bound $\De$, we need to develop a better
understanding of the central Fourier coefficients $\hat{f}(j)$ for small $j$. In
particular, we wish to apply Proposition~\ref{Main.Bound.prop} with $m=2$, i.e., we need
to develop the connections between $\ffi$ and the real part of the Fourier coefficient
$\hat{f}(1)$.

We turn now to bounding $|\hat{f}(j)|$ in terms of $\ffi$. The guiding principle is that
if $\ff$ is very concentrated then $\ffi$ will be large, and if $\ff$ is not very
concentrated then $|\hat{f}(j)|$ will be small. Green \cite[Lemma 26]{2001.Green} proves
the following lemma in a discrete setting, but since we need a continuous version we
include a complete proof.

\begin{lem}
\label{Fhat.Green.Bound.lem} Let $f$ be a pdf supported on $[-\frac14,\frac14]$. For
$j\not=0$, $$|\hat{f}(j)|^2 \leq \frac{\ffi}{\pi} \,\sin\left(\frac{\pi}{\ffi}\right).$$
\end{lem}

\begin{proof}
Let $f_1:\T\to\R$ be defined by $f_1(x):=f(x-x_0)$, with $x_0$ chosen so that
$\hat{f_1}(j)$ is real and positive (clearly $\hat{f_1}(j)=|\hat{f}(j)|$ and
$\ffi=\|f_1\ast f_1\|_\infty$). Set $h(x)$ to be the symmetric decreasing rearrangement
of $f_1\ast f_1$, and $\overline{h}(x):=\ffi I(x)$, where $I(x)$ is the indicator
function of $[-\frac{1}{2 \ffi},\frac{1}{2 \ffi}]$. We have
 \begin{align*}
 |\hat{f}(j)|^2 &=      \hat{f_1}(j)^2
                =      \widehat{f_1\ast f_1}(j)
                =      \int f_1\ast f_1(x) \cos(2\pi j x)\,dx
                \leq    \int h(x) \cos(2\pi x)\,dx
 \end{align*}
by Lemma~\ref{1988.Hardy.Littlewood.Polya.lem}. We now apply Lemma~\ref{Focus.lem} to
find
 \begin{multline*}
 |\hat{f}(j)|^2
    \leq    \int \overline{h}(x)\cos(2\pi x)\,dx \\
    =       \int_{-1/(2\ffi)}^{1/(2\ffi)} \ffi \cos(2\pi x)\,dx
    =       \frac{\ffi}{\pi}\sin\left(\frac{\pi}{\ffi}\right).
 \end{multline*}
\end{proof}

\subsection{The Full Bound}\label{full.bound.section}

With the technical result of Section~\ref{fourier.coefficients.section} in hand, we can
finally establish the lower bound on $\De$ given in Theorem~\ref{Delta.Summary.thm}(ii).

\begin{prop}
\label{the.full.De.bound.prop} $\De \ge \Dconstant \e^2$ for all $0\le\e\le1$.
\end{prop}

\begin{proof}
Let $f$ be a pdf supported on $[-\frac14,\frac14]$, and assume that
    \begin{equation}\label{Main.Bound.0.eq}
    \ffi<\twotimesDconstant.
    \end{equation}
Set $x_1=\Re\hat{f}(1)$ and $x_2=\Re\hat{f}(2)$. Since $f$ is supported on
$[-\frac14,\frac14]$, we see that $x_1>0$. By Lemma~\ref{Fhat.Green.Bound.lem},
    \begin{equation}\label{Main.Bound.1.eq}
    0< x_1 < 0.4191447.
    \end{equation}
Corollary~\ref{Cf.Main.Bound.cor} with $m^{\prime}=\infty$, gives
    \begin{equation}\label{Main.Bound.3.eq}
    \ffi\geq \|\ff\|_2 \geq
    1 +2 x_1^4 +\bigg( \frac{0.368067372 - 0.541553784 x_1}{0.239175395}\bigg)^4.
    \end{equation}
Routine calculus shows that there are no simultaneous solutions to
Inequalities~(\ref{Main.Bound.0.eq}), (\ref{Main.Bound.1.eq}), and
(\ref{Main.Bound.3.eq}). Therefore $\ffi\geq\twotimesDconstant$, whence Lemma~\ref{Fourier.Connection.lem} implies that
$\De\ge \Dconstant \e^2$.
\end{proof}

This gist of the proof of Proposition~\ref{the.full.De.bound.prop} is that if $\ffi$ is
small, then $\Re\hat{f}(1)$ is small by Lemma~\ref{Fhat.Green.Bound.lem}, and so
$\|\ff\|_2^2$ is not very small by Corollary~\ref{Cf.Main.Bound.cor}, whence $\ffi$ is not small. If $\ffi \leq \twotimesDconstant$, then we get a contradiction. We can actually prove a meaningful result about $\|\ff\|_2^2$ under the condition that $\ffi$ is not much larger than $\twotimesDconstant$. The following result will be useful in
Section~\ref{upper.bounds.Rgn.from.De.section}.

\begin{lem}\label{fft.bound.using.ffi}
Let $f$ be a pdf supported on $[-\frac14,\frac14]$. If $\twotimesDconstant \leq \ffi
\leq 1.229837$, then
    $$\|\ff\|_2^2 \geq \sum_{|j|\le m^{\prime}} |\hat{f}(j)|^4
        \ge 21.922911 - 33.711941 \ffi + 13.676987 \ffi^2-o(1)$$
as $m^{\prime}\to\infty$.
\end{lem}

\begin{proof}
The first inequality has already been shown in Proposition \ref{Main.Bound.prop}.
Set $$B(x_1):=1 +2 x_1^4 +\bigg( \frac{0.368067372 - 0.541553784
x_1}{0.239175395}\bigg)^4,$$ so that by Corollary~\ref{Cf.Main.Bound.cor},
    $$\sum_{|j|\geq m^{\prime}} |\hat{f}(j)|^4 \ge B(\Re\hat{f}(1))-o(1).$$
By hypothesis $\ffi\leq 1.229837$, so that by Lemma~\ref{Fhat.Green.Bound.lem},
    $$\Re\hat{f}(1) \leq \sqrt{\frac{\ffi}{\pi}\sin\left(\frac{\pi}{\ffi}\right)} < 0.466.$$
But $B(x_1)$ is a decreasing function for $x_1 \in [0,0.47]$, so that
    \begin{align*}
    B(\Re\hat{f}(1))
    &\ge B\left( \sqrt{\tfrac{\ffi}{\pi}\sin\left(\tfrac{\pi}{\ffi}\right)}\right) \\
    &> 21.922911 - 33.711941 \ffi + 13.676987 \ffi^2
    \end{align*}
after a straightforward computation.
\end{proof}

\subsection{The Kernel Problem}\label{kernel.problem.section}

Let $\K$ be the class of functions $K\in L^2$ satisfying $K(x)\ge1$ on
$[-\frac14,\frac14]$. Proposition~\ref{Cf.First.Result.prop} suggests the problem
of computing
\begin{equation*}
\inf_{K\in\K} \|\hat{K}\|_p =
 \inf_{K\in\K}
    \left(\sum_{j=-\infty}^\infty |\hat{K}(j)|^p\right)^{1/p}.
\end{equation*}
In Proposition~\ref{Cf.First.Result.prop} the case $p=\frac43$ arose, but using the
Hausdorff-Young inequality in place of Parseval's identity we are led to consider $1< p
\leq \frac43$. Also, we assumed in Proposition~\ref{Cf.First.Result.prop} that $K$ was
continuous, but this assumption can be removed by taking the pointwise limit of
continuous functions.

As similar problems occur in~\cite{Cilleruelo.Ruzsa.Trujillo} and in~\cite{2001.Green}, we feel it is worthwhile to detail the thoughts and experiments that led to the kernel functions chosen in Corollaries~\ref{First.cor} and~\ref{Cf.Main.Bound.cor}.

\subsubsection{{\em Ad hoc} Observations}\label{adhoc.section}

Our first observation is that if $G\in\K$, then so is
$K(x):=\frac12(G(x)+G(-x))$, and since $|\hat{K}(j)|=|\Re\hat{G}(j)|\leq
|\hat{G}(j)|$ we know that $\|\hat{K}\|_p\le\|\hat{G}\|_p$. Thus, we may
restrict our attention to the {\em even} functions in $\K$.

We also observe that $|\hat{K}(j)|$ decays more rapidly if many derivatives of $K$ are
continuous. This suggests that we should restrict our attention to continuous $K$,
perhaps even to infinitely differentiable $K$. However, computations suggest that the
best functions $K$ are continuous but {\em not} differentiable at $x=\frac14$ (see in
particular Section~\ref{piecewise.linear.section} and Figure~\ref{optimal.K.p43.fig}).

In the argument of Proposition~\ref{Cf.First.Result.prop} we used the
inequality $\int f \leq \int fK$, which is an equality if we take $K$ to be
equal to 1 on $[-\frac14,\frac14]$, instead of merely at least 1. In light of
this, we should not be surprised if the optimal functions in $\K$ are exactly
1 on $[-\frac14,\frac14]$. This is supported by our computations.

Finally, we note that if $K_i\in\K$, and $\alpha_i>0$ with $\sum_i \alpha_i = 1$, then
$\sum_i \alpha_i K_i(x)\in \K$ also. This is particularly useful with $K_1(x):=1$.
Specifically, given any $K_2\in\K$ with known $\|\hat{K_2}\|_p$ (we stipulate
$\|K_2\|_1=\hat K(0)=\le1$ to avoid technicalities), we may easily compute the
$\alpha\in[0,1]$ for which $\|\hat{K}\|_p$ is minimized, where $K(x):=\alpha K_1(x) +
(1-\alpha)K_2(x)$. We have
\begin{equation}
    \|\hat K\|_p^p = (\alpha + (1-\alpha)\hat K_2(0) )^p + (1-\alpha)^p
    \lnorm{1}{p}{\hat K}^p = (1-(1-\alpha)M)^p + (1-\alpha)^p N,
    \label{hatKpp.expression}
\end{equation}
where we have set $M := 1-\hat K_2(0)$ and $N := \lnorm{1}{p}{\hat K}^p$. Taking the
derivative with respect to $\alpha$, we obtain
\begin{equation*}
p(1-\alpha)^{p-1} \bigg( M \Big( \frac1{1-\alpha}-M \Big)^{p-1} - N \bigg),
\end{equation*}
the only root of which is $\alpha = 1 - \frac{M^{q/p}}{M^q + N^{q/p}}$ (where
$\frac1p+\frac1q=1$). It is straightforward (albeit tedious) to check by substituting
$\alpha$ into the second derivative of the expression~(\ref{hatKpp.expression}) that
this value of $\alpha$ yields a local maximum for $\|\hat K\|_p^p$. The maximum value
attained is then calculated to equal $N \big( M^q + N^{q/p} \big)^{1-p}$, which is
easily computed from the known function $K_2$.

Notice that when $p=\frac43$ (so $q=4$), applying Proposition
\ref{Cf.First.Result.prop} with our optimal function $K$ yields
\begin{equation*}
    \|\ff\|_2^2 \ge \|\hat K\|_{4/3}^{-4}
    = \Big( N \big( M^4 + N^3 \big)^{-1/3} \Big)^{-3} =
    \frac{M^4 + N^3}{N^3} = 1 + \frac{(1-\hat K_2(0))^4}{\lnorm{1}{4/3}{\hat K}^4},
\end{equation*}
whereupon we recover the conclusion of Proposition~\ref{Main.Bound.prop} with $m=1$.

\subsubsection{Trigonometric Polynomials}

We wish to identify families of functions that are at least 1 on
$[-\tfrac14,\tfrac14]$ and whose Fourier coefficients have small $\ell^p$
norm. Natural candidates are functions which have many Fourier coefficients
equal to 0. In this section we consider trigonometric polynomials
$K(x)=\sum_{j=-m}^m \hat{K}(j)e^{2\pi i j x}$ of degree $m$.

Montgomery~\cite[Chapter 1]{1994.Montgomery} defines the Selberg polynomials
$S_m^+(\alpha,\beta,x)$ and shows that $S_m^+(\alpha,\beta,x)\geq
\chi_{[\alpha,\beta]}(x)$ for all $x$, provided that $\alpha\leq\beta\leq\alpha+1$;
moreover, these functions are (in some senses) optimal $L^1$ majorants for
$\chi_{[\alpha,\beta]}(x)$ among all trigonometric polynomials of bounded degree. These
provide a natural family for investigation.

We are concerned with $[\alpha,\beta]=[-\tfrac14,\tfrac14]$. We have for instance
$$S_2^+(-\tfrac14,\tfrac14,x) = \frac56 + \big( \frac4{9\sqrt3} + \frac2{3\pi} \big) \cos
(2\pi x) - \frac29 \cos (4\pi x),$$ which satisfies
$S_2^+(-\tfrac14,\tfrac14,x)-\frac1{18} \ge 1$ when $x\in[-\frac14,\frac14]$, and
    \begin{equation*}
    \inf_{0\le \alpha \le 1}
    \| (\alpha+(1-\alpha)(S_2^+(-\tfrac14,\tfrac14,x)-\frac1{18})^{\wedge} \|_{4/3}
    >0.990.
    \end{equation*}
However, $L(x):=2\cos(2\pi x)-\cos(4\pi x)\in \K$, and
    \begin{equation*}
    \inf_{0\le \alpha \le 1}
    \| (\alpha+(1-\alpha)(L(x))^{\wedge} \|_{4/3}
    <0.989.
    \end{equation*}
Thus, even among trigonometric polynomials of degree 2, the Selberg polynomials are not
optimal. In general, we have been unable to identify the degree-$m$ trigonometric
polynomial $K(x)$ that is in $\K$ and for which $\sum_{j=-m}^m |\hat{K}(j)|^{4/3}$ is
minimized.

\subsubsection{Wavelets}

We can give an exact, finite expression for the $p$-norm of the Fourier coefficients of
some large classes of functions. Sums of Haar wavelets give the simplest theoretical
instance and the largest class of functions.

Define $$\psi(x):=
  \begin{cases}
    1 & 0\leq x < \frac12, \\
    -1 & \frac12 \leq x < 1,\\
    0& \text{otherwise}.
  \end{cases}$$
and $\psi_{m,n}(x):=2^{-m/2} \psi(2^m x-n)$. It is well-known that the Haar Wave\-lets
$\{\psi_{m,n}\colon m,n\in\Z\}$ form an orthonormal basis of the subspace of $L^2(\R)$
consisting of functions that have integral 0. By the comments in
Section~\ref{adhoc.section},
 $$\inf_{K \in \K} \|\hat{K}\|_{4/3}
 = \inf_{\substack{K\in \K \\ \int K = 0}}
    \left( 1+ \frac{1}{{\lnorm{1}{4/3}{\hat K}^4}}\right)^{-1/4},$$
so that the $\int K =0$ restriction is not a substantial restriction.

\begin{figure}
\begin{center}
 \begin{picture}(384,360)
    \put(24,108){$K(x)$}
    \put(24,228){$K(x)$}
    \put(24,348){$K(x)$}
    \put(215,108){$K(x)$}
    \put(214,228){$K(x)$}
    \put(213,348){$K(x)$}
    \put(178,64){$x$}
    \put(178,184){$x$}
    \put(178,304){$x$}
    \put(372,64){$x$}
    \put(372,184){$x$}
    \put(372,304){$x$}
    \put(85,108){$N=31$}
    \put(86,228){$N=7$}
    \put(86,348){$N=1$}
    \put(275,108){$N=63$}
    \put(275,228){$N=15$}
    \put(276,348){$N=3$}
    \includegraphics{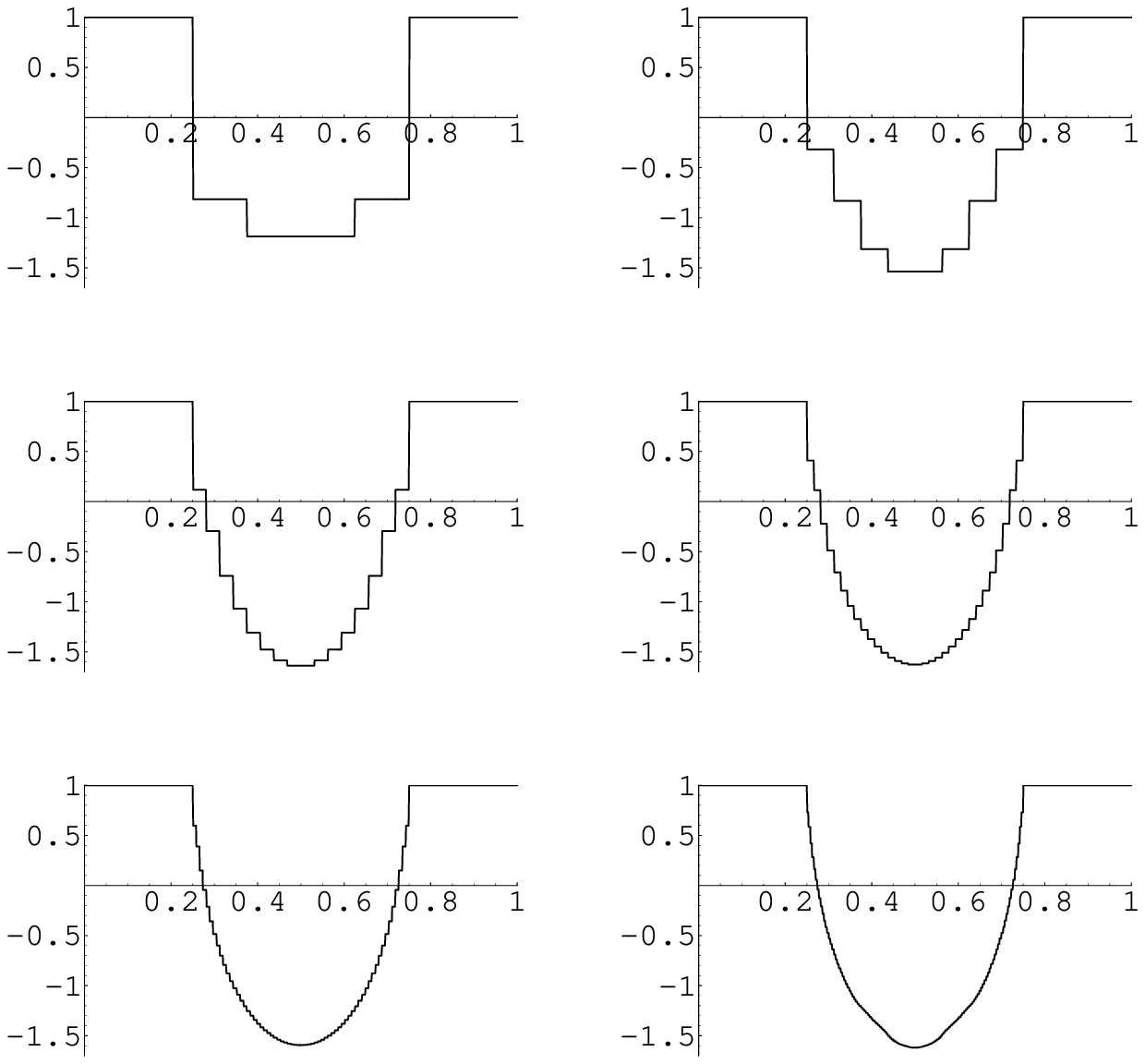}
 \end{picture}
 \begin{picture}(384,128)
    \put(25,120){$K(x)$}
    \put(186,120){$N=127$}
    \put(378,66){$x$}
    \ifpdf
        \put(0,-177){\includegraphics{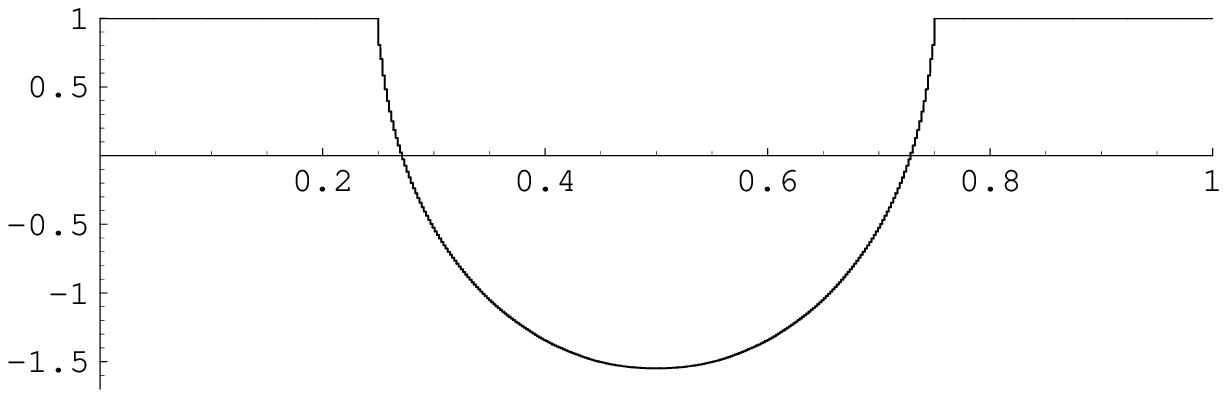}}
    \else
        \put(0,0){\includegraphics{graph7}}
    \fi
 \end{picture}
 \end{center}
    \caption{Optimal kernels generated by Haar wavelets}
    \label{optimal.K.p43.fig}
\end{figure}

It follows that every even function $K\in\K$ with $\int K = 0$ and
$x\in[-\frac14,\frac14)\Rightarrow K(x)=1$ can be written in the form
 $$K(x)=\sqrt{2}\left(\psi_{1,0}(x)+\psi_{1,0}(-x)\right)+
        \sum_{n=1}^\infty \alpha_n \left(
        \psi_{2+\floor{\log_2 n},n}(x)+\psi_{2+\floor{\log_2 n},n}(-x)
        \right).$$
Since this expression is linear, we can give the Fourier coefficients of $K$ in terms of
the easily computable Fourier coefficients of the $\psi_{m,n}$ and in terms of the
parameters $\alpha_n$. Truncating the infinite sum at $N$, we obtain a reasonably large
family of possible functions $K$ for which $\|\hat{K}\|_p$ can be computed quickly
enough to numerically optimize $\alpha_1,\dots,\alpha_{N}$. Graphs of the optimal $K(x)$
for $p=\frac43$ and various values of $N$ are displayed in
Figure~\ref{optimal.K.p43.fig}. Note that if these wavelet-based functions are converging to some limit function in $\K$, that limit function certainly does not seem to be differentiable at $\pm\frac14$.

While the exact form of this expression for $\|\hat K\|_p$ is not difficult to
compute, we suppress the details as they are very similar to the expressions
computed in the next section.

\subsubsection{Piecewise-Linear Functions}\label{piecewise.linear.section}

More useful computationally is the class of continuous piecewise-linear even functions
whose vertices all have abscissae with a given denominator. Let $\zeta(s,a)
:=\sum_{k=0}^\infty (k+a)^{-s}$ denote the Hurwitz zeta function. If ${\bf v}$ is a
vector, define $\Lambda_p({\bf v})$ to be the vector whose coordinates are the $p$th
powers of the absolute values of the corresponding coordinates of~${\bf v}$.

\begin{prop}
Let $T$ be a positive integer, $n$ a nonnegative integer, and $p\ge1$ a real
number. For each integer $0\le t\le T$, define $x_t:=\frac14+\frac{t}{4T}$, and
let $y_t$ be an arbitrary real number, except that $y_0=1$. Let $K(x)$ be the even
function on $\T$ that is linear on $[0,\frac14]$ and each of the intervals
$[x_{t-1},x_t]$ $(1\leq t \leq T)$, satisfying $K(0)=1$ and $K(x_t)=y_t$ $(0\leq
t \leq T)$. Then
    $$
    \lnorm{n}{p}{\hat{K}} =( 2 \Lambda_p(\vecd\A) \cdot\vecz )^{1/p},
    $$
where $\vecd$ is the $T$-dimensional vector $\vecd=(y_1-y_0, y_2-y_1, \dots,
y_T-y_{T-1})$, $\A$ is the $T\times4T$ matrix whose $(t,k)$-th component is
    $$
    \A_{tk} = \cos(2\pi(n+k-1)x_t)-\cos(2\pi(n+k-1)x_{t-1}),
    $$
and $\vecz$ is the $4T$-dimensional vector
    $$
    \vecz = (8T\pi^2)^{-p} \big( \zeta(2p,\tfrac{j}{4T}),
    \zeta(2p,\tfrac{j+1}{4T}), \dots, \zeta(2p,\tfrac{j+4T-1}{4T}) \big).
    $$
\label{Piecewise.Linear.prop}
\end{prop}

\begin{proof}
Note that
 \begin{align*}
   \hat{K}(-j)  =   \hat{K}(j)
                &=   \int_{-1/2}^{1/2} K(u)\cos(2\pi j u)\,du  \\
                &=   2\int_{0}^{1/4}\cos(2\pi j u)\,du+
                    2\sum_{t=1}^T\int_{x_{t-1}}^{x_t}(m_t u+b_t)\cos(2\pi j u)\,du,\\
 \end{align*}
where $m_t$ and $b_t$ are the slope and $y$-intercept of the line going
through $(x_{t-1},y_{t-1})$ and $(x_t,y_t)$. If we define $C(j):=\tfrac{\pi^2
j^2}{2T}\hat{K}(j)$, then integrating by parts we have
\begin{multline}
    C(j)= \left( \frac{\pi^2 j^2}{T}
    \left(\left. \frac{1}{2\pi j} \sin\left(2\pi j u\right)\right|_{0}^{1/4} \right)
    + \frac{\pi^2 j^2}{T}\sum_{t=1}^T \left( \left.  \frac{m_t u+b_t}{2\pi
    j}\sin(2\pi j u)\right|_{x_{t-1}}^{x_t}\right) \right) \\
    + \left( \frac{\pi^2 j^2}{T} \sum_{t=1}^T \left.\frac{m_t}{(2\pi j)^2}
    \cos(2\pi j u)\right|_{x_{t-1}}^{x_t} \right).
    \label{first.term.O1}
\end{multline}
The first term of this expression is
\begin{multline*}
    \frac{\pi j}{2T} \bigg( \sin\left( {\pi \tfrac j2} \right) + \sum_{t=1}^T
    \big(  (m_t x_t+b_t)\sin(2\pi j x_t) - (m_t x_{t-1}+b_t)\sin(2\pi
    j x_{t-1}) \big) \bigg) \\
    = \frac{\pi j}{2T} \bigg( \sin\left( {\pi \tfrac j2} \right) + \sum_{t=1}^T (m_t
    x_t+b_t)\sin(2\pi j x_t) - \sum_{t=0}^{T-1} (m_{t+1} x_t+b_{t+1})\sin(2\pi j x_t) \bigg).
\end{multline*}
Since $m_{t+1} x_t+b_{t+1} = y_t = m_t x_t+b_t$ and $x_0=\frac14$,
$x_T=\frac12$, this entire expression is a telescoping sum whose value is zero.
Eq.~(\ref{first.term.O1}) thus becomes
\begin{align}
    C(j) &= \frac{\pi^2 j^2}{T} \sum_{t=1}^T \left.\frac{m_t}{(2\pi j)^2} \cos(2\pi j
    u)\right|_{x_{t-1}}^{x_t} \notag \\
            &=   \sum_{t=1}^T    \left(y_t-y_{t-1}\right)
                                \left(\cos(2\pi j x_t)-\cos(2\pi j x_{t-1})\right)
                                \label{How.to.compute.C}
\end{align}
using $m_t=\frac{y_t-y_{t-1}}{x_t-x_{t-1}} = 4T(y_t-y_{t-1})$. Each $x_t$ is rational
and can be written with denominator $4T$, so we see that the sequence of normalized
Fourier coefficients $C(j)$ is periodic with period $4T$.

We proceed to compute $\lnorm{n}{p}{\hat{K}}$ with $n$ positive and
$p\geq1$.
\begin{align*}
    \left(\lnorm{n}{p}{\hat{K}}\right)^p
         &=  \sum_{|j|\geq n} |\hat{K}(j)|^p
         =  2\sum_{j=n}^\infty |\hat{K}(j)|^p
         =  2\sum_{j=n}^\infty \left|C(j)\frac{2T}{\pi^2 j^2}\right|^p \\
        &=  2\left(\frac{2T}{\pi^2}\right)^p
                    \sum_{j=n}^\infty \frac{|C(j)|^p}{j^{2p}}.
\end{align*}
Because of the periodicity of $C(j)$, we may write this as
\begin{align}
\left(\lnorm{n}{p}{\hat{K}}\right)^p
        &=  2\left(\frac{2T}{\pi^2}\right)^p
                    \left(\sum_{j=n}^{n+4T-1} |C(j)|^p
                            \sum_{r=0}^\infty (4Tr+j)^{-2p}\right) \notag \\
        &=  2\left(\frac{2T}{(4T)^2\pi^2}\right)^p
                    \left( \sum_{j=n}^{n+4T-1} |C(j)|^p
                            \zeta\left(2p,\tfrac{j}{4T}\right)
                            \right),\label{piecewise.linear.eq}
\end{align}
which concludes the proof.
\end{proof}

Proposition~\ref{Piecewise.Linear.prop} is useful in two ways. The first is that only
$\vecd$ depends on the chosen values $y_t$. That is, the vector $\vecz$ and the matrix
$\A$ may be precomputed (assuming $T$ is reasonably small), enabling us to compute
$\lnorm{n}{p}{\hat{K}}$ quickly enough as a function of $\vecd$ to numerically optimize
the $y_t$. The second use is through Eq.~(\ref{piecewise.linear.eq}). For a given $K$,
we set $y_t=K(x_t)$, whereupon $C(j)$ is computed for each $j$ using the formula in
Eq.~\eqref{How.to.compute.C}. Thus we can use Eq.~(\ref{piecewise.linear.eq}) to compute
$\lnorm{n}{p}{\hat{K_1}}$ with arbitrary accuracy, where $K_1$ is almost equal to $K$.
We have found that with $T=10000$ one can generally compute $\lnorm{n}{p}{\hat{K_1}}$
quickly.

In performing these numerical optimizations, we have found that ``good'' kernels
$K(x)\in \K$ have a very negative slope at $x=\frac14 ^+$ (e.g., see Figure
\ref{optimal.K.p43.fig}). Viewing graphs of these numerically optimized kernels suggests
that functions of the form
    $$
    K_{d_1,d_2}(x)=
    \begin{cases}
        1 & |x|\leq\tfrac14, \\
        1 - (1 - (4(\tfrac12 - x))^{d_1})^{d_2} & \tfrac14 < |x|\leq \tfrac12,
    \end{cases}
    $$
which have slope $-\infty$ at $x=\frac14^+$, may be very good. (Note that the graph of
$K_{2,1/2}(x)$ between $\frac14$ and $\frac34$ is the lower half of an ellipse.) More
good candidates are functions of the form
    $$
    K_{e_1,e_2,e_3}(x) =
    \begin{cases}
        1
        & |x|\leq\tfrac14, \\
        \left(\frac2\pi \tan^{-1}\left(\frac{(1-2x)^{e_1}}{(4x-1)^{e_2}}\right)\right)^{e_3}
        & \tfrac14 < |x|\leq \tfrac12,
    \end{cases}
    $$
where $e_1,e_2,$ and $e_3$ are positive. We have used a function of the form $K_{d_1,d_2}$ in the proof
of Corollary~\ref{Cf.Main.Bound.cor} and a function of the form $K_{e_1,e_2,e_3}$ in the proof of
Corollary~\ref{First.cor}.

\subsection{A Lower Bound for $\Delta(\frac12)$}\label{delta.half.section}

We begin with a fundamental relationship between $\Re \hat{f}(1)$ and $\Re \hat{f}(2)$.

\begin{lem}
    \label{Fhat.F1.F2.Domain.lem}
Let $f$ be a pdf supported on $[-\tfrac14,\tfrac14]$. Then
    $$2\big(\Re\hat{f}(1)\big)^2-1 \leq \Re\hat{f}(2) \leq 2(\Re\hat{f}(1))-1.$$
\end{lem}

\begin{proof}
Since $L_2(x):=2\cos(2\pi x)-\cos(4\pi x)$ is at least 1 for $-\tfrac14\leq x
\leq\tfrac14$, we have
    $$1 \leq \int f(x)L(x)\,dx=\sum_{j=-2}^2
    \hat{f}(j)\hat{L}(-j)=2(\Re\hat{f}(1))-\Re\hat{f}(2).$$
Rearranging, we arrive at $\Re\hat{f}(2) \leq 2(\Re\hat{f}(1))-1$.

We give two proofs of the inequality $2\big(\Re\hat{f}(1)\big)^2-1 \leq
\Re(\hat{f}(2))$, each with its own advantages.

{\em First proof:} Since $f(x)$ and $\cos(2\pi x)$ are both nonnegative on
$[-\frac14,\frac14]$, the Cauchy-Schwartz inequality gives
 \begin{align*}
    2(\Re \hat{f}(1))^2-1
    &= 2 \left(\int_{-1/4}^{1/4} f(x) \cos(2\pi x)\,dx \right)^2 -
       \int_{-1/4}^{1/4} f(x) \, dx \\
    &\leq
            2\left(\int_{-1/4}^{1/4} f(x)\,dx\right)
             \left(\int_{-1/4}^{1/4} f(x) \left(\cos(2\pi x)\right)^2\,dx\right)
             -\int_{-1/4}^{1/4} f(x) \, dx\\
    &=  2\int_{-1/4}^{1/4} f(x) \left(\cos (2\pi x)\right)^2\,dx
             - \int_{-1/4}^{1/4} f(x) \,dx \\
    &=  \int_{-1/4}^{1/4} f(x)\left(2\cos^2 (2\pi x)-1\right)\,dx \\
    &=  \int_{-1/4}^{1/4} f(x)\cos (4\pi x)\,dx  =  \Re \hat{f}(2).
 \end{align*}

{\em Second proof:} Set $L_b(x)=b \cos(2\pi x)-\cos(4\pi x)$ (with $b\ge 0$) and observe
that for $-\tfrac14\leq x \leq\tfrac14$, we have $L_b(x)\leq 1+\frac{b^2}{8}$. Thus
    $$
    1+\tfrac{b^2}{8} \geq \int f(x)L_b(x)\,dx=\sum_{j=-2}^2
    \hat{f}(j)\hat{L_b}(-j)=b\Re\hat{f}(1)-\Re\hat{f}(2).
    $$
Rearranging, we arrive at $\Re\hat{f}(2) \geq b (\Re\hat{f}(1))-1-\frac{b^2}{8}$.
Setting $b=4\Re\hat f(1)$, we find that $\Re\hat{f}(2) \geq 2(\Re\hat{f}(1))^2-1$.
\end{proof}

The first proof may be adapted to also give $4(\Re \hat{f}(1))^3-3\Re \hat{f}(1) \leq
\Re \hat{f}(3)$. The proof does not immediately extend to higher coefficients. The
second proof can be strengthened with the additional hypothesis that $f$ be an nif. We
take advantage of this in Proposition~\ref{Delta.one.half.thm}.

From the inequality $\Re \hat{f}(2) \leq 2\Re\hat{f}(1)-1$
(Lemma~\ref{Fhat.F1.F2.Domain.lem}) one easily computes that $\max\{|\hat{f}(1)|,
|\hat{f}(2)|\} \geq \frac13$, and with Lemma~\ref{Fhat.Green.Bound.lem} this gives
    $$\frac19\leq \frac{\ffi}{\pi}\sin\left(\frac{\pi}{\ffi}\right).$$
This yields $\ffi \geq 1.11$, a non-trivial bound. If one assumes that $f$ is an nif
supported on a subset of $[-\frac14,\frac14]$ with large measure, then one can do much
better than Lemma~\ref{Fhat.F1.F2.Domain.lem}. The following proposition establishes the
lower bound on $\De$ given in Theorem \ref{Delta.Summary.thm}(iii).

\begin{prop}\label{Delta.one.half.thm}
Let $f$ be an nif supported on a subset of $[-\frac14,\frac14]$ with measure $\e/2$.
Then
    $$\ffi \geq 1.1092 + 0.176158\,\e $$
and consequently
    $$ \De \geq 0.5546\e^2 + 0.088079\e^3.$$
\end{prop}

\begin{proof}
For $\e\ge \frac58$, this proposition is weaker than
Lemma~\ref{Delta.Trivial.Bounds.lem}, and for $\e\le \frac38$ it is weaker than
Proposition~\ref{the.full.De.bound.prop}, so we restrict our attention to
$\frac38<\e<\frac58$.

Let $b>-1$ be a parameter and set $L_b(x):=\cos(4\pi x)-b \cos(2\pi x)$. If we define
$F:=\max\{\Re\hat{f}(1),-\Re\hat{f}(2)\}$, then
    $$
    \int f(x) L_b(x)\,dx = \Re\hat{f}(2)-b\Re\hat{f}(1) \geq -(b+1)F
    $$
on the one hand, and
    $$
    \int f(x) L_b(x)\,dx \leq \int \sdr{f}(x)\sdr{L_b}(x)\,dx
        =\int_{-\e/4}^{\e/4} \tfrac 2\e \,\sdr{L_b}(x)\,dx
    $$
on the other, where $\sdr{L_b}(x)$ is the symmetric decreasing rearrangement of $K_b(x)$
on the interval $[-\frac14,\frac14]$. Thus
    $$
    F \geq \frac{-1}{b+1} \frac 2\e \int_{-\e/4}^{\e/4} \sdr{L_b}(x)\,dx.
    $$
The right-hand side may be computed explicitly as a function of $\e$ and $b$ and then
the value of $b$ chosen in terms of $\e$ to maximize the resulting expression. One finds
that for $\e<\frac58$, the optimal choice of $b$ lies in the interval $2<b<4$, and the
resulting lower bound for $F$ is
    $$
    F \geq \frac{3\cos (\frac{\pi \e }{4}) +
    \sin (\frac{\pi \e }{4}) -
    {\sqrt{3 + 4\cos (\frac{\pi \e }{2}) +
        2\cos (\pi \e ) -
        \sin (\frac{\pi \e }{2})}}}{\pi
     \e \cos (\frac{\pi \e }{4}) +
    \pi \e \sin (\frac{\pi \e }{4})}.
    $$
From Lemma~\ref{Fhat.Green.Bound.lem} we know that $F^2\leq
\frac{\ffi}{\pi}\sin\left(\frac{\pi}{\ffi}\right) $. We compare these bounds on $F$ to
conclude the proof. Specifically,
    \begin{equation}\label{F2.lower.bound}
    F^2 \leq \frac{\ffi}{\pi}\sin\left(\frac{\pi}{\ffi}\right)
        \leq \frac{3}{5\pi } + \frac{\left( 6 + 5{\sqrt{3}}\pi  \right)
     \left( \ffi-\frac65 \right) }{12\pi },
    \end{equation}
where the expression on the right-hand side of this equation is from the Taylor
expansion of $\frac x\pi\sin(\frac\pi x)$ at $x_0=\frac65$, and
    \begin{align}
    F^2 &\geq \left(\frac{3\cos (\frac{\pi \e }{4}) +
         \sin (\frac{\pi \e }{4}) -
        {\sqrt{3 + 4\cos (\frac{\pi \e }{2}) +
        2\cos (\pi \e ) -
        \sin (\frac{\pi \e }{2})}}}{\pi
         \e \cos (\frac{\pi \e }{4}) +
        \pi \e \sin (\frac{\pi \e }{4})}\right)^2
            \notag \\
     &\geq \frac{-8\left( -3 - \sqrt{2} + \sqrt{3} + \sqrt{6} \right) }{ \pi^2}
            \notag \\
     &\qquad + \frac{\left( 96\left( -3 - {\sqrt{2}} + {\sqrt{3}} + {\sqrt{6}} \right)  -
       4\left( 9{\sqrt{2}} - 10{\sqrt{3}} + {\sqrt{6}} \right) \pi  \right)
     \left( \e-\frac12 \right) }{3{\pi }^2}, \label{F2.upper.bound}
     \end{align}
where the expression on the right-hand side is from the Taylor expansion of the middle
expression at $\e_0=\frac 12$. Comparing Eqs.~(\ref{F2.lower.bound})
and~(\ref{F2.upper.bound}) gives a lower bound on $\ffi$, say $\ffi\geq c_1 + c_2 \e$
with certain constants $c_1,c_2$. It is easily checked that $c_1> 1.1092 $ and
$c_2>0.176158$, concluding the proof of the first asserted inequality. The second
inequality then follows from Lemma \ref{Fourier.Connection.lem}.
\end{proof}

\pagebreak \section{Upper Bounds for $R(g,n)$ Arising from $\De$}
\label{upper.bounds.Rgn.from.De.section}

\subsection{Inequalities Relating $\De$ and $R(g,n)$}

A symmetric set consists of pairs $(x,y)$ all with a fixed midpoint
$c=\tfrac{x+y}{2}$. If there are few pairs in $E\times E$ with a given sum
$2c$, then there will be no large symmetric subset of $E$ with center $c$. We
take advantage of the constructions of large integer sets whose pairwise sums
repeat at most $g$ times to construct large real subsets of $[0,1)$ with no
large symmetric subsets. Recall the definition \eqref{Rgn.definition} of the function $R(g,n)$.

\begin{prop}
For any integers $n\ge g\ge1$, we have $\Delta(\frac{R(g,n)}n) \leq
\frac{g}n.$ \label{Delta.from.R.Connection.prop}
\end{prop}

\begin{proof}
Let $S\subseteq\{1,2,\dots,n\}$ be a $B^\ast[g]$ set with $|S|=R(g,n)$. Let $$A(S) :=
\bigcup_{s\in S} \big[ \frac{s-1}{n},\frac{s}{n} \big)$$ as in Eq.~(\ref{ASdefinition});
it suffices to show that the largest symmetric subset of $A(S)$ has measure at most
$\frac gn$. Notice that the set $A(S)$ is a finite union of intervals, and so the
function $\lambda\big(A(S)\cap(2c-A(S))\big)$, which gives the measure of the largest
symmetric subset of $A(S)$ with center $c$, is piecewise linear. (Figure
\ref{AS.and.convolution} contains a typical example of the set $A(S)$ portrayed in dark
gray below the $c$-axis, together with the function $\lambda\big(A(S)\cap(2c-A(S))\big)$
shown as the upper boundary of the light gray region above the $c$-axis, for
$S=\{1,2,3,5,8,13\}$.) Without loss of generality, therefore, we may restrict our
attention to those symmetric subsets of $A(S)$ whose center $c$ is the midpoint of
endpoints of any two intervals $\big( \frac{s-1}{n},\frac{s}{n} \big)$. In other words,
we may assume that $2nc\in\Z$.

\begin{figure}[h]
    \label{AS.and.convolution}
\begin{center}
    \begin{picture}(384,178)
        \put(0,166){$\lambda(A(S)\cap(2c-A(S)))$}
        \put(370,22){$c$}
        \put(30,60){$\tfrac{1}{13}$}
        \put(30,100){$\tfrac{2}{13}$}
        \put(30,140){$\tfrac{3}{13}$}
        \ifpdf
            \put(10,-62){\includegraphics{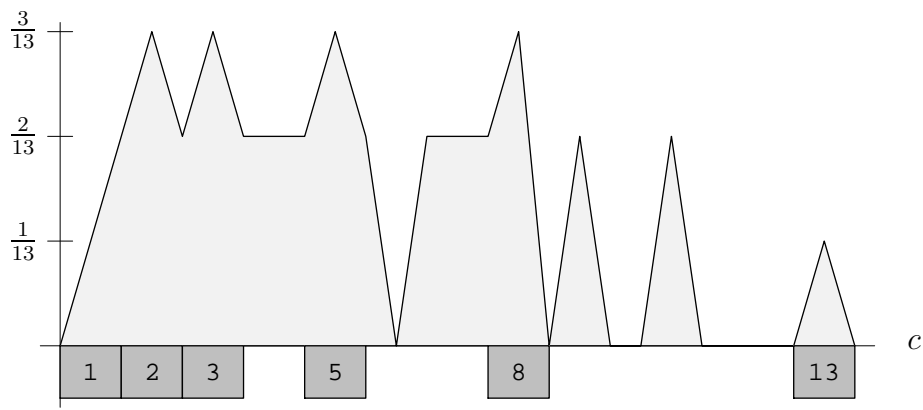}}
        \else
            \put(10,0){\includegraphics{FatSet}}
        \fi
    \end{picture}
\end{center}
    \caption{$A(S)$, and the function $\lambda(A(S)\cap(2c-A(S)))$,
        with $S=\{1,2,3,5,8,13\}$}
\end{figure}

Suppose $u$ and $v$ are elements of $A(S)$ such that $\frac{u+v}2=c$. Write
$u=\frac{s_1}n - \frac1{2n} + x$ and $v=\frac{s_2}n - \frac1{2n} + y$ for integers
$s_1,s_2\in S$ and real numbers $x,y$ satisfying $|x|,|y|<\frac1{2n}$. (We may ignore
the possibility that $nu$ or $nv$ is an integer, since this is a measure-zero event for
any fixed $c$.) Then $2nc = n(u+v) = s_1+s_2-1+n(x+y)$, and since $2nc$, $s_1$, and
$s_2$ are all integers, we see that $n(x+y)$ is also an integer. But $|n(x+y)| < 1$, so
$x+y=0$ and $s_1+s_2 = 2nc+1$.

Since $S$ is a $\Bg$ set, there are at most $g$ solutions $(s_1,s_2)$ to the equation
$s_1+s_2 = 2nc+1$. If it happens that $s_1=s_2$, the interval $\big(
\frac{s_1-1}{n},\frac{s_1}{n} \big)$ (a set of measure $\frac1n$) is contributed to the
symmetric subset with center $c$. Otherwise, the set $\big(
\frac{s_1-1}{n},\frac{s_1}{n} \big) \cup \big( \frac{s_2-1}{n},\frac{s_2}{n} \big)$ (a
set of measure $\frac2n$) is contributed to the symmetric subset with center $c$, but
this counts for the two solutions $(s_1,s_2)$ and $(s_2,s_1)$. In total, then, the
largest symmetric subset having center $c$ has measure at most $\frac gn$. This
establishes the theorem.
\end{proof}

Using Proposition~\ref{Delta.from.R.Connection.prop}, we can translate lower bounds on
$\De$ into upper bounds on $R(g,n)$, as in Corollary~\ref{Delta.R.Basic.cor}.

\begin{cor}\label{Delta.R.Basic.cor}
If $\delta\le \inf_{0<\e<1} \De/\e^2$, then $R(g,n) \le \delta^{-1/2} \sqrt{gn}$ for all
$n\ge g \ge 1$.
\end{cor}

We remark that we may take $\delta = \Dconstant$ by Proposition~\ref{Main.Bound.prop},
and so this corollary implies that $R(g,n) \le \Rconstant\sqrt{gn}$, which is one of the assertions of Theorem \ref{R.Upper.Bound.thm}.

\begin{proof}
Combining the hypothesized lower bound $\Delta(\e)\ge \delta \e^2$ with
Proposition~\ref{Delta.from.R.Connection.prop}, we find that
    $$
    \delta\bigg( \frac{R(g,n)}n \bigg)^2 \le \Delta\bigg( \frac{R(g,n)}n \bigg)
    \le \frac gn
    $$
which is equivalent to $R(g,n)\le \delta^{-1/2}\sqrt{gn}$.
\end{proof}

We have been unable to prove or disprove that
 $$\lim_{g\to\infty} \lim_{n\to\infty} \frac{R(g,n)}{\sqrt{gn}} = \left(\inf_{0<\e<1}
 \frac{\De}{\e^2} \right)^{-1/2},$$
i.e., that Corollary~\ref{Delta.R.Basic.cor} is best possible as $g\to\infty$. At any rate, for small $g$ it is possible to
do better by taking advantage of the shape of the set $A(S)$ used in the proof of
Proposition~\ref{Delta.from.R.Connection.prop}. This is the subject of
Section~\ref{upper.bounds.Rgn.section}.

Proposition~\ref{Delta.from.R.Connection.prop} provides a one-sided inequality linking
$\Delta(\e)$ and $R(g,n)$. It will also be useful for us to prove a theoretical result
showing that the problems of determining the asymptotics of the two functions are, in a
weak sense, equivalent. In particular, the following proposition implies
that the trivial lower bound $\Delta(\e)\ge\frac12\e^2$ and the trivial upper bound
$R(g,n)\le\sqrt{2gn}$ are actually equivalent. Further, any nontrivial lower bound on
$\De$ gives a nontrivial upper bound on $R(g,n)$, and vice versa.

\begin{prop}
$\Delta(\e) = \inf\{\frac gn\colon n\ge g\ge1,\, \frac{R(g,n)}n\ge\e\}$ for all
$0\le\e\le1$.
\label{Delta.as.infimum.prop}
\end{prop}

\begin{proof}
That $\Delta(\e)$ is bounded above by the right-hand side follows immediately from
Proposition~\ref{Delta.from.R.Connection.prop} and the fact that $\Delta$ is an
increasing function. For the complementary inequality, let $S\subseteq[0,1)$ with
$\lambda(S)=\e$. There exists a finite union $T$ of open intervals such that
$\lambda(S\diamond T) < \eta$, and it is easily seen that $T$ can be chosen to meet the
following criteria: $T\subseteq[0,1)$, the endpoints of the finitely many intervals
comprising $T$ are rational, and $\lambda(T)>\e$. Choosing a common denominator $n$ for
the endpoints of the intervals comprising $T$, we may write $T = \bigcup_{m\in M} \big[
\frac {m-1}n, \frac mn \big)$ (up to a finite set of points) for some set of integers
$M\subseteq \{1,\dots,n\}$; most likely we have greatly increased the number of
intervals comprising $T$ by writing it in this manner, and $M$ contains many consecutive
integers. Let $g$ be the maximal number of solutions $(m_1,m_2)\in M\times M$ to
$m_1+m_2=k$ as $k$ varies over all integers, so that $M$ is a $B^\ast[g]$ set and thus
$|M| \le R(g,n)$ by the definition of $R$. It follows that $\e < \lambda(T) = |M|/n \le
R(g,n)/n$. Now $T$ is exactly the set $A(M)$ as defined in Eq.~(\ref{ASdefinition});
hence $D(T) = \frac gn$ as we saw in the proof of
Proposition~\ref{Delta.from.R.Connection.prop}. Therefore by Lemma~\ref{Diamond.lem},
\begin{equation*}
    D(S) \ge D(T) - 2\eta   = \tfrac gn - 2\eta  \ge
    \inf\{\tfrac gn\colon n\ge g \ge1,\,
                        \tfrac{R(g,n)}n\ge\e\} - 2\eta.
\end{equation*}
Taking the infimum over appropriate sets $S$ and noting that $\eta>0$ was
arbitrary, we derive the desired inequality $\Delta(\e) \ge \inf\{\frac
gn\colon n\ge g\ge1,\, \frac{R(g,n)}n\ge\e\}$.
\end{proof}

\subsection{Upper Bounds on $R(g,n)$ and $\overline{\rho}(g)$}
\label{upper.bounds.Rgn.section}

The bulk of Theorem~\ref{R.Upper.Bound.thm} follows immediately from
Proposition~\ref{R.from.Delta.Connection.cor}, which is proved in this section.

Let $S\subseteq\{1,2,\dots,n\}$ be a \Bg set. We call the function
 $$f(x) :=
    \left\{%
    \begin{array}{ll}
        \frac{2n}{|S|},
        & \hbox{$\frac{s-1}{2n}-\frac14 \leq x < \frac{s}{2n}-\frac14$ for some $s\in S$}, \\
        0, & \hbox{otherwise} \\
    \end{array}%
    \right.$$
the nif corresponding to $S$. We think of $S$, and consequently $f$, as depending on $n$
(many readers may prefer that we write $S_n$ and $f_n$, but this is neither customary
nor sufficiently brief). For example, when we write ``Let $S\subseteq\{1,2,\dots,n\}$ be
a \Bg set with $|S|>0.7\sqrt{gn}$'', we mean ``For each $n\in\Z^+$, let
$S=S_n\subseteq\{1,2,\dots,n\}$ be a \Bg set with $|S|=|S_n|>0.7\sqrt{gn}$.''

Note that $f$ is piecewise constant, so that $\ff$ is piecewise linear. Furthermore, the corners of $\ff$ are only at the points $\frac1{4n}\Z$, and if we define
\begin{equation}
r(t):=\#\left\{(x,y)\in S^2 \colon x+y=t\right\},
\label{number.of.representations.def}
\end{equation}
then $\ff(\frac{k}{4n})=\frac{2nr(k+n+1)}{|S|^2}$. In particular, $S$ is a \Bg set if and only if $\ffi \leq \frac{2gn}{|S|^2}$.

\begin{prop}\label{R.from.Delta.Connection.cor}
Let $\phi:=\fstarftwonormconstant$ be the constant appearing in Corollary \ref{overall.twonorm.bound}, and let $\theta_0:=21.922911$, $\theta_1:=-33.711941$, and $\theta_2:=13.676987$ be the constants appearing in Lemma \ref{fft.bound.using.ffi}. Then:
\renewcommand{\theenumi}{\roman{enumi}}
\begin{enumerate}
    \item $R(g,n) \le \sqrt{2/\phi} \sqrt{(g-1)n}+\frac13$ for $n\ge g \ge 2$ and $g$ odd;
    \label{explicit.R.bound.item}
    \item $\overline{\rho}(g) \le \sqrt{2/\phi} \left(1-\frac1{g\sqrt{3}}\right)^{1/2}$ for
            $g\ge 2$, and if moreover $\overline{\rho}(g) \geq 1.275237$ then
    \begin{align*}
    \overline{\rho}(g)
    &\leq
        \frac{2\sqrt{3\theta_2g}}
        {\left({3(1-\theta_1)g-\sqrt{3}-\sqrt{\left(3(1-\theta_1)g-\sqrt{3} \right)^2
            -36\theta_0\theta_2g^2 }}\right)^{1/2}};
    \end{align*}
    \label{even.rho.item}
    \item   \label{odd.rho.item}
            If $g\geq 2$ is odd, then $\overline{\rho}(g) \le \sqrt{2/\phi}
            \left(1-\frac{1+1/\sqrt{3}}{g}\right)^{1/2}$, and if moreover $\overline{\rho}(g) \geq
            1.275237$, then
    \begin{align*}
    \overline{\rho}(g) &\le \frac{2\sqrt{3\theta_2g}}
   {{\left( 3\left( 1 - {\theta_1} \right) g - {\sqrt{3}} - 3 -
        {\sqrt{{\left( 3\left( 1 - {\theta_1} \right) g - {\sqrt{3}} - 3 \right) }^2 - 36{\theta_0}{\theta_2}g^2}}
        \right) }^{1/2}}.
   \end{align*}
\end{enumerate}
\end{prop}

We begin with a simple lemma.

\begin{lem}\label{Removal.of.Spikes.lem}
If $h,p,q$ are nonnegative functions with $h=p+q$ and $\|h\|_\infty\ge
\|p\|_\infty+\|q\|_\infty$, then
 $$\|h\|_\infty \geq \frac{\|h\|_2^2-\|p\|_2^2}{\|h\|_1+\|p\|_1}+\|p\|_\infty.$$
\end{lem}

\begin{proof}
We have
\begin{align*}
    \|h\|_2^2
    &=      \|p+q\|_2^2 \\
    &=      \|p\|_2^2+\|(2p+q)q\|_1 \\
    &\leq   \|p\|_2^2+\|2p+q\|_1 \|q\|_\infty \\
    &=      \|p\|_2^2+(\|h\|_1+\|p\|_1)\|q\|_\infty \\
    &\leq   \|p\|_2^2+(\|h\|_1+\|p\|_1)(\|h\|_\infty-\|p\|_\infty).
\end{align*}
\end{proof}

We will use this lemma with $h=\ff$ ($f$ a pdf), and $p$ chosen so that $\|p\|_1\to0$,
$\|p\|_2\to 0$, $\|p\|_\infty \not\to 0$. In this case, we have the inequality
    $$\ffi \gtrsim \|\ff\|_2^2 + \|p\|_\infty$$
which is stronger than H{\"o}lder's Inequality: $\ffi \ge \|\ff\|_2^2$.

\begin{proofof}{Proposition~\ref{R.from.Delta.Connection.cor}(\ref{explicit.R.bound.item})}
The idea of the proof is that even though $\ff$ might take
values near $\frac{2gn}{|S|^2}$, it does so only on a set of small measure, and away from that small set it is bounded by $\frac{2(g-1)n}{|S|^2}$. If the pair $(s_1,s_2)$ contributes to $r(k)$, then so does $(s_2,s_1)$, and therefore
$r(k)$ is odd if and only if $k\in \{2s \colon s\in S\}$. There are only $|S|$ such integers $k$, and no two are consecutive. Consequently, if $\ff(\frac k{4n})=\frac{2gn}{|S|^2}$, then $\ff(\frac {k-1}{4n}) \leq \frac{2(g-1)n}{|S|^2}$ and $\ff(\frac {k+1}{4n}) \leq
\frac{2(g-1)n}{|S|^2}$.

We put this idea into effect by writing $\ff(x)=p(x)+q(x)$, where $p$ represents the small contribution to $r(k)$ of pairs $(s,s)$ and $q$ is the remaining majority of $\ff$. More precisely, let $p(x) := \sum_{s\in S} T(x-\frac{s-1}{n}-\frac12)$, where $T(x)$ is
the ``tent function''
    \begin{equation*}
    T(x) := \begin{cases}
    \frac{2n}{|S|^2} (1-2n|x|), &{\rm if\ }|x| \le {\frac1{2n}}, \\
    0, &\rm{otherwise,}
    \end{cases}
    \end{equation*}
and let $q(x) := f*f(x)-p(x)$. (See Figure~\ref{minding.ps.and.qs.fig} for an
illustration of a typical example, again using the $B^*[3]$ set $S=\{1,2,3,5,8,13\}$.)
Because of our judicious choice of the peaks of $p$, both $p$ and $q$ are nonnegative
and $\ffi \ge \|q\|_\infty + \|p\|_\infty$. We compute directly from the definition of
$p$ that $\|p\|_1 = |S|\|T\|_1 = \frac1{2|S|}$, $\|p\|_2^2 = |S|\|T\|_2^2 =
\frac{2n/3}{|S|^3}$, and $\|p\|_\infty= \|T\|_\infty = \frac{2n}{|S|^2}$.

\begin{figure}
    \centerline{    \begin{picture}(192,120)
                    \put(17,109){$f \ast f(x)$}
                    \put(171,15){$x$}
                    \put(18,40){$\tfrac{13}{18}$}
                    \put(18,68){$\tfrac{13}{9}$}
                    \put(18,94){$\tfrac{13}{6}$}
                    \put(27,4){$-\tfrac12$}
                    \put(97,4){0}
                    \put(161,4){$\tfrac12$}
                    \ifpdf
                        \put(12,-120){\includegraphics{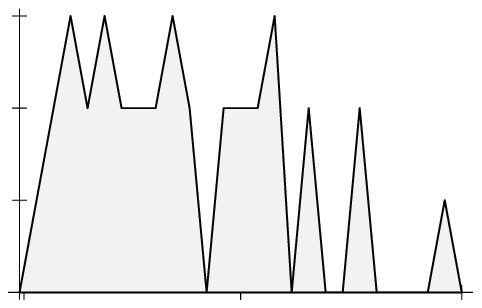}}
                    \else
                        \put(12,0){\includegraphics{fstarf}}
                    \fi
                    \end{picture}
                    }
    \centerline{    {
                    \begin{picture}(192,120)
                    \put(30,109){$p(x)$}
                    \put(174,15){$x$}
                    \put(18,40){$\tfrac{13}{18}$}
                    \put(18,68){$\tfrac{13}{9}$}
                    \put(18,95){$\tfrac{13}{6}$}
                    \put(28,4){$-\tfrac12$}
                    \put(101,4){0}
                    \put(166,4){$\tfrac12$}
                    \ifpdf
                        \put(9,-121){\includegraphics{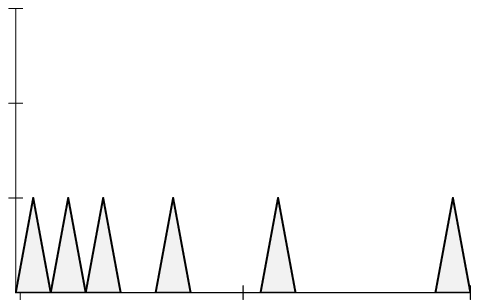}}
                    \else
                        \put(9,0){\includegraphics{fstarfppart}}
                    \fi
                    \end{picture}
                    }
                    \hfill
                    {
                    \begin{picture}(192,120)
                    \put(30,109){$q(x)$}
                    \put(174,15){$x$}
                    \put(18,40){$\tfrac{13}{18}$}
                    \put(18,68){$\tfrac{13}{9}$}
                    \put(18,95){$\tfrac{13}{6}$}
                    \put(28,4){$-\tfrac12$}
                    \put(101,4){0}
                    \put(166,4){$\tfrac12$}
                    \ifpdf
                        \put(9,-121){\includegraphics{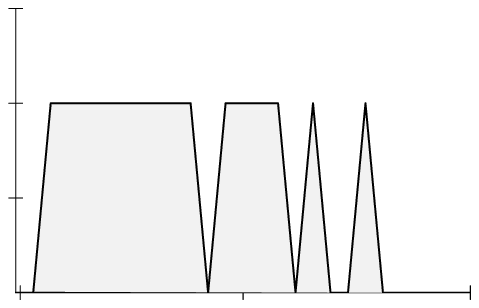}}
                    \else
                        \put(9,0){\includegraphics{fstarfqpart}}
                    \fi
                    \end{picture}
                    }
                }
    \caption{The decomposition $\ff(x)=p(x)+q(x)$, where $f$ is the nif corresponding
    to $S=\{1,2,3,5,8,13\}$}
    \label{minding.ps.and.qs.fig}
\end{figure}

Lemma~\ref{Removal.of.Spikes.lem} with $h=\ff$ gives
\begin{equation}\label{stronger.than.Holder.eq}
\ffi \geq \frac{\|\ff\|_2^2-\|p\|_2^2}{1+\|p\|_1} + \|p\|_\infty.
\end{equation}

Using the inequality $\|\ff\|_2^2 \geq \phi$ from Corollary~\ref{overall.twonorm.bound},
the inequality $\ffi \le \frac{2gn}{|S|^2}$, and the above computations of $\|p\|_1$,
$\|p\|_2$, and $\|p\|_\infty$, we can deduce from Eq.~(\ref{stronger.than.Holder.eq})
that
    $$|S|^2 \leq \frac{2(g-1)n}{\phi}+\frac{(g-1)n+2n/3}{\phi|S|}
            \le \frac{2(g-1)n}{\phi}+\frac{gn}{\phi|S|},$$
from which it follows that $|S| \leq \sqrt{\frac{2(g-1)n}{\phi}} + \frac13$.
\end{proofof}

\begin{lem}\label{periodic.coefficients.lem}
If $f$ is the nif corresponding to $S\subseteq \{1,2,\dots,n\}$, and $m^\prime \to
\infty$ with $m^\prime=o\left(\frac{n}{|S|}\right)$, then
 $$\sum_{|j|>m^\prime} |\hat f(j)|^4 \geq \frac{2}{\sqrt{3}}\frac{n}{|S|^2}-o(1).$$
\end{lem}

\begin{proof}
Note that for $j\not=0$,
    \begin{align*}
    j\hat{f}(j)
    &=  j\int f(x) e^{-2\pi i j x}\,dx \\
    &=  \sum_{s \in S} \|f\|_\infty \,j
        \int_{\frac{s-1}{2n}-\frac14}^{\frac{s}{2n}-\frac14} e^{-2\pi i j x}\,dx \\
    &=  \sum_{s \in S} \|f\|_\infty \frac{j}{-2\pi i j}
        \left( e^{-2\pi i j (\frac{s}{2n}-\frac14)}
                        -e^{-2\pi i j (\frac{s-1}{2n}-\frac14)}\right) \\
    &=  \frac{-\|f\|_\infty}{2\pi i} \sum_{s\in S}
        \left( e^{-2\pi i j (\frac{s}{2n}-\frac14)}
                        -e^{-2\pi i j (\frac{s-1}{2n}-\frac14)}\right).
    \end{align*}

In particular, if we set $c(j):=j\,|\hat{f}(j)|$, then $c(j)=c(j+2n)$ (provided neither
$j$ nor $j+2n$ is zero) and $c(j)=c(-j)$. Exploiting this periodicity is the heart of
this lemma.

We define
    \begin{align*}
    J   :=&\{m^{\prime}+1,m^{\prime}+2,\dots,n\}\cup\{2n-m^{\prime},2n-m^{\prime}+1,\dots,2n-1\}\\
    J^{\prime}:=&[m^{\prime}+1,m^{\prime}+2n]\setminus J,
    \end{align*}
and for $j\in J$ we set
 $$j^{\prime} = \left\{%
    \begin{array}{ll}
        2n-j, & m^{\prime}+1 \le j < n; \\
        2n, &  j=n; \\
        4n-j, & 2n-m^{\prime} \le j <2n. \\
    \end{array}%
    \right.$$
Observe that $c(j)=c(j^{\prime})$ and $J^{\prime}=\{j^{\prime}\colon j \in J\}$. Also,
for almost all $j\in J$ we have $j^{\prime}=2n-j$, since $m^{\prime}=o(n)$.

For all $p>1$,
    \begin{align}
    \sum_{|j|>m^{\prime}} |\hat{f}(j)|^p
    &=  2\sum_{j=m^{\prime}+1}^\infty |\hat{f}(j)|^p \notag \\
    &=  2\sum_{j=m^{\prime}+1}^\infty \left(\frac{c(j)}{j}\right)^p \notag \\
    &=  2\sum_{j=m^{\prime}+1}^{m^{\prime}+2n} \sum_{k=0}^\infty
                                        \left(\frac{c(j)}{j+2n\cdot k}\right)^p \notag \\
    &=  2\sum_{j=m^{\prime}+1}^{m^{\prime}+2n} \left(\frac{c(j)}{2n}\right)^p \zeta(p,\tfrac{j}{2n}) \notag \\
    &=  2 \sum_{j\in J} \left(\frac{c(j)}{2n}\right)^p
            \left(\zeta(p,\tfrac{j}{2n})+\zeta(p,\tfrac{j^{\prime}}{2n})\right),
\label{Hurwitz.zeta.appeared}
    \end{align}
where $\zeta(s,a):=\sum_{k=0}^\infty (k+a)^{-s}$ is the Hurwitz zeta function.

We shall use Cauchy's inequality in the form $\sum a_j^2 \geq \frac{(\sum a_j
b_j)^2}{\sum b_j^2}$, with
 $$a_j=\left(\frac{c(j)}{2n}\right)^2\left(\zeta(4,\tfrac{j}{2n}) +
 \zeta(4,\tfrac{j^{\prime}}{2n})\right)^{1/2}$$
and
    $$b_j = \frac{\zeta(2,\tfrac{j}{2n}) +
 \zeta(2,\tfrac{j^{\prime}}{2n})}{\left(\zeta(4,\tfrac{j}{2n}) +
 \zeta(4,\tfrac{j^{\prime}}{2n})\right)^{1/2}}.$$
For notational convenience set
    $$Z(m^{\prime},n):=\frac{1}{2n}\sum_{j\in J} b_j^2 = \frac1{2n}{\sum_{j\in J}
    \frac{\left(\zeta(2,\frac{j}{2n}) + \zeta(2,\frac{j^{\prime}}{2n})\right)^2
    }{\zeta(4,\frac{j}{2n}) + \zeta(4,\frac{j^{\prime}}{2n})}},$$
and note that since $m^{\prime}=o(n)$,
\begin{align}
    \lim_{n\to\infty} Z(m^{\prime},n) &= \lim_{n\to\infty} \frac1{2n}{\sum_{j=o(n)}^n
    \frac{\left(\zeta(2,\frac{j}{2n}) + \zeta(2,\frac{2n-j}{2n})\right)^2
    }{\zeta(4,\frac{j}{2n}) + \zeta(4,\frac{2n-j}{2n})}} \notag \\
&= \int_0^{1/2}
    \frac{(\zeta(2,a)+\zeta(2,1-a))^2}{\zeta(4,a)+\zeta(4,1-a)}\,da
    =\int_0^{1/2} \frac{3}{2 + \cos (2\pi a )}\,da
    =\frac{\sqrt{3}}{2}.
\label{zeta.integral}
\end{align}

Therefore from Eq. \eqref{Hurwitz.zeta.appeared} with $q=4$ and Cauchy's inequality,
    \begin{align}
    \sum_{|j|>m^{\prime}} |\hat{f}(j)|^4
    &=      2 \sum_{j\in J} \left(\frac{c(j)}{2n}\right)^4
                \left(\zeta(4,\frac{j}{2n})+\zeta(4,\frac{j^{\prime}}{2n})\right) \notag \\
    &\ge    2 \frac{\left(\sum_{j\in J} \left(\frac{c(j)}{2n}\right)^2
                \left(\zeta(2,\frac{j}{2n})
                    +\zeta(2,\frac{j^{\prime}}{2n})\right) \right)^2}
            {\sum_{j\in J} {\left(\zeta(2,\frac{j}{2n}) +
                    \zeta(2,\frac{j^{\prime}}{2n})\right)^2}
                    \big( {\zeta(4,\frac{j}{2n}) + \zeta(4,\frac{j^{\prime}}{2n})} \big)^{-1}} \notag \\
    &=      2 \frac{\left(\frac12 \sum_{|j|>m^{\prime}} |\hat{f}(j)|^2 \right)^2}{2n \,Z(m^{\prime},n)} \notag \\
    &=       \frac{ \left( \|f\|_2^2 - \sum_{|j|\le m^{\prime}} |\hat{f}(j)|^2
                     \right)^2}{4n\, Z(m^{\prime},n)},
\label{asdfgh}
    \end{align}
the last equality following from Parseval's formula. Trivially, $|\hat{f}(j)|\leq\hat
f(0)=1$, and since $f$ is the nif corresponding to $S$, we have $\|f\|_2^2 =
\|f\|_\infty^2 \lambda(\supp(f))=\|f\|_\infty=\frac{2n}{|S|}$. Since $m'=o(n/|S|)$,
Eq.~\eqref{asdfgh} becomes
    \begin{equation*}
    \sum_{|j|>m^{\prime}} |\hat{f}(j)|^4
    \ge \frac{1}{4n} \frac{\left(\frac{2n}{|S|}-2m'-1\right)^2}{Z(m^{\prime},n)}
    = \frac{2}{\sqrt{3}}\frac{n}{|S|^2}-o(1).
    \end{equation*}
\end{proof}

\begin{proofof}{Proposition~\ref{R.from.Delta.Connection.cor}(\ref{even.rho.item})}
Let $S\subseteq\{1,2,\dots,n\}$ be a \Bg set with $|S|=R(g,n)$ and let $f$ be the
corresponding nif. We shall use the inequality
    $$\ffi \geq \|\ff\|_2^2 =\sum_{j} |\hat{f}(j)|^4=
    \sum_{|j|\le m^{\prime}} |\hat{f}(j)|^4 + \sum_{|j|>m^{\prime}} |\hat{f}(j)|^4.$$
We use Proposition~\ref{overall.twonorm.bound} and Proposition~\ref{Main.Bound.prop} to
bound the sum over small $|j|$, and we use Lemma~\ref{periodic.coefficients.lem} to
bound the sum over large $|j|$.

We are now prepared to prove
Proposition~\ref{R.from.Delta.Connection.cor}(\ref{even.rho.item}). We have
    \begin{multline*}
    \frac{2}{\overline{\rho}^2(g)}
    =  \liminf_{n\to\infty} \frac{2gn}{|S|^2}
    \ge  \liminf_{n\to\infty} \ffi
    \ge    \liminf_{n\to\infty} \|f\ast f\|_2^2 \\
    =  \liminf_{n\to\infty} \sum_j |\hat{f}(j)|^4
    =  \liminf_{n\to\infty} \bigg( \sum_{|j|\le m^{\prime}}|\hat{f}(j)|^4 + \sum_{|j|>
            m^{\prime}}|\hat{f}(j)|^4 \bigg).
    \end{multline*}
We see from Corollary~\ref{overall.twonorm.bound} and
Lemma~\ref{periodic.coefficients.lem} that
    \begin{equation*}
    \frac{2}{\overline{\rho}^2(g)}
    \ge \liminf_{n\to\infty}
            \Big(\phi-o(1)+
            \frac{2}{\sqrt{3}}\frac{n}{|S|^2}-o(1)\Big)\\
    =  \phi+ \frac{2}{\sqrt{3}}\frac{1}{\overline{\rho}^2(g) g}
    \end{equation*}
Solving this inequality for $\overline{\rho}(g)$ yields
    $$
        \overline{\rho}(g)\le \sqrt{\frac 2\phi} \left(1-\frac1{g\sqrt{3}}\right)^{1/2}.
    $$

If $\overline{\rho}(g) \ge 1.275237$, then $\ffi \le 1.229837$ for infinitely many $n$,
and we can use Lemma~\ref{fft.bound.using.ffi} instead of
Corollary~\ref{overall.twonorm.bound}. Setting
    $$F:=\liminf_{n\to\infty} \ffi \leq \liminf_{n\to\infty} \frac{2gn}{|S|^2},$$
we have
    \begin{align*}
    F
    &\ge  \bigg( \liminf_{n\to\infty} \sum_{|j|\le m^{\prime}}|\hat{f}(j)|^4 + \sum_{|j|>
            m^{\prime}}|\hat{f}(j)|^4 \bigg) \\
    &\ge  \liminf_{n\to\infty}
                \Big(\theta_0+\theta_1\ffi +\theta_2\ffi^2
                -o(1) + \frac{2}{\sqrt{3}}\frac{n}{|S|^2}-o(1) \Big)\\
    &\ge    \theta_0+\theta_1 F+\theta_2 F^2 + \frac{F}{g\sqrt{3}},
    \end{align*}
whence
    $$
    \overline{\rho}(g) \le \sqrt{\frac{2}{F}} \le \frac{2{\sqrt{3\theta_2g}}}
      {\left({3(1-\theta_1)g-\sqrt{3}-\sqrt{\left(3(1-\theta_1)g-\sqrt{3} \right)^2
      -36\theta_0\theta_2g^2 }}\right)^{1/2}}.
    $$
\end{proofof}

\begin{proofof}{Proposition~\ref{R.from.Delta.Connection.cor}(\ref{odd.rho.item})}
We combine the ideas of parts (\ref{explicit.R.bound.item}) and
(\ref{even.rho.item}). As in the proof of
Proposition~\ref{R.from.Delta.Connection.cor}(\ref{even.rho.item}), we let
$S\subseteq\{1,2,\dots,n\}$ be a \Bg set with $|S|=R(g,n)$ and $f$ be the
corresponding nif.

If $g$ is odd, then (defining $p$ and $q$ as in the proof of
Proposition~\ref{R.from.Delta.Connection.cor}(\ref{explicit.R.bound.item}))
Eq.~\eqref{stronger.than.Holder.eq} is valid, and $\|p\|_1\to0,\|p\|_2 \to 0$, and
$\|p\|_\infty =\frac{2n}{|S|^2}\to \frac{2}{\overline{\rho}^2(g)\,g}$.

We find
    \begin{align*}
    \frac{2}{\overline{\rho}^2(g)}
    &=  \lim_{n\to\infty} \frac{2gn}{|S|^2}\\
    &=  \lim_{n\to\infty} \ffi \\
    &\ge    \lim_{n\to\infty} \frac{\|f \ast f\|_2^2 -\|p\|_2^2}{1+\|p\|_1}
                        +\|p\|_\infty\\
    &=  \lim_{n\to\infty} \frac{\left( \sum_{|j|\le m^{\prime}}|\hat{f}(j)|^4
            + \sum_{|j|>m^{\prime}}|\hat{f}(j)|^4\right) -\|p\|_2^2}{1+\|p\|_1}
                        +\|p\|_\infty\\
    &=  \lim_{n\to\infty} \frac{\left( \phi-o(1) + \frac{2}{\sqrt{3}}\frac{n}{|S|^2}
            -o(1)\right) -\|p\|_2^2}{1+\|p\|_1} +\|p\|_\infty\\
    &=  \phi + \frac{2}{\sqrt{3}} \frac{1}{\overline{\rho}^2(g) g} +
             \frac{2}{\overline{\rho}^2(g)g}\\
    &=  \phi + 2\frac{1+1/\sqrt{3}}{\overline{\rho}^2(g)g}.
    \end{align*}
Solving this inequality for $\overline{\rho}^2(g)$ yields
    $$\overline{\rho}(g) \le \sqrt{\frac 2\phi}
            \left(1-\frac{1+1/\sqrt{3}}{g}\right)^{1/2}.$$

If $\overline{\rho}(g) \ge 1.275237$, then $\ffi \le 1.229837$ for infinitely many $n$,
and we can use Lemma~\ref{fft.bound.using.ffi}. We have:
    \begin{align*}
    \|p\|_1 &= \frac{1}{2|S|} \to 0\\
    \|p\|_2^2 &= \frac{2n/3}{|S|^3} \to 0\\
    \|p\|_\infty &= \frac{2n}{|S|^2} =\frac 1g \ffi\\
    \sum_{|j|\le m^{\prime}}|\hat{f}(j)|^4
            & \ge \theta_0 +\theta_1 \ffi
                    + \theta_2 \ffi^2 -o(1)\\
    \sum_{|j|>m^{\prime}}|\hat{f}(j)|^4
            &= \frac{2}{\sqrt{3}}\frac{n}{|S|^2}-o(1)
                \geq \frac{\ffi}{g\sqrt{3}}-o(1).
    \end{align*}
Thus, again writing $F=\liminf_{n\to\infty}\ffi$, we now know
    \begin{align*}
    F&:= \liminf_{n\to\infty} \ffi \\
    &\ge  \liminf_{n\to\infty} \frac{\left( \sum_{|j|\le m^{\prime}}|\hat{f}(j)|^4
    + \sum_{|j|> m^{\prime}}|\hat{f}(j)|^4\right)-\|p\|_2^2}{1+\|p\|_1} +\|p\|_\infty\\
    &\ge \theta_0 +\theta_1 F+\theta_2 F^2+\frac{F}{g\sqrt{3}}+\frac{F}{g}.
    \end{align*}
Isolating $F$, we obtain
    $$F\ge \frac{\theta_1 + \frac{1 + {\sqrt{3}}}{{g\sqrt{3}}} - 1 +
     {\sqrt{{\left( \theta_1 + \frac{1 + {\sqrt{3}}}{{g\sqrt{3}}} - 1 \right) }^
          2 - 4\theta_0\theta_2}}}{-2\theta_2},
     $$
and so
    $$\overline{\rho}(g) \le \sqrt{ \frac 2F}
    \le \frac{2{\sqrt{3\theta_2g}}}
   {{\left( 3\left( 1 - {\theta_1} \right) g - {\sqrt{3}} - 3 -
        {\sqrt{{\left( 3\left( 1 - {\theta_1} \right) g - {\sqrt{3}} - 3 \right) }^2 - 36{\theta_0}{\theta_2}g^2}}
        \right) }^{1/2}},$$
which concludes the proof.
\end{proofof}

\subsection{Ubiquity of Repeated Sums in $B^*[g]$ Sets}
\label{ubiquity.section}

The method of proof of
Proposition~\ref{R.from.Delta.Connection.cor}(\ref{explicit.R.bound.item}) can be
adapted to yield more information about the number of representations of integers as
sums of pairs of elements from $\Bg$ sets. The following theorem gives a quantitative
statement of the fact that, if $S$ is a dense enough $\Bg$ set, then there is a
substantial number of integers $t$ such that $r(t)$ is large, where $r(t)$ is defined in
Eq.~\eqref{number.of.representations.def}.

\begin{thm}
Let $S\subseteq\{1,\dots,n\}$ be a $B^*[g]$ set, and let $0<L<g$ be a real number. The
number of integers $t$ such that $r(t)>L$ is at least
\begin{equation*}
\frac{\phi |S|^4 - 2Ln|S|^2}{n(g-L)(g+2L)},
\end{equation*}
where $\phi:=\fstarftwonormconstant$ is the constant appearing in Corollary
\ref{overall.twonorm.bound}. \label{ubiquity.thm}
\end{thm}

We give an illustration of this theorem after the proof.

\begin{proof}
Let $f$ be the nif corresponding to $S$. We note again that $\|f*f\|_1 = \|f\|_1^2 = 1$
and $\|f*f\|_2^2 \ge \phi$ by
Corollary~\ref{overall.twonorm.bound}.

Now let $Q$ be the number of integers $t$ between 1 and $2n$ such that $r(t)>L$. Our aim
is to show that $Q$ is large if $|S|$ is large enough and $L$ is small enough. Define
the functions $p(x)$ and $q(x)$ by
\begin{equation*}
q(x) = \min\left\{f*f(x),\frac{2Ln}{|S|^2}\right\}, \quad p(x) = f*f(x)-q(x).
\end{equation*}
We have
\begin{align}
    \phi \leq \|\ff\|_2^2=\|p+q\|_2^2
        &=  \|p\|_2^2+\|(2p+q)q\|_1 \notag \\
        &\leq   \|p\|_2^2+\|2p+q\|_1 \|q\|_\infty \notag \\
        &=      \|p\|_2^2+(1+\|p\|_1) \|q\|_\infty,\label{spikes}
\end{align}
where the latter inequality is by H{\"o}lder. By its definition, we have $\|q\|_\infty\leq
\frac{2Ln}{|S|^2}$. To bound $\|p\|_1$ and $\|p\|_2$, we look more closely at the
function $p$.

Suppose $r(t+1),r(t+2),\dots,r(t+K)$ are all greater than $L$, but that $r(t)$ and
$r(t+K+1)$ are at most $L$. Then the graph of $f*f$ between $\frac{t-n-1}{4n}$ and
$\frac{t+K-n}{4n}$ lies under the trapezoid with vertices $\big(
\frac{t-n-1}{4n},\frac{2Ln}{|S|^2} \big)$, $\big( \frac{t-n}{4n},\frac{2gn}{|S|^2}
\big)$, $\big( \frac{t+K-n-1}{4n},\frac{2gn}{|S|^2} \big)$, and $\big(
\frac{t+K-n}{4n},\frac{2Ln}{|S|^2} \big)$. The area of this trapezoid, which bounds the
contribution to $\|p\|_1$ from $x$ in the interval $\big[
\frac{t-n-1}{4n},\frac{t+K-n}{4n} \big]$, equals $\frac{K(g-L)}{2|S|^2}$. Since the
numbers $K$ from all such intervals sum to $Q$, we have established the bound $\|p\|_1
\le \frac{Q(g-L)}{2|S|^2}$.

Similarly, the graph of $(f*f)^2$ between $\frac{t-n-1}{4n}$ and $\frac{t+K-n}{4n}$ lies
under the trapezoid with vertices $\big( \frac{t-n-1}{4n},\frac{4L^2n^2}{|S|^4} \big)$,
$\big( \frac{t-n}{4n},\frac{4g^2n^2}{|S|^4} \big)$, $\big(
\frac{t+K-n-1}{4n},\frac{4g^2n^2}{|S|^4} \big)$, and $\big(
\frac{t+K-n}{4n},\frac{4L^2n^2}{|S|^4} \big)$ (since the graph of $\ff$ is piecewise
convex). The area of this trapezoid, which bounds the contribution to $\|p\|_2^2$ from
$x$ in the interval $\big[ \frac{t-n-1}{4n},\frac{t+K-n}{4n} \big]$, equals
$\frac{Kn(g^2-L^2)}{|S|^4}$. Since the numbers $K$ from all such intervals sum to $Q$,
we have established the bound $\|p\|_2^2 \le \frac{Qn(g^2-L^2)}{|S|^4}$.

Inserting these bounds into Eq.~(\ref{spikes}), we conclude that
    \begin{equation*}
    \phi \le \frac{Qn(g^2-L^2)}{|S|^4} + \left( 1 + \frac{Q(g-L)}{2|S|^2} \right)
    \frac{2Ln}{|S|^2}
    \end{equation*}
or equivalently
    \begin{equation}
    Q \ge \frac{\phi|S|^4 - 2Ln|S|^2}{n(g-L)(g+2L)}
    \label{unscaled}
\end{equation}
as desired.
\end{proof}

We remark that, using the techniques from the proof of
Proposition~\ref{R.from.Delta.Connection.cor}(\ref{even.rho.item}), the constant
$\fstarftwonormconstant$ in the statement of Theorem~\ref{ubiquity.thm} can be replaced
with $\twotimesDconstant$ (provided one is concerned only with the case $n\to\infty$).

To give an illustration of this theorem, define $\gamma = |S|/\sqrt{gn}$, $\alpha=L/g$,
and $\kappa = {Q}/{2n}$. Then the inequality (\ref{unscaled}) becomes
\begin{equation}
\kappa \ge \frac{\gamma^2(0.574575\gamma^2 - \alpha)}{(1-\alpha)(1+2\alpha)}.
\label{complicated.bound}
\end{equation}
For example, take $\gamma=0.7$ and $\alpha=0.25$, so that the right-hand side of this
inequality is greater than $0.0137382 > \frac1{73}$. Theorem~\ref{R.Lower.Bounds.thm}
tells us that for every $g\ge2$ except $g=3$ and possibly $g=5$ and $g=7$, there exists
a $B^*[g]$ set $S$ contained in $\{1,\dots,n\}$ with at least $0.7\sqrt{gn}$ elements,
at least when $n$ is large. For every such set $S$, the inequality (\ref{complicated.bound}) asserts that at
least $\frac{2n}{73}$ of the integers $t$ between 1 and $2n$ have at least $\frac g4$
representations $t=s_1+s_2$ with $s_1,s_2\in S$.

\begin{figure}[t]
\begin{center}
    \begin{picture}(384,240)
    \put(220,2){$\gamma$}
    \put(2,130){$\alpha$}
    \includegraphics{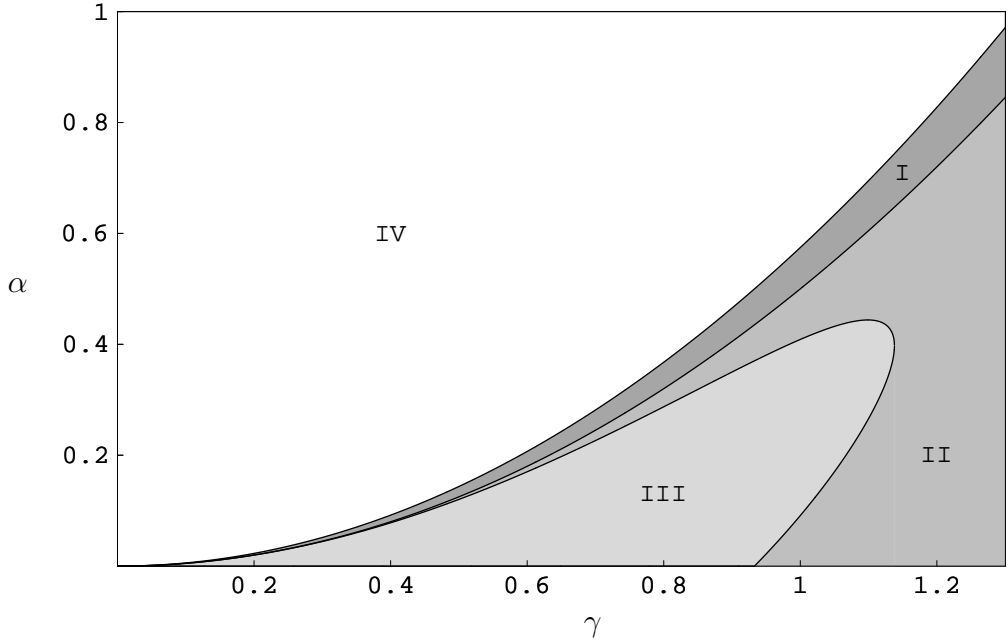}
    \end{picture}
    \caption{Comparing the bounds of Eq.~(\ref{simple.bound}) and
            Eq.~(\ref{complicated.bound})}
    \label{ubiquity.bounds.fig}
\end{center}
\end{figure}

To give a basis for comparison for Theorem~\ref{ubiquity.thm}, we note that the simple argument
\begin{equation*}
|S|^2 = \sum_t r(t) = \sum_{t\colon r(t)\le L} r(t) + \sum_{t\colon r(t)>L} r(t) \le
(2n-Q)L + Qg
\end{equation*}
gives $Q \ge \frac{|S|^2 - 2nL}{g-L}$, which translates in the above notation into
\begin{equation}
\kappa \ge \frac{\gamma^2-2\alpha}{2-2\alpha}.
\label{simple.bound}
\end{equation}
In Figure~\ref{ubiquity.bounds.fig}, Region IV corresponds to the pairs
$(\gamma,\alpha)$ for which neither Eq.~(\ref{complicated.bound}) nor
Eq.~(\ref{simple.bound}) gives a nontrivial bound on $\kappa$ (the trivial bound is
$\kappa\geq 0$). In Region I, which contains our example point $(0.7,0.25)$,
Eq.~(\ref{complicated.bound}) gives a nontrivial bound while Eq.~(\ref{simple.bound}) does not. Eq.~(\ref{simple.bound}) gives a nontrivial bound in Regions II and III, but in Region II the bound in Eq.~(\ref{complicated.bound}) is better. That is, the simple bound in Eq.~(\ref{simple.bound}) is better than Eq.~(\ref{complicated.bound}) only in Region III.

\subsection{The Uniform Distribution Hypothesis}

We say that a sequence $(S_n)_{n=1}^\infty$ of sets of positive integers {\em
becomes uniformly distributed} if the discrepancy of $S_n$ goes to 0 as
$n\to\infty$. That is, $(S_n)_{n=1}^\infty$ becomes uniformly distributed if
 $$\limsup_{n\to\infty}\quad\sup_{0\leq\alpha<\beta\leq1} \left|
    \frac{|S_n \cap [\alpha M_n,\beta M_n]|}{|S_n|}-(\beta-\alpha)\right|=0,$$
where $M_n$ denotes the largest element of $S_n$.

It has long been known~\cite{1944.Erdos,1944.Chowla} that $R(2,n) \sim \sqrt{n}$. In
1991, Erd\H{o}s and Freud~\cite{1991.Erdos.Freud} proved that if
$S_n\subseteq\{1,\dots,n\}$ is a sequence of $B^\ast[2]$ sets with $|S_n| \sim
\sqrt{n}$, then $(S_n)$ becomes uniformly distributed. We are led to make the following
conjecture.

\begin{cnj}
Let $g\ge2$ be an integer. Suppose that $S_n\subseteq\{1,\dots,n\}$ is a sequence of
$\Bg$ sets with $|S_n| \sim R(g,n)$. Then $(S_n)$ becomes uniformly distributed.
\label{uniformly.distributed.cnj}
\end{cnj}

Independent of the conjecture, we are able to prove a strong result on the
cardinality of uniformly distributed $\Bg$ sets.

\begin{thm}
\label{UD.Hypothesis.thm} Let $g\ge2$ be an integer. Let $S_n \subseteq \{1,\dots,n\}$
be a sequence of $\Bg$ sets that becomes uniformly distributed as $n\to\infty$. Then
$|S_n| \lesssim 1.15988 \sqrt{gn}$.
\end{thm}

\begin{proof}
Let $f_n$ be the nif corresponding to $S_n$. Since the $S_n$ become uniformly
distributed, the functions $f_n$ converge in measure to the function that is identically
2 on $[-\frac14,\frac14]$. Therefore
 $$\hat{f_n}(1)\to \int_{-1/4}^{1/4} 2 e^{-2\pi ix}\,dx = \frac{2}{\pi}.$$
Now, apply Lemma~\ref{Fhat.Green.Bound.lem} to find $\liminf_n \|f_n\ast
f_n\|_\infty \geq 1.486634$. Since $\|f_n\ast f_n\|_\infty = \frac{g}{2n}
\big( \frac{2n}{|S_n|} \big)^2$, we have shown $|S_n| \lesssim
1.15988 \sqrt{gn}$.
\end{proof}

In light of Theorem \ref{UD.Hypothesis.thm}, we see that Conjecture \ref{uniformly.distributed.cnj} implies that $\overline\rho(g) \le 1.15988$, improving Theorem \ref{R.Upper.Bound.thm} for all sufficiently large values of $g$.

\pagebreak \section{Constructions of \Bg Sets and Lower Bounds for $R(g,n)$}

We begin by considering a modular version of \Bg sets. A set $S$ is a \Bgm set if for
any given $m$ there are at most $g$ ordered pairs $(s_1,s_2)\in S \times S$ with
$s_1+s_2\equiv m\pmod n$ (equivalently, if the coefficients of the least-degree
representative of $\left(\sum_{n\in S} z^n\right)^2 \pmod{z^n-1}$ are bounded by $g$).
For example, the set $\{0,1,2,4\}$ is a $B^\ast[3]\pmod{7}$ set, and $\{0,1,3,7\}$ is a
$B^\ast[2]\pmod{12}$ set. Note that $7+7\equiv 1+1 \pmod{12}$, so that $\{0,1,3,7\}$ is
not a ``modular Sidon set'' as defined by some authors,
e.g.,~\cite{1980.2.Graham.Sloane} or \cite[Problem C10]{1994.Guy}. It is, however, the
natural companion to the study of \Bg sets, as evidenced by the clean forms of
Propositions~\ref{Modular.Construction.prop} and~\ref{HalfModular.Construction.prop}
below.

Just as we defined $R(g,n)$ to be the largest possible cardinality of a \Bg set
contained in $[0,n)$, we define $C(g,n)$ to be the largest possible cardinality of a
\Bgm set. The mnemonic is ``R'' for the $\R$eal problem and ``C'' for the Circular
problem. After exhibiting some basic bounds for this new function $C(g,n)$ in the next
section, we construct some explicit families of large \Bgm sets in
Section~\ref{circular.constructions.section}, which are used in turn to construct large
\Bg sets themselves in Section~\ref{constructions.bg.section}. We also demonstrate the
existence of large \Bgm sets via a probabilistic construction in Section
\ref{circular.probabilistic.section}, and we give a similar probabilistic construction
of large \Bg sets in Section~\ref{probabilistic.bg.section}. These constructions will
yield lower bounds on the size of \Bg sets, which we collect in
Section~\ref{lower.bounds.bg.section}.

\subsection{Upper Bounds on $C(g,n)$}

Since \Bgm sets have not been rigorously developed in the literature, we begin
this section by giving some simple {\em upper} bounds on $C(g,n)$.

We can obtain $$C(g,n)\leq
  \begin{cases}
    \sqrt{gn} & \text{$g$ even}, \\
    \sqrt{1-\tfrac1g}\sqrt{gn}+1 & \text{$g$ odd},
  \end{cases}$$
from the pigeonhole principle as follows. There are $|S|^2$ pairs of elements
from $S$, and there are just $n$ possible values for the sum of two elements,
and a possible value is realized at most $g$ times. Thus $|S| \leq \sqrt{gn}$.
The only way a sum can occur an odd number of times is if it is twice an
element of $S$, so for odd $g$, $|S|^2 \leq (g-1)n+|S|$.

For $g=2$ we can be more precise: $\binom{C(2,n)}{2}\leq \floor{\frac n2}$. To establish
this, let $S\subseteq[0,n)$ be a witness and note that there are $\tbinom{|S|}{2}$ pairs
of distinct elements from $S$, and each such pair $s_1,s_2$ leads to a pair of
differences $\{s_1-s_2,s_2-s_1\} \in \left\{\{i,n-i\} \colon 1\leq i < n\right\}$. If each of $(s_1,s_2)$ and $(s_3,s_4)$ is a pair of incongruent elements and if $s_1-s_2\equiv s_3-s_4 \pmod{n}$, then $s_1+s_4\equiv s_4+s_1\equiv s_2+s_3\equiv s_3+s_2 \pmod{n}$. The fact that $S$ is a {\mbox{$B^{\ast}[2]\pmod{n}$} } set forces $\{s_1,s_4\}\equiv\{s_2,s_3\}\pmod n$, and since $s_1\not\equiv s_2\pmod n$ by assumption we conclude that $s_1\equiv s_3\pmod n$ and $s_2\equiv s_4\pmod n$. Therefore distinct pairs of incongruent elements lead to distinct
sets of differences, of which there are at most $\floor{n/2}$, establishing
$\binom{C(2,n)}{2}\leq \floor{\frac n2}$. This bound is actually achieved for
$n=p^2+p+1$ when $p$ is prime (see Proposition~\ref{Modular.Constructions.prop}(iii)
below) and implies that $C(2,n)<\sqrt{n}+1$.

\subsection{Explicit Constructions of \Bgm Sets}
\label{circular.constructions.section}

We turn now to the problem of constructing large \Bgm sets. The literature contains
several examples of families of $B^\ast[2]\pmod{n}$ sets which we can generalize to
families of \Bgm sets. The following proposition collects several lower bounds for $C(g,x)$
corresponding to various constructions given in the proofs thereafter.

\begin{prop}
    \label{Modular.Constructions.prop}
    \renewcommand{\theenumi}{\roman{enumi}}
    Let $p$ be a prime power, and let $1\le k < p$.
    \begin{enumerate}
    \item If $p$ is a prime, then $C(2k^2,p^2-p) \geq k (p-1)$;
    \item $C(2k^2,p^2-1) \geq k p$;
    \item $C(2k^2,p^2+p+1) \geq kp+1$.
    \end{enumerate}
\end{prop}

The $k=1$ cases of Proposition~\ref{Modular.Constructions.prop}(i), (ii), and
(iii) are due to Ruzsa~\cite{1993.Ruzsa}, Bose~\cite{1942.Bose}, and
Singer~\cite{1938.Singer}, respectively. Part (i) uses the existence of a
primitive root modulo a prime $p$ and the Chinese Remainder Theorem (the fact
that the orders of the additive and multiplicative groups modulo $p$ are
relatively prime is important, which is why the construction does not
generalize to all finite fields). Part (ii) uses the existence of the finite
fields $\field{p^t}$ with cyclic multiplicative groups and their vector space
structures. The third part also uses a fact (Lemma~\ref{finite.field.lem})
connecting the multiplicative group of a finite field to its vector space
structure.

All three constructions make use of the following lemma. Although the lemma is
a special case of unique factorization, we give here a simple, elementary
proof of the special case that we require.

\begin{lem}
\label{distinct.pairs.lem} In any field, there are at most 2 ordered pairs
$(a,b)$ of solutions to the polynomial equation $x^2-c_1 x +c_2=(x-a)(x-b)$.
\end{lem}

\begin{proof}
Suppose that there are three solutions $(a_m,b_m)$, $1\leq m \leq 3$.
Obviously we have $a_1+b_1=c_1=a_2+b_2$ and $a_1b_1=c_2=a_2b_2$. This leads to
$0=a_1(a_1+b_1-a_2-b_2)=a_1^2+a_1b_1-a_1a_2-a_1b_2 =
a_1^2+a_2b_2-a_1a_2-a_1b_2=(a_1-a_2)(a_1-b_2)$, whence $a_1 \in \{a_2,b_2\}$.
If $a_1=a_2$, then $b_1=b_2$ and the three solutions are not distinct.
Otherwise $a_1=b_2$, and likewise $a_1=b_3$, whence $b_2=b_3$, again a
contradiction.
\end{proof}

\begin{proofof}{Proposition~\ref{Modular.Constructions.prop}(i)}
Let $g$ be a primitive root modulo $p$. Using the Chinese Remainder Theorem,
define $a_{t,i}$ for $1\leq t<p$ and $1\leq i \leq k$ by the pair of
congruences
 \begin{equation}\label{Ruzsa.def}
    a_{t,i} \equiv t \pmod{p-1} \quad\text{ and }\quad
    a_{t,i} \equiv i g^t \pmod{p}.
 \end{equation}
Set $S_i:=\{a_{t,i} : 1\leq t<p\}$; clearly $|S_i|=p-1$. we shall show that
$\bigcup_{i=1}^k S_i$ witnesses $C(2k^2,p(p-1))\geq k(p-1)$.

Suppose $a_{m_1,i}=a_{m_2,j}$, with $m_1,m_2 \in [1,p)$. We have $m_1 \equiv a_{m_1,i}
=a_{m_2,j}\equiv m_2 \pmod{p-1}$, so $m_1=m_2$. Reducing the equation
$a_{m_1,i}=a_{m_2,j}$ modulo $p$, we find $ig^{m_1}\equiv jg^{m_2}=jg^{m_1} \pmod{p}$,
so $i=j$. Thus $S_i$ and $S_j$ are distinct for distinct $i,j$, and $\left|
\bigcup_{i=1}^k S_i \right| = k(p-1)$.

First assume that, given any $i,j\in[1,k]$ and any $n\in [1,p(p-1)]$, there
are at most two pairs $(r,v)\in S_i \times S_j$ with $r+v \equiv n
\pmod{p(p-1)}$. Since there are $k^2$ such choices of $i,j$, and each such
choice leads to at most two pairs with a given sum, this shows that each $n$
arises as the sum of at most $2k^2$ pairs, concluding the proof.

Now suppose if possible that there are three pairs $(a_{r_m,i},a_{v_m,j})\in
S_i \times S_j$ satisfying $a_{r_m,i}+a_{v_m,j}\equiv n \pmod{p(p-1)}$. Then
obviously $a_{r_m,i}+a_{v_m,j}\equiv n \pmod{p}$. Also, since
$a_{r_m,i}+a_{v_m,j}\equiv n \pmod{p-1}$, we have $a_{r_m,i}\cdot a_{v_m,j} =
ig^{r_m} jg^{v_m}=ijg^{r_m+v_m} \equiv ij g^n \pmod{p}$. Thus the pairs
$(a_{r_m,i},a_{v_m,j})$ are solutions to $x^2-nx+ijg^n \pmod{p}$. Applying
Lemma~\ref{distinct.pairs.lem}, we find that two of the pairs
$(a_{r_m,i},a_{v_m,j})$ are equal, say $a_{r_1,i}\equiv a_{r_2,i} \pmod{p}$.
By Eq.~(\ref{Ruzsa.def}), $i g^{r_1}\equiv a_{r_1,i}\equiv a_{r_2,i} \equiv i
g^{r_2}\pmod{p}$. Since $g$ has order $p-1$, this tells us that $a_{r_1,i}
\equiv r_1 \equiv r_2 \equiv a_{r_2,i} \pmod{p-1}$. Since $a_{r_1,i}\equiv
a_{r_2,i} \pmod{p}$, $a_{r_1,i}\equiv a_{r_2,i} \pmod{p-1}$, and $a_{t,i}\in
[1,p(p-1))$ (by definition), we see that $a_{r_1,i}=a_{r_2,i}$, whence the
three pairs are not distinct.
\end{proofof}

\begin{proofof}{Proposition~\ref{Modular.Constructions.prop}(ii)}
Let $\theta$ generate the multiplicative group of the finite field
$\field{p^2}$, and observe that $\{1,\theta\}$ is a basis of $\field{p^2}$ as
a vector space over $\field{p}$. For $i\in\field{p}$, define
 $$S_i:=\{ s' \in [1,p^2-1] : \theta^{s'}= i \theta+ s, \quad
                                            s\in\field{p}\}.$$
Bose~\cite{1942.Bose} showed that each $S_i$ is a Sidon set.  we shall show
that $S=\bigcup_{i=1}^k S_i$ witnesses $C(2k^2,p^2-1)\geq kp$.

First, note that for each $s\in\field{p}$, there is an integer $s' \in
[1,p^2-1]$ with $\theta^{s'}= i \theta+ s$, so that $|S_i|=p$. Since
$1,\theta$ is a basis for $\field{p^2}$ over $\field{p}$, we know that $i
\theta+s_1=j\theta+s_2$ implies that $i=j$ and $s_1=s_2$. In particular, if
$i\not=j$, then $S_i \cap S_j=\emptyset$. Thus $|\bigcup_{i=1}^k S_i|= k p$.

It is sufficient to show that for each $n$ and any $i,j$ (for which there are
$k^2$ choices) there are at most 2 pairs $(r_m',v_m')$ in $S_i \times S_j$
with $r_m'+v_m' \equiv n \pmod{p^2-1}$. Define $c_1,c_2 \in \field{p}$ by
$(ij)^{-1} \theta^n-\theta^2=c_1\theta+c_2$. Since $r_m'+v_m'\equiv n
\pmod{p^2-1}$, we have
 \begin{multline*}
        c_1\theta+c_2 = (ij)^{-1} \theta^n-\theta^2 =
        (ij)^{-1} \theta^{r_m'+v_m'}-\theta^2 =
        (ij)^{-1} \theta^{r_m'}\theta^{v_m'}-\theta^2 = \\
        (ij)^{-1} (i \theta+r_m)(j \theta+v_m)-\theta^2 =
        (i^{-1}r_m+j^{-1}v_m)\theta+i^{-1}r_mj^{-1}v_m.
 \end{multline*}
This means that $(a,b)=(i^{-1}r_m,j^{-1}v_m)$ is a solution to
$x^2-c_1x+c_2=(x-a)(x-b)$. By Lemma~\ref{distinct.pairs.lem}, there are at
most two such pairs.
\end{proofof}

To extend Singer's~\cite{1938.Singer} construction, we shall need the following lemma.

\begin{lem}
\label{finite.field.lem} Let $p$ be a prime power, and let $\theta \in
\field{p^3}$ generate the multiplicative group. Then $\theta^a$ and $\theta^b$
are linearly dependent over $\field{p}$ iff $a\equiv b \pmod{p^2+p+1}$.
\end{lem}

\begin{proof}
The multiplicative group of $\field{p}$ is a subgroup of the multiplicative group of
$\field{p^3}$, i.e., $\field{p}=\{\theta^k : k\equiv 0 \pmod{\frac{p^3-1}{p-1}}\}$.
Since two elements of $\field{p^3}$ are linearly dependent over $\field{p}$ exactly if
their ratio is in $\field{p}$, we see that $\theta^a$ and $\theta^b$ are linearly
dependent exactly if $a-b \equiv 0 \pmod{\frac{p^3-1}{p-1}}$.
\end{proof}

\begin{proofof}{Proposition~\ref{Modular.Constructions.prop}(iii)}
Let $\theta \in \field{p^3}$ generate the multiplicative group. Since $\theta$
is algebraic with degree 3 over $\field{p}$, the elements $1, \theta,
\theta^2$ are a basis for $\field{p^3}$ over $\field{p}$. Define
 $$ T_i :=  \{0\} \cup
            \{ s' \in [1,p^3-1] : \theta^{s'} = i \theta + s,\quad
                                                         s\in\field{p}\} $$
for $1\leq i \leq k$, and define $S_i$ to be the set of congruence classes
modulo $q:=\tfrac{p^3-1}{p-1}=p^2+p+1$ which intersect $T_i$.

Since each nonzero $s'\in T_i$ corresponds to a unique $s\in\field{p}$, and
for each $s\in\field{p}$ there is an $s'\in[1,p^3-1]$ with
$\theta^{s'}=i\theta+s$, we see that $|T_i|=|\field{p}|+1=p+1$, and so by
virtue of Lemma~\ref{finite.field.lem}, $|S_i|=p+1$. Also, each $s'\not=0$
occurs in at most one of the $T_i$, so that $|\bigcup_{i=1}^k S_i|=kp+1$.

We wish to show that, given any $n$, there are at most two pairs $(r',v')\in
T_h \times T_j$ with $r'+v' \equiv n \bmod{q}$. Since there are $k^2$ choices
of $h,j$, this will establish that $\bigcup_{i=1}^k S_i$ witnesses
$C(2k^2,p^2+p+1)\geq kp+1$. Define $L_0=1$, and for each $x'\in T_i$ set
$L_{x'}=\theta^{x'}=i \theta+x$.

Suppose that $n$ is fixed, and $(r'_m,v'_m)\in T_h \times T_j$ ($1\leq m \leq
3$) are distinct pairs with $r_m'+v'_m \equiv n \bmod{q}$ for $1\leq m \leq
3$. It follows from Lemma~\ref{finite.field.lem} that each pair of
$L_{r'_1}L_{v'_1}, L_{r'_2}L_{v'_2}, L_{r'_3}L_{v'_3}$ are linearly dependent,
i.e., they are multiples of each other. If $r'_1=0$, then $L_{r'_1}L_{v'_1}$
is linear. This means that $L_{r'_2}L_{v'_2}$ is also linear, and so one of
$r'_2,v'_2$ is zero and the other is equal to $v'_1\in T_j$. If $r'_2=0$, then
$(r'_1,v'_1)=(r'_2,v'_2)$, contradicting distinctness, and if $v'_2=0$, then
$v'_1\equiv r'_2 \pmod{q}$, which is only possible if $v'_1=r'_2=0$, whence
again $(r'_1,v'_1)=(r'_2,v'_2)$.

Thus each $L_{r'_m}L_{v'_m}$ is a quadratic in $\theta$. The coefficient of
$\theta^2$ in each $\theta^n=L_{r'_m}L_{v'_m}=(h
\theta+r_m)(j\theta+v_m)=hj\theta^2+(r_m j+v_m h)\theta+r_mv_m$ is $hj$. Since
the $L_{r'_m}L_{v'_m}$ are multiples of each other with the same lead
coefficients, they must in fact be equal. Set $c_1,c_2$ by
$(hj)^{-1}\theta^n-\theta^2=c_1\theta+c_2$, and observe that
$(a,b)=(h^{-1}r_m,j^{-1}v_m)$ is a solution to $x^2-c_1x+c_2=(x-a)(x-b)$;
hence, by Lemma~\ref{distinct.pairs.lem},  there are only two such pairs.
\end{proofof}

We now show how to combine two \Bgm sets to construct another.

\begin{prop}
\label{Modular.Construction.prop} Let $\gcd(x,y)=1$, and let $S$ be a $B^\ast[g]
\pmod{x}$ set and $M$ be a $B^\ast[h] \pmod{y}$ set. Then the set
$M+yS := \{m+ys\colon m\in M, s\in S\}$ is a $B^\ast[gh] \pmod{xy}$ set. In
particular, $C(gh,xy)\geq C(g,x) C(h,y)$ for any positive integers $g,h,x,y$ with
$\gcd(x,y)=1$.
\end{prop}

\begin{proof}
Consider $m_i,n_i\in M$ and $s_i,t_i\in S$ with
    \begin{equation}
        \label{mainsupposition}
    (m_1+ys_1)+(n_1+yt_1)\equiv
    \dots\equiv(m_{gh+1}+ys_{gh+1})+
                            (n_{gh+1}+yt_{gh+1}) \pmod{xy}.
    \end{equation}
We need to show that $m_i=m_j$, $s_i=s_j$, $n_i=n_j$, and $t_i=t_j$, for some $i,j$.
Reducing these Eq.~(\ref{mainsupposition}) modulo $y$, we see that $m_1+n_1\equiv
m_2+n_2\equiv \dots\equiv m_{gh+1}+n_{gh+1}\pmod{y}$. Since $M$ is a $B^\ast[h]
\pmod{y}$ set, we can reorder the $m_i,n_i,s_i,t_i$ so that $m_1=m_2=\dots=m_{g+1}$ and
$n_1=n_2=\dots=n_{g+1}$. Reducing Eq.~(\ref{mainsupposition}) modulo $x$ we arrive at
$$
y s_1+y t_1\equiv y s_2+ yt_2\equiv \dots\equiv y s_{g+1}+y t_{g+1} \pmod{x}
$$
whence, since $\gcd(x,y)=1$,
$$
s_1+t_1\equiv s_2+t_2\equiv \dots \equiv s_{g+1}+t_{g+1} \pmod{x}.
$$
The $s_i$ and $t_i$ are from a $B^\ast[g]\pmod{x}$ set, so that for some
$i,j$, $s_i=s_j$ and $t_i=t_j$.
\end{proof}

We have computed $C(g,n)$ for small $g$ and $n$ by exhaustive search. The results are
summarized in Table~\ref{C.Table}. The entry for $(k,g)=(10,5)$ is 28; this means that
$C(5,28)\ge10$ (a witness is $\{0, 1, 2, 4, 5, 8, 12, 15, 23, 24\}$), while there is no
$n<28$ for which $C(5,n)\ge10$.

It is straightforward to verify that $C(5,28)\geq 10$: one simply verifies that the
witness has 10 elements and is indeed a $B^\ast[5]\pmod{28}$ set. It is not
straightforward, however, to verify that $C(5,n) <10$ for $n<28$. We have made these
verifications for each entry given in the table by a long computer search using {\em
Mathematica}.

\begin{table}\footnotesize\begin{center}
$g$\vskip4pt $k$\quad
\begin{tabular}{|c||r|r|r|r|r|r|r|r|r|r}
\hline
    &         2 &  3 &  4 &  5 &  6 &  7 &  8 &  9 & 10 & 11 \\ \hline \hline
  3 &         6 &    &    &    &    &    &    &    &    &    \\ \hline
  4 &        12 &  7 &    &    &    &    &    &    &    &    \\ \hline
  5 &        21 & 11 &  8 &    &    &    &    &    &    &    \\ \hline
  6 &        31 & 19 & 11 &  9 &    &    &    &    &    &    \\ \hline
  7 &        48 & 29 & 14 & 13 & 10 &    &    &    &    &    \\ \hline
  8 &        57 & 43 & 22 & 17 & 12 & 11 &    &    &    &    \\ \hline
  9 &        73 & 57 & 28 & 19 & 16 & 13 & 12 &    &    &    \\ \hline
 10 &        91 &    & 36 & 28 & 19 & 17 & 14 & 13 &    &    \\ \hline
 11 &           &    &    & 35 & 22 & 21 & 18 & 15 & 14 &    \\ \hline
 12 &           &    &    &    & 30 & 23 & 21 & 19 & 16 & 15 \\ \hline
 13 &           &    &    &    &    & 31 & 24 & 22 & 19 & 17 \\ \hline
 14 &           &    &    &    &    &    & 28 & 25 &    & 20 \\ \hline
\end{tabular}
\caption{$\min\{n\colon C(g,n) \geq k\}$}\label{C.Table}
\end{center}\end{table}

\subsection{Probabilistic Constructions of \Bgm Sets}
\label{circular.probabilistic.section}

The algebraic methods of the previous section provide effective and
completely explicit constructions of large \Bgm sets. However, we can
establish the existence of even larger \Bgm sets using a probabilistic
construction. We rely upon the following two lemmas, which are quantitative
statements of the Law of Large Numbers for sums of many independent random
variables.

%

\begin{lem}
Let $p_1,\dots,p_n$ be real numbers in the range $[0,1]$, and set
$p=(p_1+\dots+p_n)/n$. Define mutually independent random variables
$X_1,\dots,X_n$ such that $X_i$ takes the value $1-p_i$ with probability
$p_i$ and the value $-p_i$ with probability $1-p_i$ (so that the expectation
of each $X_i$ is zero), and define $X=X_1+\dots+X_n$. Then for any positive
number $a$,
\begin{equation*}
\Prob [ X>a ] < \exp \Big( \frac{-a^2}{2pn} + \frac{a^3}{2p^2n^2} \Big)
\quad\text{and}\quad
\Prob [ X<-a ] < \exp \Big( \frac{-a^2}{2pn} \Big).
\end{equation*}
\label{prob.method.lem}
\end{lem}

\begin{proof}
These assertions are Theorems A.11 and A.13 of~\cite{2000.Alon.Spencer}.
\end{proof}

\begin{lem}
Let $p_1,\dots,p_n$ be real numbers in the range $[0,1]$, and set
$E=p_1+\dots+p_n$. Define mutually independent random variables
$Y_1,\dots,Y_n$ such that $Y_i$ takes the value 1 with probability
$p_i$ and the value 0 with probability $1-p_i$, and define $Y=Y_1+\dots+Y_n$
(so that the expectation of $Y$ equals $E$). Then $\Prob [ Y>E+a ] < \exp
\big( \frac{-a^2}{3E} \big)$ for any real number
$0<a<E/3$, and $\Prob [ Y<E-a ] < \exp \big( \frac{-a^2}{2E} \big)$ for any
positive real number $a$.
\label{our.application.lem}
\end{lem}

\begin{proof}
This follows immediately from Lemma~\ref{prob.method.lem} upon defining $X_i = Y_i -
p_i$ for each $i$ and noting that $E=pn$ and that $\frac{a^3}{2E^2} < \frac{a^2}{6E}$
under the assumption $0<a<E/3$.
\end{proof}

We now give the probabilistic construction of large \Bgm sets.

\begin{prop}
For every $0<\e\le1$, there is a sequence of ordered pairs $(n_j,g_j)$ of positive
integers such that $\frac{C(g_j,n_j)}{n_j} \gtrsim \e$ and $\frac{g_j}{n_j}
\lesssim \e^2$.
\label{circle.probabilistic.prop}
\end{prop}

\begin{proof}
Let $n$ be an odd integer. We define a random subset $S$ of $\{1,\dots,n\}$ as follows:
for every $1\leq i \leq n$, let $Y_i$ be 1 with probability $\e$ and 0 with probability
$1-\e$ with the $Y_i$ mutually independent, and let $S := \{i \colon Y_i = 1\}$. We see that $|S|=\sum_{i=1}^n Y_i$ has
expectation $E=\e n$. Setting $a = \sqrt{\e n\log 4}$, Lemma~\ref{our.application.lem}
gives
$$\Prob \big[|S| < \e n - \sqrt{\e n\log 4}\, \big] < \frac12.$$

Now for any integer $k$, define the random variable
    \begin{align*}
    R_k &:= \# \{1\le c,d\le n \colon c+d\equiv k \pmod{n},\, Y_c=Y_d=1\}\\
        &= \sum_{c+d\equiv k \pmod{n}} Y_c Y_d,
    \end{align*}
so that $R_k$ is the number of representations of $k\pmod{n}$ as the sum of two elements
of $S$. Observe that $R_k$ is the sum of $n-1$ random variables taking the value 1 with
probability $\e^2$ and the value 0 otherwise, plus one random variable (corresponding to
$c\equiv d\equiv 2^{-1}k \pmod n$) taking the value 1 with probability $\e$ and the
value 0 otherwise. Therefore the expectation of $R_k$ is $E=(n-1)\e^2+\e$. Setting $a =
\sqrt{3((n-1)\e^2+\e)\log 2n}$, and noting that $a<E/3$ when $n$ is sufficiently large
in terms of $\e$, Lemma \ref{our.application.lem} gives
$$
\Prob \big[R_k > (n-1)\e^2+\e - \sqrt{3((n-1)\e^2+\e)\log 2n}\, \big] <
\frac1{2n}
$$
for each $1\le k\le n$. Therefore, there exists a $\Bgm$ set $S\subseteq \{1,\dots,n\}$,
with $g \le (n-1)\e^2+\e - \sqrt{3((n-1)\e^2+\e)\log 2n} \lesssim \e^2n$, such that $|S|
\ge \e n - \sqrt{\e n\log 4} \gtrsim \e n$. This establishes the proposition.
\end{proof}

Define $\Delta_\T(\e)$ to be the supremum of those real numbers $\delta$ such
that every subset of $\T$ with measure $\e$ has a subset with measure $\delta$
that is fixed by a reflection $t\mapsto c-t$. The function $\Delta_\T(\e)$
stands in relation to $C(g,n)$ as $\De$ stands to $R(g,n)$. However, it turns
out that $\Delta_\T$ is much easier to understand:

\begin{cor}
Every subset of $\T$ with measure $\e$ contains a symmetric subset
with measure $\e^2$, and this is best possible for every $\e$. In particular,
$\Delta_\T(\e) = \e^2$ for all $0\le\e\le1$.
\label{circle.answer.cor}
\end{cor}

\begin{proof}
In the proof of the trivial lower bound for $\De$ (Lemma
\ref{Delta.Trivial.Lower.Bound.lem}), we saw that every subset of $[0,1]$ with measure
$\e$ contains a symmetric subset with measure at least $\frac12\e^2$. The proof is
easily modified to show that every subset of $\T$ with measure $\e$ contains a symmetric
subset with measure $\e^2$. This shows that $\Delta_\T(\e) \ge \e^2$ for all $\e$. On
the other hand, the proof of Proposition~\ref{Delta.from.R.Connection.prop} is also
easily modified to show that $\Delta_\T\big( \frac{C(g,n)}n \big) \le \frac gn$, as is
the proof of Lemma~\ref{Lipschitz.Condition.lem} to show that $\Delta_\T$ is continuous.
Then, by virtue of Proposition~\ref{circle.probabilistic.prop} and the monotonicity of
$\Delta_\T$, we have $\Delta_\T(\e) \le \e^2$.
\end{proof}

\subsection{Explicit Constructions of \Bg Sets}
\label{constructions.bg.section}

If $S\subseteq[0,100)$ is a $B^\ast[g]\pmod{200}$ set, then the modular sums are the
same as the real sums, and so $S$ is a \Bg set as well. This observation is the
fundamental idea behind using the method of Proposition~\ref{Modular.Construction.prop}
to construct \Bg sets.

\begin{prop}
\label{HalfModular.Construction.prop} Let $g,h,x,y$ be positive integers. Then
 $$R(gh,xy) \geq R(gh,xy+1-\ceiling{\tfrac{y}{C(h,y)}})
            \geq R(g,x)C(h,y).$$
\end{prop}

\begin{proof}
Let $M\subseteq[1,m]\subseteq[1,y]$ witness the value of $C(h,y)$, and let
$S\subseteq[0,s)$ witness the value of $R(g,s)$. Take $x>2s$, and relatively prime to
$y$, and note that $S$ is a $B^\ast[g] \pmod{x}$ set. By
Proposition~\ref{Modular.Construction.prop}, the set $M+y S$ is a $B^\ast[gh] \pmod{xy}$
set. But by taking $x$ to be sufficiently large, we see that $M+yS$ is actually a
$B^\ast[gh]$ set.

We now compute the smallest and largest element of $M+yS$. Clearly the smallest element
is 1, and the largest is $m+y(s-1)$. Since $M$ is a $B^\ast[h]\pmod{y}$ set, we may
shift it modulo $y$ so as to minimize $m$. $M$ has $C(h,y)$ elements, so there must be
two consecutive elements whose difference $\!\pmod{y}$ is at least
$\ceiling{\tfrac{y}{C(h,y)}}$, i.e., we may take $m \leq y+1-
\ceiling{\tfrac{y}{C(h,y)}}$. Thus $M+yS \subseteq
[1,y+1-\ceiling{\tfrac{y}{C(h,y)}}+y(s-1)] = [1,ys+1-\ceiling{\tfrac{y}{C(h,y)}}]$, and
$|M+yS|=R(g,s)C(h,y)$.
\end{proof}

The reader might feel that the part of the argument concerning the largest gap in $M$ is
more trouble than it is worth. We include this for two reasons. First,
Erd\H{o}s~\cite[Problem C9]{1994.Guy} offered \$500 for an answer to the question, ``Is
$R(2,n)=\sqrt{n}+\bigO{1}$?'' This question would be answered in the negative if one
could show, for example, that the $B^\ast[2]\pmod{p^2-1}$ sets constructed by Bose (the
$k=1$ case of Proposition~\ref{Modular.Constructions.prop}(ii)) contain a gap which is
not $\bigO{p}$, as seems likely from the experiments of Zhang~\cite{1994.Zhang} and
Lindstr\"{o}m~\cite{1998.2.Lindstrom}. Second, there is some literature (e.g.,
\cite{1995.Erdos.Sarkozy.Sos} and \cite{1996.Ruzsa}) concerning the possible size of the
largest gap in a maximal Sidon set contained in $\{1,\dots,n\}$.  In short, we include
this argument because there is some reason to believe that this is a significant source
of the error term in at least one case, and because there is some reason to believe that
improvement is possible.

Our plan is to employ the inequality of
Proposition~\ref{HalfModular.Construction.prop} when $y$ is large, $h=2$, and
$x\approx \frac83 g$. In other words, we need nontrivial lower bounds for
$C(2,n)$ for $n\to\infty$ and for $R(g,n)$ for values of $n$ that are not much
larger than $g$. The first need is filled by
Proposition~\ref{Modular.Constructions.prop}, while the second need is filled by
the following lemma.

\begin{lem}
For all $g\ge1$ we have $R(g,3g-\floor{g/3}+1) \geq
g+2\floor{g/3}+\floor{g/6}$.
\label{R.Construction.Small.gn.lem}
\end{lem}

\begin{proof}
One can verify that a witness is
\begin{equation*}
\Big[ 0, \floor{\frac g3} \Big) \cup \bigg\{ g - \floor{\frac g3} + 2 \Big[ 0,
\floor{\frac g6} \Big) \bigg\} \cup \Big[ g, g+\floor{\frac g3} \Big) \cup
\Big( 2g-\floor{\frac g3}, 3g-\floor{\frac g3} \Big].
\end{equation*}
\end{proof}

We remark that this family of examples was motivated by the finite sequence
$S=\{1,0,\frac12,1,0,1,1,1\}$, which has the property that its
autocorrelations are small relative to the sum of its entries. In other
words, the ratio of the $\ell^\infty$-norm of $S*S$ to the $\ell^1$-norm of $S$
itself is small. If we could find a finite sequence of rational numbers for
which the corresponding ratio were smaller, we could convert it directly into
a family of examples that would improve the lower bound for
$\underline\rho(2g)$ in Theorem~\ref{R.Lower.Bounds.thm} for large $g$ (see
the proof of the theorem in Section~\ref{lower.bounds.bg.section}).

In addition to these parametric results, we have established by direct
(exhaustive) computation the exact value of $R(g,n)$ for small values of $g$
and $n$. Table~\ref{R.Table} records, for given values of $g$ and
$k$, the smallest possible value of $\max S$ among all $B^\ast[g]$ sets $S$
consisting of exactly $k$ positive integers; in other words, the entry
corresponding to $k$ and $g$ is $\min\{n\colon R(g,n) \geq k\}$. For example,
the $(k,g)=(8,2)$ entry records the fact that there exists an 8-element Sidon
set of integers from $[1,35]$ but no 8-element Sidon set of integers from
$[1,34]$.

To show that $R(2,35)\geq 8$, for instance, it is only necessary to observe that the
witness $\{1, 3, 13, 20, 26, 31, 34, 35\}$ has 8 elements and is a $B^\ast[2]$ set. To
show that $R(2,35)\leq 8$, however, seems to require an extensive search.

\begin{table}\footnotesize\begin{center}
$g$\vskip4pt $k$\quad
\begin{tabular}{|c||r|r|r|r|r|r|r|r|r|r|}
\hline
    &          2 &          3 &         4 &         5 &         6 & 7 & 8 & 9 &    10     &    11     \\ \hline\hline
 3  &          4 &            &           &           &           & & & & &           \\ \hline
 4  &          7 &          5 &           &           &           & & & & &           \\ \hline
 5  &         12 &          8 &         6 &           &           & & & & &           \\ \hline
 6  &         18 &         13 &         8 &         7 &           & & & & &           \\ \hline
 7  &         26 &         19 &        11 &         9 &         8 & & & & &           \\ \hline
 8  &         35 &         25 &        14 &        12 &        10 & 9 & & & &           \\ \hline
 9  &         45 &         35 &        18 &        15 &        12 & 11 & 10 &           &           &           \\ \hline
 10 &         56 &         46 &        22 &        19 &        14 & 13 & 12 &        11 &           &           \\ \hline
 11 &         73 &         58 &        27 &        24 &        17 & 15 & 14 &        13 &    12     &           \\ \hline
 12 &  $\leq 92$ &  $\leq 72$ &        31 &        29 &        20 & 18 & 16 &        15 &    14     &        13 \\ \hline
 13 & $\leq 143$ & $\leq 101$ &        37 &        34 &        24 & 21 & 18 &        17 &    16     &        15 \\ \hline
 14 &            & $\leq 128$ &        44 &        40 &        28 & 26 & 21 &        19 &    18     &        17 \\ \hline
 15 &            &            & $\leq 52$ & $\leq 47$ &        32 & 29 & 24 &        22 &    20     &        19 \\ \hline
 16 &            &            &           &           &        36 & 34 & 27 &        24 &    22     &        21 \\ \hline
 17 &            &            &           &           & $\leq 42$ & $\leq 38$ &        30 &        28 &    24     &        23 \\ \hline
 18 &            &            &           &           &           & & 34 & 32 &    27     &        25 \\ \hline
 19 &            &            &           &           &           & & $\leq 38$ & $\leq 36$ &    30     &        28 \\ \hline
 20 &            &            &           &           &           & & & & 33     &        31 \\ \hline
 21 &            &            &           &           &           & & & & $\leq 37$ &        35 \\ \hline
 21 &            &            &           &           &           & & & & & $\leq 38$ \\ \hline
\end{tabular}
\caption{$\min\{n\colon R(g,n) \geq k\}$}\label{R.Table}
\end{center}\end{table}

\begin{table}\footnotesize\begin{center}
 $
\begin{array}{ccccr}
   g    &  x & R(g,x) &                  \text{Witness}                   & R(g,x)/\sqrt{2gx}             \\
\hline
   2    &  7 &   4    &                    \{1,2,5,7\}                    & \tfrac{2}{\sqrt7}            \approx 0.756    \\
   3    &  5 &   4    &                    \{1,2,3,5\}                    & \tfrac{2\sqrt2}{\sqrt{15}}   \approx 0.730    \\
   4    & 31 &  12    &         \{1,2,4,10,11,12,14,19,25,26,30,31\}      & \frac2{\sqrt7}               \approx 0.756    \\
   5    &  9 &   7    &               \{1,2,3,4,5,7,9\}                   & \tfrac7{3\sqrt{10}}          \approx 0.738    \\
   6    & 20 &  12    &         \{1,2,3,4,5,6,9,10,13,15,19,20\}          & \tfrac{\sqrt3}{\sqrt5}       \approx 0.775    \\
   7    & 15 &  11    &       \{1, 2, 3, 7, 8, 9, 10, 11, 12, 13, 15\}    & \tfrac{11}{\sqrt{210}}       \approx 0.759    \\
   8    & 30 &  17    &  \{1,2,5,7,8,9,11,12,13,14,16,18,26,27,28,29,30\} & \tfrac{17}{4\sqrt{30}}       \approx 0.776    \\
   9    & 24 &  16    &     \{1,2,3,4,5,6,7,8,9,13,14,15,17,22,23,24\}    & \tfrac{4}{3\sqrt{3}}         \approx 0.770    \\
   10   & 33 &  20 &\{1,2,4,5,6,7,8,9,10,11,13,15,20,21,22,23,30,31,32,33\}& \tfrac{2\sqrt{5}}{\sqrt{33}} \approx 0.778    \\
   11   & 25 &  18 &\{ 1,2,3,4,5,11,12,13,14,15,16,17,18,19,20,21,23,25\} & \tfrac{18}{5\sqrt{22}} \approx 0.768 \\
\end{array}
$ \caption{Important values of $R(g,x)$ and witnesses}\label{Witness.Table}
\end{center}\end{table}

\subsection{Probabilistic Constructions of \Bg Sets}
\label{probabilistic.bg.section}

We can use the probabilistic methods employed in Section
\ref{circular.probabilistic.section} to construct large \Bg sets in $\Z$. The proof is
more complicated because it is to our advantage to endow different integers with
different probabilities of belonging to our random set. Although all of the constants in
the proof could be made explicit, we are content with inequalities having error terms
involving big-O notation.

\begin{prop}
Let $\gamma\ge\pi$ be a real number and $n\ge\gamma$ be an integer. There exists a \Bg
set $S\subseteq\{1,\dots,n\}$, where $g = \gamma + O(\sqrt{\gamma\log n})$, with $|S|
\ge 2\sqrt{\frac{\gamma n}\pi} + O(\gamma + (\gamma n)^{1/4})$.
\label{probabilistic.bg.prop}
\end{prop}

\begin{proof}
Define mutually independent random variables $Y_k$, taking only the values 0 and 1, by
\begin{equation}
\Pr\{Y_k=1\} = p_k := \begin{cases}
1, & 1\le k < \frac\gamma\pi, \\
\sqrt{\!\frac\gamma{\pi k}}, & \frac\gamma\pi \le k \le n, \\
0, & k > n.
\end{cases}
\label{schinzel.pk.def}
\end{equation}
(Notice that $p_k \le \sqrt{\!\frac\gamma{\pi k}}$ for all $k\ge1$.) These
random variables define a random subset $S = \{k\colon Y_k=1\}$ of the integers
from 1 to~$n$. We shall show that, with positive probability, $S$ is a large
$B^*[g]$ set with $g$ not much bigger than $\gamma$.

The expected size of $S$ is
\begin{align}
E_0 := \sum_{1\le j\le n} p_j &= \sum_{1\le j<\gamma/\pi} 1 +
\sum_{\gamma/\pi\le j\le n} \sqrt{\frac\gamma{\pi j}} \notag \\
&= \frac\gamma\pi + \int_{\gamma/\pi}^n \sqrt{\frac\gamma{\pi t}} \,dt + O(1)
= 2\sqrt{\frac{\gamma n}\pi} - \frac\gamma\pi + O(1).
\label{second.main.term.eq}
\end{align}
If we set $a_0 := \sqrt{2E_0\log3}$, then Lemma~\ref{our.application.lem} tells us that
    $$\Prob[|S| < E_0 - a_0] < \exp\big( \frac{-a_0^2}{2E_0} \big) = \frac13.$$

Now for any integer $k\in[\gamma,2n]$, let
\begin{equation*}
R_k := \sum_{1\le j\le n} Y_jY_{k-j} = 2\sum_{1\le j<k/2} Y_jY_{k-j} + Y_{k/2},
\end{equation*}
the number of representations of $k$ as $k=s_1+s_2$ with $s_1,s_2\in S$. (Here
we adopt the convention that $Y_{k/2}=p_{k/2}=0$ if $k$ is odd). Notice that
in this latter sum, $Y_{k/2}$ and the $Y_jY_{k-j}$ are mutually independent
random variables taking only values 0 and 1, with $\Pr[Y_jY_{k-j}=1] =
p_jp_{k-j}$. Thus the expectation of $R_k$ is
\begin{align}
E_k := 2\sum_{1\le j<k/2} p_jp_{k-j} + p_{k/2} &\le 2\sum_{1\le j<k/2}
\sqrt{\frac\gamma{\pi j}} \, \sqrt{\frac\gamma{\pi(k-j)}} +
\sqrt{\frac{\gamma}{\pi k/2}} \notag \\
&\le \frac{2\gamma}\pi \int_0^{k/2} \sqrt{\frac1{t(k-t)}} \,dt +
\sqrt{\frac{2\gamma}{\pi k}} = \gamma + \sqrt{\frac{2\gamma}{\pi k}} < \gamma+1
\label{pi.over.2.appearance}
\end{align}
using the inequalities $p_k \le \sqrt{\!\frac\gamma{\pi k}}$ and $k\ge\gamma$.

If we set $a = \sqrt{3(\gamma+1)\log3n}$, then Lemma~\ref{our.application.lem} tells us
that
\begin{equation*}
\Prob[R_k > \gamma + 1 + a] < \Prob[R_k > E_k + a] < \exp\Big(
\frac{-a^2}{3E_k} \Big) < \exp\Big( \frac{-a^2}{3(\gamma+1)} \Big) =
\frac1{3n}
\end{equation*}
for every $k$ in the range $\gamma\le k\le2n$. Note that $R_k\le\gamma$
trivially for $k$ in the range $1\le k\le\gamma$. Therefore, with probability
at least $1-\frac13-(2n-\gamma)\frac1{3n} = \frac\gamma{3n}>0$, the set $S$
has at least $E_0 - a_0 = 2\sqrt{\frac{\gamma n}\pi} + O(\gamma + (\gamma
n)^{1/4})$ elements and satisfies $R_k \le \gamma + 1 + a$ for all $1\le
k\le2n$. Setting $g:=\gamma+1+a = \gamma + O(\sqrt{\gamma\log n})$, we conclude
that any such set $S$ is a \Bg set. This establishes the proposition.
\end{proof}

Schinzel conjectured that among all pdfs supported on $[0,\frac12]$, the function
\begin{equation*}
f(x) = \begin{cases}
\tfrac1{\sqrt{2x}}, & x\in[0,\tfrac12], \\
0, & \text{otherwise}
\end{cases}
\end{equation*}
has the property that $\ffi$ is minimal. We have
\begin{equation*}
\ff(x) = \begin{cases}
\tfrac\pi2, & x\in[0,\tfrac12], \\
\tfrac\pi2 - 2\mathop{\rm tan}^{-1}\sqrt{2x-1}, & x\in[\frac12,1], \\
0, & \text{otherwise}
\end{cases}
\end{equation*}
and so $\ffi = \frac\pi2$. We have adapted this function in our definition
(\ref{schinzel.pk.def}) of the probabilities $p_k$; the constant $\frac\pi2$ appears as
the value of the last integral in Eq.~(\ref{pi.over.2.appearance}). If Schinzel's
conjecture were false, then we could immediately incorporate any better function $f$
into the proof of Proposition~\ref{probabilistic.bg.prop} and improve the lower bound on
$|S|$. Indeed, Schinzel's conjecture is one of the motivations for our Conjecture
\ref{uniformly.distributed.cnj}, which by the above discussion is logically stronger.

\medskip\noindent
{\bf Theorem~\ref{g.growing.lower.bound.thm}.} {\it For any $\delta>0$, we have
$R(g,n)>\big( \frac2{\sqrt\pi} -\delta \big) \sqrt{gn}$ if both $\frac g{\log n}$ and $\frac
ng$ are sufficiently large in terms of~$\delta$.}
\medskip

\begin{proof}
In the proof of Proposition~\ref{probabilistic.bg.prop}, we saw that $\gamma\le g$ and
$g=\gamma+O(\sqrt{\gamma\log n})$; this implies that $\gamma=g+O(\sqrt{g\log n}) = g
\big( 1 + O\big( \sqrt{\frac{\log n}g} \big) \big)$. Therefore the size of the
constructed set $S$ was at least
\begin{align*}
2\sqrt{\tfrac{\gamma n}\pi} + O(\gamma + (\gamma n)^{1/4}) &= 2\sqrt{\tfrac
{gn}\pi \Big( 1 + O\Big( \sqrt{\tfrac{\log n}g} \,\Big) \Big)} +
O(g+(gn)^{1/4}) \\
&= 2\sqrt{\tfrac{gn}\pi} \Big( 1 + O \Big( \sqrt{\tfrac{\log n}g} +
\sqrt{\tfrac gn} \, \Big) \Big).
\end{align*}
This establishes the theorem.
\end{proof}

\subsection{Lower Bounds on $R(g,n)$}
\label{lower.bounds.bg.section}

We are now ready to prove Theorem~\ref{R.Lower.Bounds.thm}, which we restate here for
the reader's convenience. Recall that $\underline{\rho}(g) = \liminf_{n\to\infty}
\frac{R(g,n)}{\sqrt{gn}}$.

\medskip\noindent
{\bf Theorem~\ref{R.Lower.Bounds.thm}.} {\it We have
$$
\begin{array}{r@{{}\ge{}}l@{{}>{}}l}
    \underline{\rho}(4)  & \tfrac{2}{\sqrt7}          & 0.755, \\ \vspace{1mm}
    \underline{\rho}(6)  & \tfrac{2\sqrt2}{\sqrt{15}} & 0.730, \\ \vspace{1mm}
    \underline{\rho}(8)  & \tfrac2{\sqrt7}            & 0.755, \\ \vspace{1mm}
    \underline{\rho}(10) & \tfrac7{3\sqrt{10}}        & 0.737, \\ \vspace{1mm}
    \underline{\rho}(12) & \tfrac{\sqrt3}{\sqrt5}     & 0.774,
\end{array}
\qquad
\begin{array}{r@{{}\ge{}}l@{{}>{}}l}
    \underline{\rho}(14) & \tfrac{11}{\sqrt{210}}     & 0.759, \\ \vspace{1mm}
    \underline{\rho}(16) & \tfrac{17}{4\sqrt{30}}     & 0.775, \\ \vspace{1mm}
    \underline{\rho}(18) & \tfrac{4}{3\sqrt{3}}       & 0.769, \\ \vspace{1mm}
    \underline{\rho}(20) & \tfrac{2\sqrt{5}}{\sqrt{33}}& 0.778,\\ \vspace{1mm}
    \underline{\rho}(22) & \tfrac{18}{5\sqrt{22}}     & 0.767,
\end{array}
$$
and for any $g\ge12$,
$$
\underline{\rho}(2g) \ge \frac{g+2\floor{g/3}+\floor{g/6}}
{\sqrt{6g^2-2g\floor{g/3}+2g}}\,.
$$
In particular, for any $\delta>0$ we have
$R(g,n)>(\frac{11}{8\sqrt3}-\delta)\sqrt{gn}$ if both $g$ and $\frac ng$ are
sufficiently large in terms of~$\delta$.}
\medskip

\begin{proof}
For any positive integers $x$ and $m\le\sqrt{n/x}$, the monotonicity of $R$ in the
second variable gives $R(2g,n) \ge R(2g,x(m^2-1)) \ge R(g,x)C(2,m^2-1)$ by
Proposition~\ref{HalfModular.Construction.prop}.  If we choose $m$ to be the largest
prime not exceeding $\sqrt{n/x}$ (so that $m\gtrsim\sqrt{n/x}$ by the Prime Number
Theorem), then Proposition~\ref{HalfModular.Construction.prop} gives $R(2g,n) \ge
R(g,x)\cdot m \gtrsim R(g,x)\sqrt{\frac nx}$ for any fixed positive integer $g$, and
hence $$\underline{\rho}(2g) = \liminf_{n\to\infty} \frac{R(2g,n)}{\sqrt{2gn}} \ge
\frac{R(g,x)}{\sqrt{2gx}}.$$

The problem now is to choose $x$ so as to make $\frac{R(g,x)}{\sqrt{2gx}}$ as
large as we can for each $g$. For $g=2,3,\dots,11$, we use Table~\ref{R.Table}
to choose $x=7,5,31,9,20,15,30,24,33,25$, respectively (see
Table~\ref{Witness.Table} for witnesses to the values claimed for $R(g,x)$).
This yields the first group of assertions in Theorem~\ref{R.Lower.Bounds.thm}.
For $g\ge12$, we set $x=3g-\floor{g/3}+1$ and appeal to
Theorem~\ref{R.Construction.Small.gn.lem}, giving the second assertion of
Theorem~\ref{R.Lower.Bounds.thm}.

We remark that the above proof gives the more refined result
\begin{equation*}
R(2g,n) \ge \frac{11}{8\sqrt3}\sqrt{2gn} \,\Big( 1 + O\Big( g^{-1} + \Big(
\frac ng \Big)^{\!(\alpha-1)/2} \,\Big) \Big)
\end{equation*}
as $\frac ng$ and $g$ both go to infinity, where $\alpha<1$ is any number such that for
sufficiently large $y$, there is always a prime between $y-y^\alpha$ and $y$. For
instance, we can take $\alpha=0.525$ by~\cite{2001.Baker.Harman.Pintz}. This
clarification implies the final assertion of the theorem for even $g$, and the obvious
inequality $R(2g+1,n)\ge R(2g,n)$ implies the final assertion for odd $g$ as well.
\end{proof}

Habsieger and Plagne~\cite{Habsieger.Plagne} have proven that
$R(2,x)/\sqrt{4x}$ is maximized at $x=7$. For $g>2$, we have chosen $x$ based
solely on the computations reported in Table~\ref{R.Table}. For general $g$,
it appears that $R(g,x)/\sqrt{2gx}$ is actually maximized at a fairly small
value of $x$, suggesting that this construction suffers from ``edge effects''
and is not best possible.

\pagebreak\section{Upper Bounds for $\De$} \label{upper.bounds.for.De.section}

\subsection{Upper Bounds Derived from Constructions of \Bg Sets}

In Section 5 we used the connection between $\Bg$ sets and measurable sets
with small symmetric subsets to deduce upper bounds for $R(g,n)$ from lower
bounds for $\De$. In this section we exploit this relationship in the opposite
direction, converting the lower bounds on $R(g,n)$ established in Section 6 into
upper bounds for $\De$. Our first proposition verifies the statement of Theorem
\ref{Delta.Summary.thm}(i).

\begin{prop}
$\De=2\e-1$ for $\tfrac{11}{16}\leq \e \leq 1$, and $\De\ge2\e-1$ for all $0<\e\le1$. \label{Line.2e-1.prop}
\end{prop}

\begin{proof}
Recall from Lemma \ref{Lipschitz.Condition.lem} that the function $\Delta$ satisfies the Lipschitz condition $|\Delta(x)-\Delta(y)|\leq 2|x-y|$. Therefore the inequality $\De\ge2\e-1$ for all $0<\e\le1$ follows easily from the trivial value $\Delta(1)=1$.  To prove that $\De=2\e-1$ for $\tfrac{11}{16}\leq \e \leq 1$, then, it suffices to prove that $\De\le2\e-1$ in that range; and again by the Lipschitz condition, it suffices to prove simply that $\Delta\left(\frac{11}{16}\right) \le\frac38$.

For any positive integer $g$, we combine Proposition~\ref{Delta.from.R.Connection.prop} and Lemma~\ref{R.Construction.Small.gn.lem} and the monotonicity of $\Delta$ to
see that
\begin{equation*}
\frac{g}{3g-\floor{g/3}+1} \ge \Delta\bigg(
\frac{R(g,3g-\floor{g/3}+1)}{3g-\floor{g/3}+1} \bigg) \ge \Delta\bigg(
\frac{g+2\floor{g/3}+\floor{g/6}}{3g-\floor{g/3}+1} \bigg).
\end{equation*}
Since $\Delta$ is continuous by Lemma~\ref{Lipschitz.Condition.lem}, we may
take the limit of both sides as $g\to\infty$ to obtain
$\Delta\left(\frac{11}{16}\right) \le \frac38$ as desired.
\end{proof}

\begin{prop}
The function $\frac{\De}{\e^2}$ is increasing on $(0,1]$.
\label{De.over.e2.increasing.prop}
\end{prop}

\begin{proof}
Choose $0<\e<\e_0$. By Proposition~\ref{HalfModular.Construction.prop}, we have
$$
\frac{R(g,x)}{x} \frac{C(h,y)}{y} \leq \frac{R(gh,xy)}{xy}.
$$
With the monotonicity of $\De$ and Proposition~\ref{Delta.from.R.Connection.prop},
this gives
$$
\Delta \left( \frac{R(g,x)}{x} \frac{C(h,y)}{y} \right) \leq \Delta\left(
\frac{R(gh,x y)}{x y}  \right) \leq\frac{gh}{xy}.
$$
Let $g_i,x_i$ be such that $\frac{R(g_i,x_i)}{x_i}\to \e_0$ and
$\frac{g_i}{x_i}\to \Delta(\e_0)$, which is possible by
Proposition~\ref{Delta.as.infimum.prop}. By Proposition
\ref{circle.probabilistic.prop}, we may choose sequences of integers $h_j$ and
$y_j$ such that $\frac{C(h_j,y_j)}{y_j} \gtrsim \frac\e{\e_0}$ and
$\frac{h_j}{y_j} \lesssim \big( \frac\e{\e_0} \big)^2$ as $j\to\infty$. This implies
$$
\frac{R(g_i,x_i)}{x} \frac{C(h_j,y_j)}{y_j} \gtrsim \e
 \quad \text{and}\quad
\frac{g_i}{x_i} \frac{h_j}{y_j} \lesssim \Delta(\e_0) \Big( \frac{\e}{\e_0}
\Big)^2,
$$
so that, again using the monotonicity and continuity of $\Delta$,
\begin{equation*}
\Delta(\e_0) \frac{\e^2}{\e_0^2} \gtrsim \frac{g_ih_j}{x_iy_j} \ge \Delta
\left( \frac{R(g_i,x_i)}{x_i} \frac{C(h_j,y_j)}{y_j} \right) \gtrsim \De
\end{equation*}
as $j\to\infty$. This shows that $\frac{\De}{\e^2} \le \frac{\Delta(\e_0)}{\e_0^2}$ as desired.
\end{proof}

We can immediately deduce two nice consequences of this proposition.

\begin{cor}
$\lim_{\e\to0^+} \frac{\De}{\e^2}$ exists.
\end{cor}

\begin{proof}
This follows from the fact that the function $\frac{\De}{\e^2}$ is
increasing and bounded below by $\frac12$ on $(0,1]$ by the trivial lower
bound (Lemma \ref{Delta.Trivial.Lower.Bound.lem}).
\end{proof}

\begin{cor}
$\De \leq \tfrac{96}{121}\e^2$ for $0\leq \e\leq \frac{11}{16}$.
\label{part.iii.cor}
\end{cor}

\begin{proof}
This follows from the value $\Delta\big( \frac{11}{16} \big) = \frac38$
calculated in Proposition~\ref{Line.2e-1.prop} and the fact that the function
$\frac{\De}{\e^2}$ is increasing.
\end{proof}

The corollary above proves part (iv) of Theorem~\ref{Delta.Summary.thm}, leaving only part (v) yet to be established. The following proposition finishes the proof of Theorem~\ref{Delta.Summary.thm}.

\begin{prop}
$\frac\De{\e^2} \le \frac\pi{(1+\sqrt{1-\e})^2}$ for all $0<\e\le1$.
\label{part.iv.prop}
\end{prop}

\begin{proof}
Define $\alpha:=1-\sqrt{1-\e}$, so that $2\alpha-\alpha^2=\e$. If we set
$\gamma=\pi\alpha^2 n$ in the proof of Proposition \ref{probabilistic.bg.prop}, then the
sets constructed are \Bg sets with $g = \pi\alpha^2n + O(\sqrt{n\log n})$ and have size
at least
\begin{equation*}
E_0 - a_0 = 2\sqrt{\frac{\pi\alpha^2n^2}\pi} - \frac{\pi\alpha^2n}\pi + O(1 +
a_0) = (2\alpha-\alpha^2)n + O((\gamma n)^{1/4}) = \e n + O(\sqrt n)
\end{equation*}
from Eq.~(\ref{second.main.term.eq}).

Therefore, for these values of $g$ and $n$,
\begin{equation*}
\Delta \big( \frac{R(g,n)}n \big) \ge \Delta \big( \frac{\e n + O(\sqrt n)}n
\big) \to \De
\end{equation*}
as $n$ goes to infinity, by the continuity of $\Delta$. On the other hand, we see by
Proposition~\ref{Delta.from.R.Connection.prop} that
\begin{multline*}
\e^{-2}\Delta \big( \frac{R(g,n)}n \big) \le \frac {\e^{-2}g}n =
\frac{\pi\alpha^2n + O(\sqrt{n\log n})}{\e^2 n} \\
= \frac{\pi \alpha^2}{(2\alpha-\alpha^2)^2} + O\Big( \sqrt{\frac{\log n}{\e^2
n}} \,\Big) =  \frac{\pi}{(2-\alpha)^2} + o(1) \to
\frac{\pi}{(1+\sqrt{1-\e})^2}
\end{multline*}
as $n$ goes to infinity. Combining these two inequalities yields
$\frac{\De}{\e^2} \le \frac{\pi}{(1+\sqrt{1-\e})^2}$ as desired.
\end{proof}

\subsection{Upper Bounds Derived from Finite Unions of Intervals}

Another way to approach bounding $\De$ is to compute precisely $\Delta_k(\e)$,
the supremum of those real numbers $\delta$ such that every subset of $[0,1)$
with measure $\e$ that is the union of $k$ intervals has a symmetric subset
with measure $\delta$. From the definition it is easy to see that
$\Delta_k(\e)$ is a decreasing function of $k$. In fact, it follows directly
from the proof of Proposition~\ref{Delta.as.infimum.prop} that $\De =
\inf_k \Delta_k(\e) = \lim_{k\to\infty} \Delta_k(\e)$. We also have the
following formula.

\begin{lem}
    \label{DeltaK.Computation.Tool.lem}
Let $E=\bigcup_{i=1}^k (\alpha_i,\beta_i)\subseteq[0,1)$ be the union of $k$ disjoint
intervals. The largest measure of a symmetric subset of $E$ is
    $$D(E) = \max_{0\leq c\leq 1}\left\{\sum_{i=1}^k\sum_{j=1}^k \max\left\{0,
        \min\{c-\alpha_i,\beta_j-c\} -
        \max\{c-\beta_i,\alpha_j-c\}\right\}\right\}.
    $$
\end{lem}

\begin{proof}
This follows immediately from the fact that, for any real numbers $0\leq
\alpha_i<\beta_i\leq1$ and $c$, the measure of the set $\{x\colon c-x\in
(\alpha_1,\beta_1),\, c+x\in(\alpha_2,\beta_2)\}$ is $\max\big\{0,
\min\{c-\alpha_1,\beta_2-c\}$ $-$ $\max\{c-\beta_1,\alpha_2-c\}\big\}$.
\end{proof}

Since the maximum over $c$ is achieved for some $c$ that is the midpoint of
endpoints of the intervals (i.e., $c=(\alpha_i+\beta_j)/2$,
$c=(\alpha_i+\alpha_j)/2$, or $c=(\beta_i+\beta_j)/2$), this theorem provides
an effective method for the computation of $D(E)$ for any particular union of
$k$ intervals. Indeed, this formula reduces finding a particular value
$\Delta_k(\e)$ to a finite number of linear programming problems, though in
practice this computation becomes unmanageably large even for small values of
$k$. In principal, the entire function $\Delta_k$ could be calculated by
solving these linear programming problems with the constant $\e$ remaining
unspecified, branching finitely many times depending on various simple
inequalities for $\e$.

\begin{figure}
 \begin{center}
 \begin{picture}(384,240)
    \put(370,12){$\e$}
    \put(22,228){$\Delta_k(\e)$}
    \includegraphics{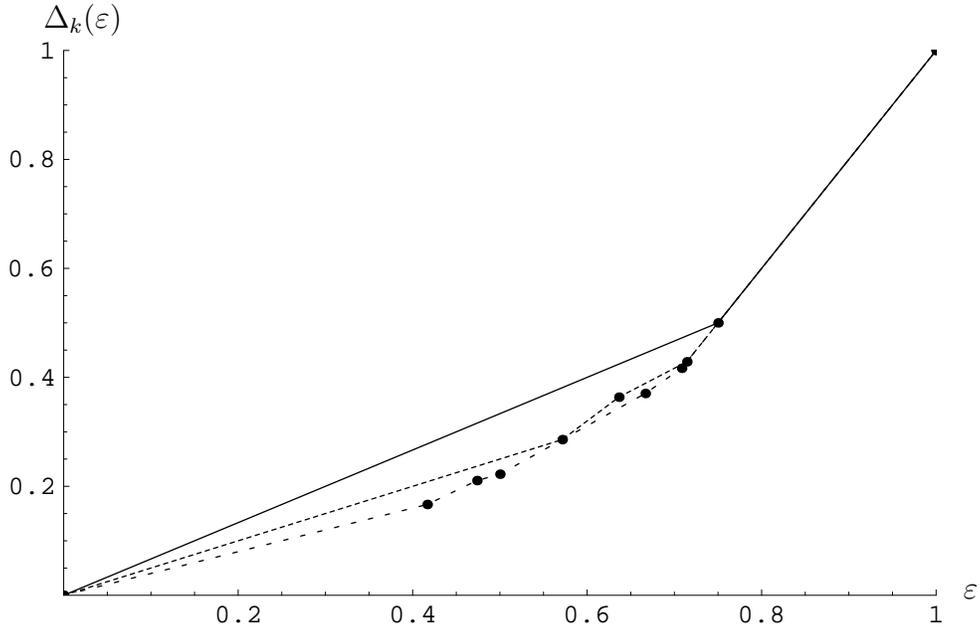}
 \end{picture}
 \end{center}
    \caption{The graph of $\Delta_2(\e)$ and the conjectured graphs of
                    $\Delta_3(\e)$ and $\Delta_4(\e)$}
\label{Delta.234.fig}
\end{figure}

The case $k=2$ is simple enough to deal with directly; we state the result in
Proposition~\ref{Delta2.prop} but omit the proof. We have shown computationally that the
graphs of $\Delta_3(\e)$ and $\Delta_4(\e)$ lie on or below the polygonal paths
described in Conjecture~\ref{Delta3.cnj} (see Figure~\ref{Delta.234.fig}), but we have
not verified that these upper bounds are in fact sharp.

\begin{prop}
    \label{Delta2.prop} The graph of the function $(\e,\Delta_2(\e))$ is
the polygonal path connecting $(0,0)$, $(3/4,1/2)$, and $(1,1)$.
\end{prop}

\begin{cnj}
\label{Delta3.cnj}
The graph of the function $\Delta_3(\e)$ is the polygonal path connecting
$(0,0)$, $(\frac47,\frac27)$, $ (\frac7{11},\frac4{11})$,
$(\frac57,\frac37)$, and $(1,1)$. The graph of the function $\Delta_4(\e)$ is
the polygonal path connecting $(0,0)$, $(\frac5{12},\frac16)$,
$(\frac9{19},\frac4{19})$, $(\frac12,\frac29)$, $(\frac23,\frac{10}{27})$,
$(\frac{17}{24},\frac5{12})$, and $(1,1)$.
\end{cnj}

It seems likely that the graph of the function $\Delta_k(\e)$ always contains
the line segment connecting $(0,0)$ and $(\frac{k+1}n,\frac2n)$, where $n$ is
the least integer for which $R(2,n)=k+1$. It is easy to show that for every
$k\ge2$, the graph of the function $\Delta_k(\e)$ contains the line segment
connecting the points $(\frac34,\frac12)$ and $(1,1)$.

\pagebreak \section{Some Remaining Questions}

We group the problems in this section into three categories, although some problems do
not fit clearly into any of the categories and others fit into more than one. (We also
refer the reader to the
Conjectures~\ref{infinity.two.cnj},~\ref{uniformly.distributed.cnj} and~\ref{Delta3.cnj}
already propounded.)

\subsection{Properties of the Function $\De$}

The first open problem on the list must of course be the exact determination
of $\De$ for all values $0\le\e\le1$. In the course of our investigations, we
have come to believe the following assertion.

\begin{cnj}
$\De = \max\{2\e-1,\frac\pi4\e^2\}$ for all $0\le\e\le1$.
\label{bold.De.cnj}
\end{cnj}

Notice that the upper bounds given in Theorem~\ref{Delta.Summary.thm} are not too far
from this conjecture, the difference between the constants $\frac{96}{121}=0.7934$ and
$\frac\pi4=0.7854$ in the middle range for $\e$ being the only discrepancy. In fact, we
believe it might be possible to prove that the expression in Conjecture
\ref{bold.De.cnj} is indeed an upper bound for $\De$ by a more refined application of
the probabilistic method employed in Section~\ref{probabilistic.bg.section}. The key
would be to show that the various events $R_k > \gamma + 1 + a$ are more or less
independent of one another (as it stands we have to assume the worst---that they are all
mutually exclusive---in obtaining our bound for the probability of obtaining a ``bad''
set), so that we could omit the $\log3n$ factor from the chosen value of $a$.

There are some intermediate qualitative results about the function $\De$ that
might be easier to resolve. It seems likely that $\De$ is convex, for example,
but we have not been able to prove this. A first step towards clarifying the
nature of $\De$ might be to prove that
$$\frac{|\Delta(x)-\Delta(y)|}{|x-y|} \ll \max\{x,y\}.$$
Also, we would not be surprised to see accomplished an exact computation of
$\Delta(\tfrac12)$, but we have been unable to make this computation ourselves.

We do not believe that there is always a set with measure $\e$ whose largest symmetric
subset has measure $\De$. In fact, we do not believe that there is a set with measure
$\e_0:=\inf\{\e\colon \De = 2\e-1\}$ whose largest symmetric subset has measure
$\Delta(\e_0)$, but we do not even know the value of $\e_0$. In
Proposition~\ref{Line.2e-1.prop}, we showed that $\e_0 \leq \tfrac{11}{16}$, but this
was found by rather limited computations and is unlikely to be sharp. The quantity
$\frac{11}{8\sqrt{3}}$ in Theorem~\ref{R.Lower.Bounds.thm} is of the form
$\frac{\e_0}{\sqrt{4\e_0-2}}$, and similarly the quantity $\frac{96}{121}$ in Theorem
\ref{Delta.Summary.thm}(iv) is of the form $\frac{2\e_0-1}{\e_0^2}$. Thus any
improvement in the bound $\e_0\leq \frac{11}{16}$ would immediately result in
improvements to Theorem~\ref{R.Lower.Bounds.thm} and
Theorem~\ref{Delta.Summary.thm}(iv). We remark that Conjecture \ref{bold.De.cnj} implies
that $\e_0 = \frac2{2+\sqrt{4-\pi}} = 0.6834$, which in turn would allow us to replace
the constant $\frac{11}{8\sqrt{3}}$ in Theorem~\ref{R.Lower.Bounds.thm} by
$\sqrt{\frac2\pi}$.

\subsection{Artifacts of our Proof}\label{artifacts.section}

Let $\K$ be the class of functions $K\in L^2(\T)$ satisfying $K(x)\ge1$ on
$[-\frac14,\frac14]$. How small can we make $\|\hat{K}\|_p$ for $1\le p\le2$?
We are especially interested in $p=\frac43$, but a solution for any $p$ may be
enlightening.

To give some perspective to this problem, note that a trivial upper bound for
$\inf_{K\in\K} \|\hat{K}\|_p$ can be found by taking $K$ to be identically equal to 1,
which yields $\|\hat K\|_p=1$. One can find functions that improve upon this trivial
choice; for example, the function $K$ defined in Eq.~(\ref{easy.K.def}) is an example
where $\|\hat{K}\|_{4/3} = 0.96585$. On the other hand, since the $\ell^p$-norm of a
sequence is a decreasing function of $p$, Parseval's identity immediately gives us the
lower bound $\|\hat K\|_p \ge \|\hat K\|_2 = \|K\|_2 \ge \big( \int_{-1/4}^{1/4} 1^2\,dt
\big)^{1/2} = \frac1{\sqrt2} = 0.707107$, and of course this is the exact minimum for
$p=2$.

We remark that Proposition~\ref{Cf.First.Result.prop} and the function $b(x)$ defined
after the proof of Corollary~\ref{First.cor} provide a stronger lower bound for $1\le p\le\frac43$. By direct
computation we have $1.14939> \|b\ast b\|_2^2$, and by
Proposition~\ref{Cf.First.Result.prop} we have $\|b\ast b\|_2^2 \geq
\|\hat{K}\|_{4/3}^{-4}$ for any $K\in\K$. Together these imply that
$\|\hat{K}\|_{p}\geq \|\hat{K}\|_{4/3}> 0.96579.$  In particular, for $p=\frac43$ we know the value of $\inf_{K\in\K} \|\hat K\|_{4/3}$ to within one part in ten thousand.
The problem of determining the actual
infimum for $1<p<2$ seems quite mysterious. We remark that Green~\cite{2001.Green}
considered the discrete version of a similar optimization problem, namely the
minimization of $\|\hat K\|_p$ over all pdfs $K$ supported on $[ {-\frac14},\frac14 ]$.

As mentioned at the end of Section~\ref{basic.argument.section}, we used the inequality
$\| g \|_2^2 \leq \|g\|_\infty \|g\|_1$ which is exact when $g$ takes on one non-zero
value, i.e., when $g$ is an nif. We apply this
inequality when $g=\ff$ with $f$ supported on an interval of length $\frac12$, which usually looks very different from an nif. In this circumstance, the inequality does not seem to be best possible, although the corresponding inequality in the exponential sums approach of~\cite{Cilleruelo.Ruzsa.Trujillo} and in the discrete Fourier approach of \cite{2001.Green} clearly is best possible. Specifically, we ask for a lower bound on
 $$ \sup_{\substack{\text{$f \ge 0$} \\
        \supp(f)\subseteq[-\frac14,\frac14]}}
        \frac{\ffi \|\ff\|_1}{\| \ff \|_2^2}$$
that is strictly greater than 1.

\subsection{The Analogous Problem for Other Sets}

More generally, for any subset $E$ of an abelian group endowed with a measure, we can
define $\Delta_E(\e):=\inf\{D(A) \colon A \subseteq E,\, \lambda(A)=\e\}$, where $D(A)$
is defined in the same way as in Eq.~(\ref{Ddef}). For example, $\Delta_{[0,1]}(\e)$ is
the function $\De$ we have been considering throughout this paper, and $\Delta_\T(\e)$
was considered in Section \ref{circular.probabilistic.section}.

We believe that for each $E\subseteq\R$, there is a positive constant $c$ such that $\De
\geq \Delta_E(\e)\geq c \e^2$ for $0<\e\leq 1$. The work of Abbott~\cite{1990.Abbott}
seems relevant. If we are concerned only with $\e\to 0$, and we normalize by considering
only sets $E$ for which $\lambda(E)=1$, then it may be possible to take an absolute
constant. In other words, is it true that
$$ \inf_{\substack{E\subseteq\R \\ \lambda(E)=1}} \liminf_{\e\to0}
        \frac{\Delta_E(\e)}{\e^2} > 0\,?$$

Most of the work in this paper generalizes easily from $E=[0,1]$ to
$E=[0,1]^d$. We have had difficulties, however, in finding good kernel
functions in higher dimensions. That is, we need functions $K(\bar{x})$ such
that
$$\sum_{\bar{j}\in\Z^d} \big| \hat{K} \big( \bar{j} \big) \big|^{4/3}$$
is as small as possible, while $K(\bar{x})\geq 1$ if all components of
$\bar{x}$ are less than $\frac14$ in absolute value. This restricts $K$ on
one-half of the space in 1 dimension, one-quarter of the space in 2 dimensions,
and only $2^{-d}$ of the space in $d$ dimensions. For this reason one might
expect that {\em better} kernels exist in higher dimensions, but the
computational difficulties have prevented us from finding them.

\bigskip{\small
{\it Acknowledgements.} The authors thank Heini Halberstam for thoughtful
readings of this man\-u\-script and helpful suggestions. The first author was
supported in part by grants from the Natural Sciences and Engineering Research Council.}

\end{document}